\documentstyle[german,12pt,leqno]{article}
\newtheorem{lemma}{Lemma}[section]
\newtheorem{satz}[lemma]{Satz}
\newtheorem{cor}[lemma]{Korollar}
\newtheorem{beh}[lemma]{Behauptung}
\begin{document}

 \begin{titlepage}
  \vspace*{2cm}
   \Large\bf
   \begin{center}
     (FP)$_{\bf n}$-Eigenschaften der verallgemeinerten\\
            Houghton-Gruppen
   \end{center}
  \vfill
   \large
   \begin{center}
      Diplomarbeit\\
         von\\
      Heike Sach\\
    \end{center}
  \vfill
  \begin{center}
    Fachbereich Mathematik\\
   Johann Wolfgang Goethe-Universit"at Frankfurt/Main\\
    Januar 1992
   \end{center}
 \end{titlepage}
   
\section*{Einleitung}

Eine Gruppe $G$ hei\3t {\em vom Typ\/} {\sc (FP)}$_n$ (bzw. {\sc (FP)}
$_{\infty}$), wenn ${\sl Z}\!\!{\sl Z}$ als ${\sl Z}\!\!{\sl Z}G$-Modul
 eine projektive Aufl"osung besitzt, die in jeder Dimension$\,\le n$ (bzw.
 in allen Dimensionen) endlich erzeugt ist (vgl. [Bi] oder [Br2]). Zum Beispiel ist eine Gruppe
 $G$ genau dann vom Typ {\sc (FP)}$_1$, wenn sie endlich erzeugt ist, und
 $G$ ist vom Typ {\sc (FP)}$_2$, falls sie endlich pr"asentiert ist. (Die
 Umkehrung der letzten Behauptung konnte bisher nicht gezeigt werden). 
 Weiterhin ist jede endliche Gruppe vom Typ {\sc (FP)}$_{\infty}$.\\

Ausgangspunkt dieser Diplomarbeit bildet eine Folge von Gruppen $H_n$
 ($n\ge 1$), die von C. Houghton eingef"uhrt wurde (vgl. [Hou]). Betrachtet
 man die Menge $S={\rm I}\!{\rm N}\times \{1,\ldots,n\}$, so ist $H_n$ die
 Gruppe aller Permutationen $h$ von $S$, f"ur die ein $x_0\in{\rm I}\!{\rm N}$
 sowie ein $n$-Tupel $(m_1,\ldots,m_n)\in{\sl Z}\!\!{\sl Z}^n$ existiert,
 soda\3 $(x,i)h=(x+m_i,i)$ f"ur alle $x\ge x_0$ erf"ullt ist. Stellt man sich 
 die Menge ${\rm I}\!{\rm N}\times\{i\}$ f"ur $i=1,\ldots,n$ jeweils auf 
 einem Zahlenstrahl angeordnet vor, so bedeutet diese Forderung, da\3 die
 Elemente des $i$-ten Zahlenstrahles ab einer Stelle $x_0$ um $m_i$
 verschoben werden. Die Permutationen aus $H_n$ verhalten sich somit f"ur
 hinreichend gro\3e $x$ wie eine Translation.\\

Die Gruppe $H_n$ besitzt besondere Endlichkeitseigenschaften: Wie von K.S.
 Brown in [Br1] gezeigt wurde, ist $H_n$ vom Typ {\sc (FP)}$_{n-1}$, jedoch
 nicht vom Typ {\sc (FP)}$_n$ und f"ur $n\ge 3$ endlich pr"asentiert.
Beispiele f"ur Gruppen mit solchen Eigenschaften liegen nicht auf der Hand,
und wir wollen in dieser Arbeit ein weiteres vorstellen, das durch eine 
 Verallgemeinerung von Houghton's Gruppen entstanden ist.\\

 Den Ansatzpunkt der Verallgemeinerung bildet die der Gruppe $H_n$
 zugrundeliegende Menge ${\rm I}\!{\rm N}\times\{1,\ldots,n\}$. F"ur
 $i=1,\ldots,n$ ersetzen wir ${\rm I}\!{\rm N}\times\{i\}$
  durch $({\rm I}\!{\rm N}
 \times{\rm I}\!{\rm N})\times\{i\}$. Die Elemente aus $({\rm I}\!{\rm N}
 \times{\rm I}\!{\rm N})\times\{i\}$ kann man sich als ganzzahlige
  Gitterpunkte eines Quadranten innerhalb eines kartesischen Koordinatensystemes
 vorstellen, und die neue Menge $S=({\rm I}\!{\rm N}\times{\rm I}\!{\rm N})
 \times\{1,\ldots,n\}$ befindet sich gem"a\3 dieser Auffassung in der 
 disjunkten Vereinigung von $n$ Quadranten.

 Betrachtet werden nun spezielle Permutationen dieser Menge $S=({\rm I}\!{\rm N}\times{\rm I}\!{\rm N})
 \times\{1,\ldots,n\}$, deren Hauptmerkmal -- in Analogie zu den Houghton-Gruppen --
 darin besteht, ab einer Stelle $p_0=(x_0,y_0)$ (d.h. f"ur alle $((x,y),i)$
 aus $S$ mit $x\ge x_0$ und $y\ge y_0$) in jedem der $n$ Quadranten eine 
 Translation zu sein.

 Im Unterschied zur Menge ${\rm I}\!{\rm N}\times\{i\}$ 
 gibt es jedoch beliebig viele Translationsrichtungen, die das "`Gitter"'
$({\rm I}\!{\rm N}\times{\rm I}\!{\rm N})\times\{i\}$ invariant lassen. Aus
 diesem Grunde stellen wir zwei verschieden Typen von Gruppen als m"ogliche
 Verallgemeinerung der Houghton-Gruppen vor:

 Zum einen die Gruppen $G_n$, bei denen die ab einer Stelle $p_0$ beginnende
 Translation in jedem der Quadranten nur in {\em einer\/} Richtung, n"amlich in 
 Richtung der Winkelhalbierenden, stattfindet, und zum anderen die 
 Gruppen $\widetilde{G}_n$, bei denen hinsichtlich dieser Translationen jede
 (ganzzahlige) Richtung zugelassen wird.\\

Es stellt sich die Frage, welche der beiden Gruppentypen die "`nat"urlichere"'
 Verallgemeinerung der Gruppen $H_n$ ist. F"ur $G_n$ spricht, da\3 bei den 
 Houghton-Gruppen ebenfalls nur eine m"ogliche Translationsrichtung (entlang
 des Zahlenstrahles) vorhanden ist. Andererseits liegt es nahe, bei 
 Erweiterung der Menge ${\rm I}\!{\rm N}\times\{i\}$ auf $({\rm I}\!{\rm N}
 \times{\rm I}\!{\rm N})\times\{i\}$ alle dadurch zus"atzlich entstehende
 Translationsrichtungen auszusch"opfen.\\

Wie dem auch sei ist es als Hauptresultat dieser Arbeit gelungen, mit der
 von K.S. Brown im Falle der Houghton-Gruppen angewandten Methode nachzuweisen,
 da\3 $G_n$ vom Typ {\sc (FP)}$_{n-1}$, nicht vom Typ {\sc (FP)}$_n$ und
 f"ur $n\ge 3$ endlich pr"asentiert ist (vgl. Satz 4.9).

Dagegen ist der Versuch einer "Ubertragung des Beweises auf die Gruppen
 $\widetilde{G}_n$ gescheitert, was seine Ursache im wesentlichen in der 
 Notwendigkeit von zwei (anstelle von nur einer) Erzeugenden f"ur die
 Translationen von $({\rm I}\!{\rm N}\times{\rm I}\!{\rm N})$
 hatte. Allerdings kann man zeigen, da\3 es eine exakte Sequenz von Gruppen
    \[ 0\longrightarrow G_n\longrightarrow 
         \widetilde{G}_n\longrightarrow
           {\sl Z}\!\!{\sl Z}^{n-1}\longrightarrow 0  \]
  gibt, woraus sich schlie\3en
l"a\3t, da\3 $\widetilde{G}_n$ zumindest vom Typ {\sc (FP)}$_{n-1}$ und f"ur
 $n\ge 3$ endlich pr"asentiert ist. Ob $\widetilde{G}_n$ ebenfalls die 
 Eigenschaft besitzt, nicht vom Typ {\sc (FP)}$_n$ zu sein, bleibt offen
 und kann vielleicht mit anderen Methoden gekl"art werden.\\\\

Es folgt ein "Uberblick:

Im ersten Abschnitt stellen wir ein Kriterium von K.S. Brown f"ur die 
 {\sc (FP)}$_n$- Eigenschaft sowie f"ur die endliche Pr"asentierbarkeit
 einer Gruppe $G$ vor, und geben eine auf unsere Anwendungen zugeschnittene
 kombinierte Version beider Kriterien an:

 \begin{quote}
    Zu festem $n\ge 1$ sei ein zusammenziehbarer $CW$-Komplex $X$ gegeben, auf dem
   die Gruppe $G$ durch Hom"oomorphismen operiert, welche die Zellen von
   $X$ permutieren, 
   soda\3 der Stabilisator jeder $p$-Zelle $\sigma$ aus $X$ vom 
   Typ {\sc (FP)}$_n$ und im Falle $n\ge 3$ endlich pr"asentiert ist. Weiterhin sei 
   $\{X_j\}_{j\in {\rm I}\!{\rm N}}$ eine Filtrierung von $X$, soda\3 jedes $X_j$ endlich ist
    mod $G$.\\ Entsteht f"ur alle hinreichend gro\3en $j$ $X_{j+1}$ aus $X_j$ bis auf Homotopie 
   durch Hinzuf"ugen von $n$-Zellen, so ist $G$ vom Typ {\sc (FP)}$_{n-1}$ und nicht vom Typ {\sc (FP)}$_n$.
   Im Falle $n\ge 3$ ist $G$ endlich pr"asentiert.
 \end{quote}

 Als Grundlage f"ur die weitere Arbeit wird die Anwendung dieses Kriteriums
 im Falle der Houghton-Gruppen erl"autert.\\

Der zweite Abschnitt enth"alt eine formale Beschreibung der Gruppen $G_n$ 
 und $\widetilde{G}_n$, sowie den Nachweis, da\3 die Houghtongruppe $H_n$
 sowohl als Untergruppe wie auch als Faktorgruppe von $G_n$ und $\widetilde{G}
 _n$ vorkommt. Dar"uberhinaus zeigen wir, da\3 $G_n$ ein Normalteiler von
 $\widetilde{G}_n$ ist, und da\3 die Faktorgruppe von $\widetilde{G}_n$ nach
 $G_n$ isomorph zur Gruppe ${\sl Z}\!\!{\sl Z}^{n-1}$ ist.\\

Im dritten Abschnitt beginnen wir mit der Herstellung der Voraussetzungen
 des obigen Kriteriums f"ur die Gruppen $G_n$. Es wird ein zusammenziehbarer
 Simplizialkomplex konstruiert, auf dem $G_n$ simplizial operiert, wobei
 die Stabilisatoren dieser Operation jeweils eine Untergruppe von endlichem
  Index enthalten, die isomorph zu einer der Gruppen
 $H_m$ f"ur ein $m\ge n+1$ ist. (Somit sind die Stabilisatoren vom Typ {\sc (FP)}$_n$ und
 im Falle $n\ge 3$ zus"atzlich endlich pr"asentiert, vgl. Lemma 1.4).\\

Im letzten Abschnitt geben wir eine Filtrierung dieses Simplizialkomplexes
 an, und weisen die im Kriterium hinsichtlich der Filtrierung geforderten
 Eigenschafen nach.\newpage

\tableofcontents\newpage

\pagenumbering{arabic}
\section{Vorbereitungen}

\subsection{Endlichkeitskriterien}

Wir stellen in diesem Abschnitt zwei Endlichkeitskriterien vor, zum einen f"ur die
({\sc FP})$_n$-Eigenschaft und zum anderen f"ur die endliche Pr"asentierbarkeit einer Gruppe $G$.
Dabei setzen wir stets $n\ge 1$ voraus.\\
In beiden F"allen handelt es sich um notwendige und hinreichende Kriterien, die gewisse topologische 
 Eigenschaften eines $CW$-Komplexes $X$, auf dem die Gruppe in geeigneter Weise operiert, miteinbeziehen. Sie wurden von 
 K.S. Brown in seiner Arbeit "`Finetness properties of groups"' ([Br1]) entwickelt, und 
 wir geben daraus eine kurze Zusammenfassung der Abschnitte 1-3. 
 Es werden im wesentlichen nur die zur Formulierung der Kriterien n"otigen
 Begriffe erl"autert, ohne auf Beweise n"aher einzugehen. F"ur eine ausf"uhrliche
 Darstellung der Resultate sei auf die oben erw"ahnte Arbeit verwiesen.\\\\
{\bf Def\/inition:}$\:$Sei $G$ eine Gruppe. Dann versteht man unter einem {\em G-Komplex\/}
 einen $CW$-Komplex $X$, auf dem $G$ durch Hom"oomorphismen operiert, welche die
 Zellen von $X$ permutieren.\\\\
 Zum Beispiel ist jeder Simplizialkomplex $X$, auf dem $G$ simplizial operiert, ein $G$-Komplex.
 Wir ben"otigen $G$-Komplexe mit speziellen Eigenschaften:\\\\
{\bf Def\/inition:}$\:$ Ein $G$-Komplex $X$ hei\3t {\em$n$-gut\/} f"ur die Gruppe $G$, wenn die beiden
 folgenden Bedingungen erf"ullt sind:
   \begin{itemize}
     \item[(i)]$\;$ Die reduzierten Homologiegruppen $\widetilde{H}_i(X)=0\quad$f"ur alle $i<n$
     \item[(ii)]$\;$ F"ur alle $0\le p\le n$ ist der Stabilisator Stab$_G$($\sigma$) jeder $p$-Zelle
      $\sigma$ aus $X$ vom \linebreak \hspace*{1.5 ex} Typ ({\sc FP})$_{n-p}\:$.
   \end{itemize}
 Gibt es  einen $n$-guten $G$-Komplex $X$ mit 
  endlichem  $n$-dimensionalen Ger"ust  
  mod $G$, so folgt die ({\sc FP})$_n$-Eigenschaft von $G$ (vgl. Proposition 1.1, [Br]).\\ Ein "ahnliches Ergebnis
 liegt f"ur die endliche Pr"asentierbarkeit vor: Bei speziellen $G$-Komplexen $X$ kann man schlie\3en, 
 da\3 $G$ endlich pr"asentiert ist, falls $X$ ein endliches $2$-Ger"ust mod $G$ besitzt
 (s. Proposition 3.1, [Br]). \\
 Diese Kriterien versagen im Falle eines unendlichen $n$- bzw. $2$-Ger"ustes mod $G$.
 Es ist aber m"oglich, einen $G$-Komplex $X$ in Unterkomplexe zu zerlegen,  
 die alle ein endliches $n$-Ger"ust mod $G$ haben. 
 In der folgenden Definition wird beschrieben, von welcher Form eine solche
 Zerlegung sein mu\3, um daraus Informationen
 "uber die Endlichkeitseigenschaften von $G$ gewinnen zu k"onnen.\\

{\bf Def\/inition:}$\;$ Eine {\em Filtrierung\/} eines $G$-Komplexes $X$ besteht aus einer Familie $G$-invarianter
 Unterkomplexe $\{X_{\alpha}\}_{\alpha\in D}$ mit einer gerichteten Indexmenge $D$
 (d.h. $\forall$ $\alpha,\beta\in D\quad\exists$ $\gamma\in D\quad$mit $\alpha\le\gamma$ und $\beta\le \gamma$), soda\3 f"ur 
 $\alpha\le \beta$\quad$X_{\alpha}\subseteq X_{\beta}$ und $X=\bigcup_{\alpha\in D}X_{\alpha}$ ist.
 Die Filtrierung hei\3t von {\em endlichem n-Typ\/}, falls das $n$-Ger"ust
 jedes Unterkomplexes $X_{\alpha}$ endlich ist mod $G$.
 \vfill
   {\bf Bemerkung:}$\:$ Es existiert immer eine Filtrierung von endlichem $n$-Typ, zum Beispiel die Familie aller Unterkomplexe von
 $X$ mit endlichem $n$-Ger"ust mod $G$, geordnet durch Inklusion. Dar"uberhinaus
 kann erreicht werden, da\3 $\bigcap_{\alpha\in D}X_{\alpha}\not=\emptyset$ ist.
 Andernfalls ersetzt man f"ur ein festes $v\in  X$ die Indexmenge $D$ durch $D':\,=\{\alpha\in D\mid v\in X_{\alpha}\}$ (vgl. [Br], Abschnitt 3).
 Wir k"onnen also o.B.d.A. annehmen, da\3 es in allen $X_{\alpha}$ einen gemeinsamen Basispunkt $v$ gibt.\\\\
 Aus einer solchen Filtrierung von $X$ ergeben sich direkte Systeme von Gruppen, zum einen f"ur jedes $i\ge$$0$
 das System der reduzierten Homologiegruppen $\{\widetilde{H}_i(X_{\alpha})\}_{\alpha\in D}$ und zum 
 anderen das System der Fundamentalgruppen $\{\pi_1(X_{\alpha},v)\}_{\alpha\in D}$ bzgl. gemeinsamem Basispunkt $v$.
 Die ({\sc FP})$_n$-Eigenschaft bzw. die endliche Pr"asentierbarkeit einer Gruppe $G$ sind
 unter bestimmten Voraussetzungen "aquivalent
 zu gewissen Eigenschaften dieser Systeme von Gruppen.
 \vfill
    {\bf Definition:}$\;$ Ein direktes System von Gruppen $\{A_{\alpha}\}_{\alpha\in D}$ hei\3t 
 {\em im wesentlichen trivial\/}, wenn es zu jedem $\alpha$ aus $D$ ein $\beta\ge \alpha$ gibt, 
 soda\3 die Abbildung $\,A_{\alpha}\longrightarrow A_{\beta}\,$ trivial ist.
 \vfill
Mit diesem Begriff lassen sich nun die Endlichkeitskriterien formulieren:
 \vfill
\begin{satz} Ist $X$ ein $n$-guter $G$-Komplex mit einer Filtrierung von endlichem
  $n$-Typ, so ist $G$ genau dann vom Typ {\sc (FP)}$_n$, wenn f"ur jedes $i<n$ das direkte System
  $\{\widetilde{H}_i(X_{\alpha})\}$ der reduzierten Homologiegruppen im wesentlichen trivial ist.
\end{satz}
 \vfill
{\sc Beweis}:$\:$ [Br1], Theorem 2.2 \hfill $\Box$
 \vfill
\begin{satz} Sei $X$ ein einfach zusammenh"angender $G$-Komplex, soda\3 die Stabilisatoren der Ecken endlich pr"asentiert
  und die der Kanten endlich erzeugt sind. Sei weiterhin $\{X_{\alpha}\}$ eine Filtrierung
  von endlichem $2$-Typ von $X$ mit gemeinsamem Basispunkt $v\in\bigcap X_{\alpha}$. Ist $G$ endlich erzeugt, so ist $G$ genau dann
  endlich pr"asentiert, wenn das direkte System $\{\pi_1(X_{\alpha},v)\}$ im wesentlichen
  trivial ist.
 \end{satz}
 \vfill
{\sc Beweis}:$\:$ [Br1], Theorem 3.2 \hfill $\Box$
 \vfill
    Wir benutzen in dieser Arbeit die Aussagen von Satz 1.1 und Satz 1.2 in einer kombinierten 
 Form, was unter gewissen Voraussetzungen m"oglich ist. (vgl. auch [Br1], Corollar 3.3, jedoch sind dort die
 Voraussetzungen sch"arfer. Der Beweis l"a\3t sich dennoch unver"andert "ubernehmen).
 \pagebreak

\begin{cor} Zu festem $n\ge 1$ sei ein zusammenziehbarer $G$-Komplex $X$ gegeben 
   mit der Eigenschaft, da\3 der Stabilisator jeder $p$-Zelle $\sigma$ aus $X$ vom 
   Typ {\sc (FP)}$_n$ und im Falle $n\ge 3$ endlich pr"asentiert ist. Weiterhin sei 
   $\{X_j\}_{j\in {\rm I}\!{\rm N}}$ eine Filtrierung von $X$, soda\3 jedes $X_j$ endlich ist
   {\em mod} $G$.\\ Entsteht dann f"ur alle hinreichend gro\3e $j$ $X_{j+1}$ aus $X_j$ bis auf Homotopie 
   durch Hinzuf"ugen von $n$-Zellen, so ist $G$ vom Typ {\sc (FP)}$_{n-1}$ und nicht vom Typ {\sc (FP)}$_n$.
   Im Falle $n\ge 3$ ist $G$ endlich pr"asentiert.
 \end{cor}
{\sc Beweis}:$\;$ $G$ ist vom Typ ({\sc FP})$_{n-1}$:\\
   F"ur $n=1$ ist die Behauptung trivial. Es sei also $n\ge 2$. Aus der Zusammenziehbarkeit von $X$
   ergibt sich 
     \begin{equation}  
       \begin{array}[t]{c} \lim\limits_{\to}\\ \stackrel{j\in{\rm I}\!{\rm N}}{}\end{array}
       \widetilde{H}_i(X_j)=0\quad\forall\; i\ge 0\:.   
     \end{equation}        
 Da sich durch Hinzuf"ugen von $n$-Zellen die $i$-dimensionalen 
 Homologiegruppen f"ur alle $i\le n-2$ 
 nicht ver"andern, erh"alt man f"ur alle hinreichend gro\3e $j$ und $i=0,\ldots, n-2$
    \[ \widetilde{H}_i(X_j)\cong \widetilde{H}_i(X_{j+1})\:. \]
Mit (1.1) folgt daraus $\widetilde{H}_i(X_j)=0$ f"ur alle hinreichend
  gro\3e $j$, d.h.
 das direkte System der reduzierten Homologiegruppen $\{\widetilde{H}_i(X_j)\}_{j\in{\rm I}\!{\rm N}}$ 
 ist f"ur alle $0\le i<n-1$ im wesentlichen trivial. Da $X$ und die Filtrierung $\{X_j\}_{j\in{\rm I}\!{\rm N}}$ von $X$
 offensichtlich die Voraussetzungen von Satz 1.1 erf"ullen, folgt damit die Behauptung.\\\\
  $G$ ist endlich pr"asentiert:\\
 Die endliche Pr"asentierbarkeit von $G$ im Falle $n\ge 3$ ergibt sich "ahnlich, indem man bei der Argumentation
 anstelle der reduzierten Homologiegruppen die Fundamentalgruppe $\pi_1(X_j)$ benutzt und mit
 Satz 1.2 schlie\3t.\\\\
 $G$ ist nicht vom Typ ({\sc FP})$_n$:\\
Da f"ur gro\3e $j$ $X_{j+1}$ bis auf Homotopie aus $X_j$ durch Hinzuf"ugen von $n$-Zellen entsteht, 
 ist die Abbildung $\widetilde{H}_{n-1}(X_j)\longrightarrow \widetilde{H}_{n-1}(X_{j+1})$ surjektiv.
 Der Kern dieser Abbildung ist nicht trivial, denn angenommen durch das Hinzuf"ugen von $n$-Zellen 
 verschwinden keine ($n-1$)-dimensionalen Homologieklassen.
 Dann entstehen neue $n$-dimensionale Homologieklassen,
 die durch weiteres Hinzuf"ugen von $n$-Zellen nicht mehr "`ausgel"oscht"'  werden
 k"onnen. Dies widerspricht der Zusammenziehbarkeit von $X$.
 Somit ist f"ur alle gro\3en $j$ Ker($\widetilde{H}_{n-1}(X_j)\longrightarrow \widetilde{H}_{n-1}(X_{j+1})$)$\not=0$, also insbesondere
 $\widetilde{H}_{n-1}(X_j)\not=0$. Wegen der Surjektivit"at obiger Abbildung
 ergibt sich daraus,da\3 das System $\{\widetilde{H}_{n-1}(X_j)\}$ 
 nicht im wesentlichen trivial ist, d.h.
 $G$ ist nicht vom Typ ({\sc FP})$_n$.\hfill $\Box$\\\\
Zum Schlu\3 dieses Abschnittes weisen wir auf einige n"utzliche Zusammenh"ange
 zwischen Endlichkeitseigenschaften von Gruppen und Untergruppen hin.
\begin{lemma}
    Sei $G$ eine Gruppe, $U\le G$ eine Untergruppe von endlichem Index in $G$.
    Dann $G$ ist genau dann vom Typ {\sc (FP)}$_n$, wenn $U$ es ist. 
\end{lemma}
 {\sc Beweis}:$\;$ [Bi], Proposition 2.5\hfill $\Box$\\
 \begin{lemma}
      Sei $0\longrightarrow N\longrightarrow G\longrightarrow Q\longrightarrow 0$ eine exakte Sequenz von Gruppen.
      Sind $N$ und $Q$ vom Typ {\sc (FP)}$_n$ (bzw. endlich pr"asentiert), so ist 
       auch $G$ vom Typ {\sc (FP)}$_n$ (bzw. endlich pr"asentiert).
 \end{lemma}
 {\sc Beweis}:$\;$Die  
 Aussage "uber die ({\sc FP})$_n$-Eigenschaft l"a\3t sich analog zu
  Proposition 2.7 in [Bi] zeigen. Zur endlichen Pr"asentiertheit vgl. z.Bsp.
  [Ro], Abschnitt 2.2$\,$. \hfill $\Box$\\
              
\subsection{Topologische Eigenschaften partiell geordneter Mengen}

Einer partiell geordneten Menge $M$ l"a\3t sich in folgender Weise ein Simplizialkomplex
 $\mid M\mid$ zuordnen: Jede Kette der L"ange $k, k\ge 0$ in $M$ 
    \[ a_0<a_1<\ldots<a_k \]
 bildet ein $k$-Simplex $(a_0,\ldots,a_p)$ von $\mid M\mid$.
 Die Eckenmenge von $\mid M\mid$ besteht somit aus
 $M$ selbst.\\
 Unter einem Morphismus von partiell geordneten Mengen $f:M\longrightarrow N$ 
 versteht man eine Abbildung, die ordnungserhaltend ist, d.h. f"ur alle $a,b$ aus $M$ mit 
 $a\le b$ mu\3 auch $f(a)\le f(b)$ sein. Ein solches $f$ induziert eine simpliziale Abbildung
       \[ \mid f\mid\,:\:\mid M\mid\,\longrightarrow\:\mid N\mid\:,\]
 die durch $f(a_0,\ldots,a_k):\,=(f(a_0),\ldots,f(a_k))$ gegeben ist.\\\\
{\bf Bemerkung 1.1:}$\:$ Ist $G$ eine Gruppe, die auf einer partiell geordneten Menge $M$ durch 
 ordnungserhaltende Abbildungen operiert, so induziert dies eine simpliziale Operation von $G$ auf
 $\left| M\right|$. $\left| M\right|$ ist also ein $G$-Komplex (vgl. Abschnitt 1.1).\\\\
 Es bestehen  gewisse  Zusammenh"ange zwischen den
 Eigenschaften der gegebenen Ordnung
 auf $M$ und den topologischen Eigenschaften von $\left|  M\right|$. Bez"uglich der Zusammenziehbarkeit von $\mid M\mid$
 gilt (s. [Br1], Abschnitt 4):\\
\begin{lemma} Sei $M$ eine partiell geordnete Menge mit der Eigenschaft, da\3 f"ur alle $a,b$ aus $M$ ein $c$ existiert mit
 $a\le c$ und $b\le c$ ($M$ gerichtet). Dann ist $\left| M\right|$ zusammenziehbar.
\end{lemma}
 Eine weitere Aussage l"a\3t sich hinsichtlich des Nerves von "Uberdeckungen
 machen. Wir benutzen dazu folgendes allgemeine Resultat:

\begin{lemma}
   Sei $K$ ein Simplizialkomplex, der durch eine Familie $\{K_i\}_{i\in I}$ von Unterkomplexen "uberdeckt wird.
   Ist jeder endliche Durchschnitt $K_{i_1}\cap K_{i_2}\cap\ldots\cap K_{i_n}$ 
   entweder leer oder zusammenziehbar, so besitzt $K$ denselben Homotopietyp wie
   der Nerv $\cal N$ dieser "Uberdeckung.
 \end{lemma}
 {\sc Beweis}:$\;$ [Bj], Lemma 1.1 \hfill $\Box$\\

 Definieren wir nun f"ur ein $b$ aus $M$ die Menge
    \[ M_{\le b}:\,=\{\,a\in M\mid a\le b\,\} \]
 und betrachten zu einer nichtleeren Teilmenge $J$ von $M$ 
 den Simplizialkomplex \mbox{$S_J:\,=\bigcup_{b\in J}\mid M_{\le b}\mid$,}
 so ergibt sich aus Lemma 1.7\\
\begin{cor}
   Sei $M$ eine partiell geordnete Menge, $J\subseteq M$ eine Teilmenge von $M$.
   Gibt es zu jeder endlichen Teilmenge $\{b_1,\ldots,b_n\}$ von $J$, die in $M$ nach unten beschr"ankt ist ,
    eine gr"o\3te untere Schranke $d$ in $M$, so ist $S_J=\bigcup_{b\in J}\mid M_{\le b}\mid$ homotopie"aquivalent zu
  dem Nerv $\cal N$ der Familie $\{\mid M_{\le b}\mid\}_{b\in J}$.\\
 \end{cor}
{\sc Beweis}:$\;$ Ist zu $b_1,\ldots,b_n$ der Durchschnitt $\bigcap_{i=1}^n \mid M_{\le b_i}\mid\not=\emptyset$,
 so gibt es ein $c\in M$ mit $c\le b_i$ f"ur alle $i=1,\ldots,n$. Die Menge $\{b_1,\ldots,b_n\}$ ist somit nach unten
 beschr"ankt, und nach Vorraussetzung existiert eine gr"o\3te gemeinsame untere 
 Schranke $d$ der $b_1,\ldots,b_n$ in $M$. Damit ist $\bigcap_{i=1}^n\mid M_{\le b_i}\mid=\mid M_{\le d}\mid$,
 also ein Kegel und daher zusammenziehbar. \hfill $\Box$\\

\subsection{Erl"auterung der Beweismethode anhand der Hough\-ton-Gruppen}

Wir erinnern kurz an die Definition der Houghton-Gruppen. F"ur eine nat"urliche Zahl $n\ge 1$
ist $H_n$ die Gruppe aller Permutationen $h$ der Menge $S={\rm I}\!{\rm N}\times\{1,\ldots,n\}$,
f"ur die ein $x_0\in {\rm I}\!{\rm N}$ sowie ein $n$-Tupel $(m_1,\ldots,m_n)\in {\sl Z}\!\!{\sl Z}^n$ 
existiert, soda\3 f"ur alle $i=1,\ldots,n$ gilt:
   \begin{eqnarray}  (x,i)\,h=(x+m_i,i)\quad\mbox{f"ur alle }x\ge x_0\,.
   \end{eqnarray}

In [Br1] wurde von K.S. Brown gezeigt, da\3 $H_n$ vom Typ {\sc (FP)}$_{n-1}$, jedoch nicht vom Typ {\sc (FP)}$_n$ und f"ur
 $n\ge 3$ endlich pr"asentiert ist. Zum Nachweis dieser Endlichkeitseigenschaften hat er das in Abschnitt 1.1 vorgestellte
 Kriterium benutzt (Korollar 1.3).\\
 Die von Brown angewandte Methode zur Herstellung der Voraussetzungen von Korollar 1.3
 wird sp"ater auf die in Abschnitt 2 definierten Gruppen $G_n$ "ubertragen. Um einen Einblick
 in das Verfahren zu erhalten, beschreiben wir im folgenden die wesentlichen Schritte. Einige hierbei aufgestellten 
 Behauptungen werden im weiteren Verlauf der Arbeit (zum Teil in einer allgemeineren Form) bewiesen.
 In diesen F"allen weisen wir auf die entsprechenden Stellen hin.
 \pagebreak

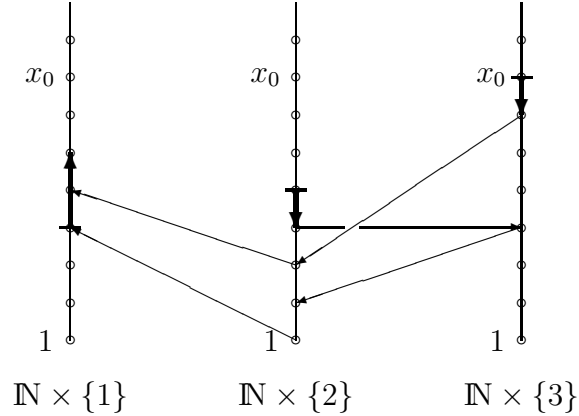
\begin{figure}[ht]
\unitlength 0.5cm
 \parbox[t]{15.3cm}
       {\small Die fettgedruckten Pfeile geben die Translation
         in der jeweiligen Menge $\scriptstyle {\rm I}\!{\rm N}\times \{i\},\: i=1,2,3$ an. Der Pfeilanfang
         befindet sich an der Stelle, wo die Translation beginnt. Die d"unngezeichneten Pfeile kennzeichnen 
         die restliche Permutation, wobei Elemente, von denen kein Pfeil ausgeht, auf sich selbst abgebildet werden.
       Hier ist $\scriptstyle m_1=2,\,m_2=m_3=-1$ und $\scriptstyle x_0=8$.}

\begin{center}
\begin{picture}(17,11)
\multiput(2,2)(6,0){3}{\begin{picture}(0,9)
                       \put(0,0){\line(0,1){9}}
                       \multiput(0,0)(0,1){9}{\circle{0.2}}
                       \put(0,0){\makebox(-1.3,0){1}}
                       \put(0,7){\makebox(-1.6,0){$x_0$}}
                       \end{picture}}
\put(8,2){\vector(-2,1){6}}
\put(8,4){\vector(-3,1){6}}
\put(9.7,5){\vector(1,0){4.3}}
\put(8,5){\line(1,0){1.3}}
\put(14,8){\vector(-3,-2){6}}
\put(14,5){\vector(-3,-1){6}}
\put(1,0.5){\makebox(2,0){${\rm I}\!{\rm N}\times \{1\}$}}
\put(7,0.5){\makebox(2,0){${\rm I}\!{\rm N}\times \{2\}$}}
\put(13,0.5){\makebox(2,0){${\rm I}\!{\rm N}\times \{3\}$}}
\thicklines
\multiput(1.975,5)(0.05,0){2}{\vector(0,1){2}}\multiput(1.725,5)(0.05,0){2}{\line(1,0){0.5}}
\multiput(7.975,6)(0.05,0){2}{\vector(0,-1){1}}\multiput(7.725,6)(0.05,0){2}{\line(1,0){0.5}}
\multiput(13.975,9)(0.05,0){2}{\vector(0,-1){1}}\multiput(13.725,9)(0.05,0){2}{\line(1,0){0.5}}
\end{picture}
\end{center}

 \caption{Beispiel f"ur ein Element aus $H_n$ ($n=3$)}
\end{figure}
{\bf Schritt 1:}$\;${\bf Konstruktion eines $H_n$-Komplexes}\\[0.4cm]
(i) {\em Die partiell geordnete Menge $\cal M$}\\[0.25cm]
Wir betrachten das Monoid $M$ aller injektiven Abbildungen von $S$ nach $S$, die (1.2) erf"ullen.
 ($M$ operiert von rechts, d.h. f"ur $\alpha,\beta$ aus $M$ bedeutet $\alpha
 \beta$ die Hintereinanderausf"uhrung von $\alpha$ gefolgt von $\beta$). Ein Element $t$ aus $M$ hei\3t 
 {\em Translation\/}, falls es ein $n$-Tupel $(m_1,\ldots,m_n)$ aus ${\sl Z}\!\!{\sl Z}^n$ gibt, soda\3 $(x,i)\,t=(x+m_i,i)$ f"ur {\em alle\/}
 $(x,i)$ aus $S$ erf"ullt ist. Mit $T$ bezeichnen wir das kommutative Untermomoid von $M$, das aus allen
 Translationen besteht. $T$ wird erzeugt von $t_1,\ldots,t_n$, wobei $t_i$ f"ur alle $x\in {\rm I}\!{\rm N}$ definiert ist durch
   \[ (x,j)\,t:\,=\left\{\begin{array}{l} (x+1,i)\quad\mbox{f"ur}\; j=i \\
                                          (x,j)\quad\mbox{f"ur}\; j\in \{1,\ldots,n\}-\{i\}
                         \end{array}\right. .   \]
Mit Hilfe von $T$ l"a\3t sich eine {\em Ordnung\/} auf $M$ erkl"aren: F"ur $\alpha$, $\beta$ aus $M$ schreiben
 wir $\alpha\le \beta$, falls ein $t$ aus $T$ existiert mit $t\alpha=\beta$
 (zur Verifikation der Ordnungsaxiome vgl. Abschnitt 3.1).\\
 Mit $\cal M$ bezeichnen wir die dem Monoid $M$ zugrundeliegende Menge zusammen mit dieser Ordnung. Man kann
 zeigen, da\3 es f"ur alle $\alpha$ aus $\cal M$ ein $t\in T$ gibt mit $\alpha\le t$. Daraus folgt
 die Gerichtetheit von $\cal M$ (vgl. Korollar 3.2), d.h. der der Menge $\cal M$ zugeordnete
 Simplizialkomplex $\mid \cal M\mid$ ist zusammenziehbar (s. Lemma 1.6).\\[0.4cm]
(ii) {\em Operation von $H_n$ auf $\mid\cal M\mid$ }\\[0.25cm]
$H_n$ operiert als Untergruppe von $\cal M$ durch Rechtsmultiplikation auf $\cal M$. Dabei gilt f"ur alle $\alpha$ aus $\cal M$, $t\in T$ und $h\in H_n$
    \[  (t\alpha)\, h=t(\alpha h)\,, \]
d.h. die Operation von $H_n$ auf $\cal M$ ist ordnungserhaltend. Nach Bemerkung 1.1 ist 
 $\mid \cal M\mid$ somit ein $H_n$-Komplex. Der Stabilisator eines Simplexes aus $\mid\cal M\mid$ ist 
 endlich (ohne Beweis), also endlich pr"asentiert und sogar vom Typ {\sc (FP)}$_{\infty}$.\\[0.8cm]
   {\bf Schritt 2:}$\;$ {\bf Filtrierung von $\mid\cal M\mid$}\\[0.4cm]
F"ur ein $\alpha$ aus $\cal M$ sei der {\em Grad\/} von $\alpha$, bezeichnet
 mit $gr(\alpha)$,
 die gr"o\3te ganze Zahl $k$, soda\3 es 
 eine Kette der Form $\alpha=\alpha_k>\alpha_{k-1}>\ldots>\alpha_0$ in $\cal M$ gibt.
 Offensichtlich ist f"ur $\alpha <\beta$ auch
 $gr(\alpha)<gr(\beta)$, d.h. die Funktion $gr:\cal M\longrightarrow {\rm I}\!{\rm N}$ ist streng wachsend.\\
Mit Hilfe der Teilmengen
   \[   {\cal M}_k:\,=\{ \alpha\in {\cal M}\mid gr(\alpha)\le k\}   \]
 f"ur ein $k$ aus ${\rm I}\!{\rm N}_0$ erh"alt man eine Folge von Unterkomplexen $\{\mid{\cal M}_k\mid\}_{k\ge 0}$
 mit der Eigenschaft $\mid{\cal M}_k\mid\subseteq \mid{\cal M}_{k+1}\mid$ und $\mid{\cal M}\mid=\bigcup_{k\ge 0}\mid{\cal M}_k\mid$.
Wegen $gr(\alpha h)=gr(\alpha)$ f"ur $\alpha\in{\cal M}$, $h\in H_n$ (vgl. Lemma 3.9), sind die einzelnen $\mid{\cal M}_k\mid$
$H_n$-invariant.\\[0.4cm]
 (iii) {\em Die Unterkomplexe der Filtrierung {\rm mod} $H_n$}\\[0.25cm]
Zum Nachweis der Endlichkeit von $\mid{\cal M}_k\mid$ mod $H_n$ vergleiche man Abschnitt 4.1; die dort
angewandte Argumentation ist im "ubertragenden Sinne auch hier m"oglich.\\[0.4cm]
(iv) {\em Homotopieeigenschaften der Inklusion $\mid{\cal M}_k\mid\,\subseteq\,\mid{\cal M}_{k+1}\mid$}\\[0.25cm]
Da zwei verschiedene Elemente $\alpha$, $\beta$ aus $\cal M$ mit $gr(\alpha)=gr(\beta)$ nicht
 vergleichbar sind ($gr$ ist streng wachsend),
 entsteht $\mid{\cal M}_{k+1}\mid$ aus $\mid {\cal M}_k\mid$, indem f"ur alle $\alpha$ aus $\cal M$ mit $gr(\alpha)=k+1$
 ein Kegel "uber $\mid{\cal M}_{<\alpha}\mid$ hinzugef"ugt wird. Gelingt es zu zeigen, da\3 f"ur $\alpha$ aus $\cal M$
 mit hinreichend gro\3em Grad $gr(\alpha)$ der Simplizialkomplex $\mid{\cal M}_{<\alpha}\mid$ homotopie"aquivalent
 zu einem Bouquet aus $(n-1)$-dimensionalen Sph"aren ist, so erh"alt man daher
$\mid{\cal M}_{k+1}\mid$ bis auf Homotopie aus $\mid {\cal M}_k\mid$ durch 
 Hinzuf"ugen von $n$-Zellen.\\
Brown hat nachgewiesen, da\3 $\mid{\cal M}_{<\alpha}\mid$ mit $gr(\alpha)=k$ homotopie"aquivalent zu folgendem kombinatorischen
 Simplizialkomplex $\Sigma_{n,k}$ ist:\\
F"ur eine Menge $W$ der Ordnung $k$ ist $\{1,\ldots,n\}\times W$ die Eckenmenge von $\Sigma_{n,k}$, und
 $((i_0,w_0),(i_1,w_1),\ldots,(i_p,w_p))$ bildet ein $p$-Simplex von $\Sigma_{n,k}$, falls
 $i_{\nu}\not=i_{\mu}$ und $w_{\nu}\not=w_{\mu}$ f"ur alle $\nu \not=\mu$, $\nu,\mu\in\{0,\ldots,p\}$ ist.\\
Das $1$-Ger"ust von $\Sigma_{n,k}$ ist ein $n$-gef"arbter Graph $\Gamma$, wobei
$\{i\}\times W$ aus den Ecken der Farbe $i$ besteht. Dabei gilt im Falle $k\ge 2n$ f"ur alle $i=1,\ldots,n$:
  \begin{quote}  Zu je $2(n-1)$ Ecken au\3erhalb von $\{i\}\times W$ gibt es in
                 $\{i\}\times W$ mindestens zwei Ecken, die zu jeder dieser $2(n-1)$ Ecken benachbart sind.
  \end{quote}
($\:2(n-1)$ Ecken beziehen $2(n-1)$ Elemente $w_1,\ldots,w_{2(n-1)}$ aus
 $W$ mit ein. Wegen \#$W=k\ge 2n$ existieren 2 verschiedene
  $w,w^{'}\in W-\{w_1,\ldots,w_{2(n-1)}\}$. Dann sind $(i,w)$ und $(i,w^{'})$
  die geforderten Ecken ). Mit Hilfe von Lemma 4.7
 ergibt sich daraus die Homotopie"aquivalenz zwischen $\Sigma_{n,k}$ und einem
 $(n-1)$-dimensionalen Sph"arenbouquet.\newpage

\section{Verallgemeinerung der Houghton-Gruppen} 
F\"ur eine nat\"urliche Zahl $n\ge1$ und f\"ur die Mengen  $X:\,={\rm I}\!{\rm N}$  
sowie $Y:\,={\rm I}\!{\rm N}$  betrachten wir die Menge  \[ S:\,=(X\times Y)\times \{1,\ldots,n\}\:,  \]
 welche aus der disjunkten Vereinigung von $n$ Kopien der Menge \,$X\times Y$\, besteht. Die i-te Kopie 
von \,$X\times Y$\, bezeichnen wir mit $Q_i\,$. Die Elemente aus $Q_i$ kann man
 sich als ganzzahlige Gitterpunkte eines Quadranten in einem kartesischen 
Koordinatensystem vorstellen (Abbildung 2).
\unitlength0.6cm
\begin{figure}[ht]
\begin{center}
\begin{picture}(11.3,11.3)
\thicklines
\put(0,0){\vector(1,0){10}}\put(0,0){\vector(0,1){10}}
\thinlines
\multiput(0.9,1)(1,0){9}
  {\begin{picture}(0,0)\multiput(0,0)(0,1){9}{.}\end{picture}}
\put(4,3){\circle*{0.2}}
\put(10,0){\makebox(0.8,0){x}}\put(0,10){\makebox(0,1){y}}
\multiput(1,0)(1,0){9}{\circle*{0.2}}\multiput(0,1)(0,1){9}{\circle*{0.2}}
\put(1,-1.2){\makebox(0,1.2){1}}\put(-0.8,0){\makebox(0.8,0){1}}
\put(2.5,3){\makebox(3,1){$\scriptstyle((x,y),\,i)$}}
\end{picture}\end{center}
\caption{Die Menge $Q_i=(X\times Y)\times\{i\}$}
\end{figure}
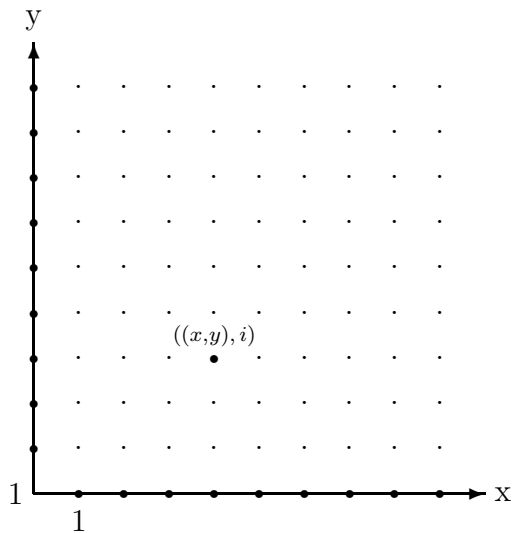
$S$ besitzt eine durch die Menge ${\rm I}\!{\rm N}$ induzierte Ordnung, d.h. es ist
$((x,y),i)\le((x',y'),i)$, falls $x\le x'$ und $y\le y'$.\\
Es werden im folgenden zwei m"ogliche Verallgemeinerungen der Houghton-Gruppen beschrieben.
In beiden F"allen handelt es sich um gewisse Permutationsgruppen der oben beschriebenen
Menge $S$, wobei analog zu den Houghton-Gruppen jede Permutation ab einer Stelle $p_0$
(d.h. f"ur alle $((x,y),i)\in S$ mit $(x,y)\ge p_0$) eine Translation ist.\\
Wir betrachten nun zum einen die Gruppen $G_n$, bei denen diese Translation in jedem Quadranten $Q_i$
 nur in einer Richtung, n"amlich in Richtung der Winkelhalbierenden von $Q_i$, stattfindet
(s. Abbildung 3, links), und zum anderen die Gruppen $\widetilde{G}_n$, bei denen als zus"atzliche
 Verallgemeinerung Translationen in jeder Richtung gestattet werden (s. Abbildung 3, rechts).\\
Neben der Forderung, ab einer gewissen Stelle eine Translation zu sein, m"ussen die Permutationen aus
$G_n$ bzw. $\widetilde{G}_n$ noch weitere Bedingungen erf"ullen, die jedoch f"ur beide Gruppentypen
"ubereinstimmen und im n"achsten Abschnitt erl"autert werden.
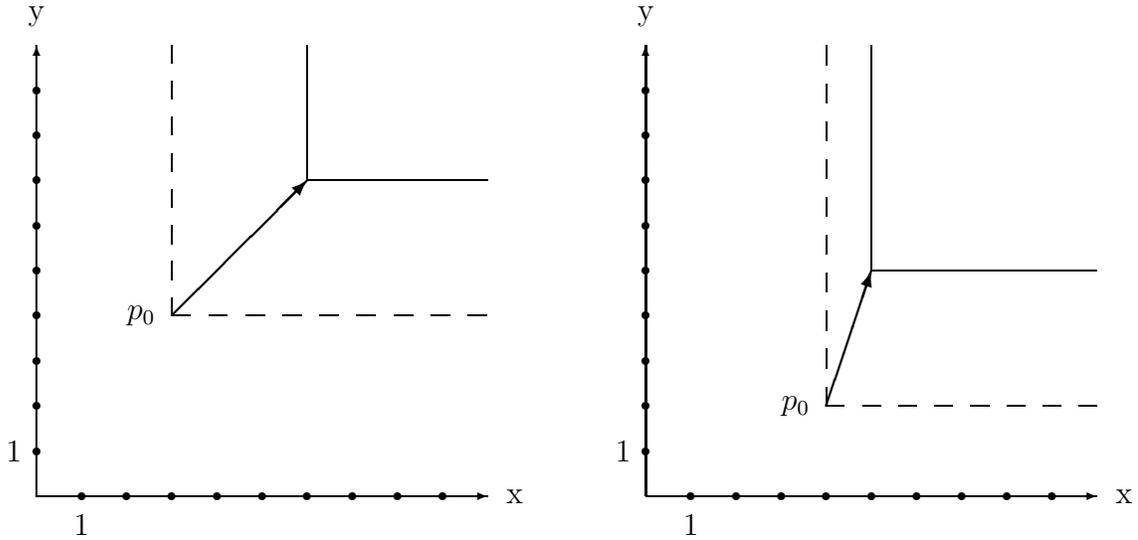
\begin{figure}[t]\unitlength 0.6cm
\begin{center}\parbox{14.4cm}{\small Die unterbrochenen Linien begrenzen den Bereich, welcher verschoben wird;
                                  der fettgedruckte Pfeil gibt die Translation an und das schraffierte Gebiet stellt
                                  den zugeh"origen Bildbereich dar. }
  \end{center}
\begin{center}\begin{picture}(24,13)
    \multiput(0,1)(13.5,0){2}{ \begin{picture}(10,10)
                             \put(0,0){\vector(1,0){10}}
                             \put(0,0){\vector(0,1){10}}
                             \multiput(1,0)(1,0){9}{\circle*{0.2}}
                             \multiput(0,1)(0,1){9}{\circle*{0.2}}
                             \put(10,0){\makebox(1.2,0){x}}
                             \put(0,10){\makebox(0,1.3){y}}
                             \put(1,-1.3){\makebox(0,1.3){1}}
                             \put(-1,1){\makebox(1,0){1}}
                             \end{picture} }
  \put(3,5){ \begin{picture}(7,6)
             \multiput(0,0)(0.82,0){9}{\line(1,0){0.41}} 
             \multiput(0,0)(0,0.8){8}{\line(0,1){0.4}}
             \put(3,3){\line(1,0){4}}
             \put(3,3){\line(0,1){3}}
             \thicklines\put(0,0){\vector(1,1){3}}
             \end{picture} }
  \put(17.5,3){ \begin{picture}(6,8)
              \multiput(0,0)(0.8,0){8}{\line(1,0){0.4}}
              \multiput(0,0)(0,0.84){10}{\line(0,1){0.42}}
              \put(1,3){\line(1,0){5}}
              \put(1,3){\line(0,1){5}}
              \thicklines\put(0,0){\vector(1,3){1}}
              \end{picture} } 
  \put(16.2,3){\makebox(1.7,0){$p_0$}}
  \put(1.7,5){\makebox(1.7,0){$p_0$}}
  \end{picture}\end{center}
  
  \caption{ Beispiel einer zul"assigen Translation bzgl. der Gruppe $G$ (links)
            und der Gruppe $\widetilde{G}_n$ (rechts) innerhalb eines Quadranten $Q_i$}
  \end{figure}
 \subsection{Beschreibung der Gruppen $G_n$ und $\protect{\widetilde{G}_n}$}
 F"ur ein $n\ge 1$ und $S=(X\times Y)\times \{1,\ldots,n\}$ sei $G_n$ die Gruppe aller Permutationen $g$ von $S$, f\"ur 
die ein $\:p_0=(x_0,y_0)\:$ aus $\,X\times Y\,$ existiert, soda\ss\ folgende
 Bedingungen erf\"ullt sind :
  \begin{itemize}
   \item[{\bf B$_1$}]:\hspace*{2ex} Es gibt ein $n$-Tupel $\:(m_1,\ldots,m_n)\:$ aus 
              ${\sl Z}\!\!{\sl Z}^n$, soda\3 f"ur alle $\,(x,y)\,\ge\,p_0\,$ und 
              \hspace*{2.35ex} $\,i=1,\ldots,n\,$ gilt:
                  \[  ((x,y),i)\,g\,=\,((x+m_i\,,y+m_i),i)   \]
   \item[{\bf B$_2$}]$\!${\bf (a)}$\,$:\hspace*{2ex} F"ur alle $(x,i)$ aus $X\times \{1,\ldots,n\}$ gibt es ein 
             $q_{(x,i)}$ aus ${\sl Z}\!\!{\sl Z}$ sowie ein \hspace*{5.9ex} $(x,i)'=(x'\!,i')$  aus $X\times \{1,\ldots,n\}$, 
             soda\3 f"ur alle $y\ge y_0$ gilt:
                 \[ ((x,y),i)\,g\,=\,((x',y+q_{(x,i)}),i') \]
   \item[{\bf B$_2$}]$\!${\bf (b)}$\,$:\hspace*{2ex} F"ur alle $(y,i)$ aus $Y\times \{1,\ldots,n\}$ gibt es ein
             $r_{(y,i)}$ aus ${\sl Z}\!\!{\sl Z}$ sowie ein \hspace*{6ex} $(y,i)'=(y'\!,i')$ aus $Y\times \{1,\ldots,n\}$,
             soda\3 f"ur alle $x\ge x_0$ gilt:
                 \[  ((x,y),i)\,g\,=\,((x+r_{(y,i)}),i')  \]
 \end{itemize}
 \pagebreak
{\bf Bemerkung:}$\:$F"ur $(x,i)\in X\times \{1,\ldots,n\}$ ist nur das Element
 $(x,i)'$ eindeutig bestimmt, jedoch nicht die einzelnen Komponenten von 
 $(x,i)'$. Es ist m"oglich, da\3 f"ur Elemente $(x_1,i_1)$,$(x_2,i_2)$ mit
 $i_1=i_2$ trotzdem $i_1{}^{'}\not=i_2{}^{'}$, bzw. f"ur $x_1=x_2$ $x_1{}^{'}
 \not=x_2{}^{'}$ ist. Dasselbe gilt nat"urlich auch f"ur die Elemente $(y,i)$
 aus $Y\times \{1,\ldots,n\}$.\\\\
Die Gruppen $\widetilde{G}_n$ erh"alt man, indem B$_1$ ersetzt wird durch die Bedingung
\begin{itemize} 
   \item[{\bf $\widetilde{\rm\bf B}_1$}]:\hspace*{2ex} Es gibt ein $n$-Tupel $(m_1,\ldots,m_n)=((m_{11},m_{12}),\ldots,(m_{n1},m_{n2}))$ 
       aus $({\sl Z}\!\!{\sl Z}^2)^n$, \hspace*{2.75ex} soda\3 f"ur alle $(x,y)\ge p_0$ und $i=1,\ldots,n$ gilt:
             \[  ((x,y),i)\,g\,=\,((x,y)+m_i,i) \:, \]
  \end{itemize}
 und B$_2$(a) bzw. B$_2$(b) unver"andert bleiben.\\\\
{\bf Erl"auterungen.} Die Forderung B$_1$ bzw. $\widetilde{\rm B}_1$ bewirkt, da\3 $g$ ab $p_0$ in jedem Quadranten $Q_i$ eine Translation
 ist (mit Translationsvektor $(m_i,m_i)$ bzw. $(m_{i1},m_{i2})$, zur Veranschaulichung s. Abbildung 3). Notwendigerweise mu\3 
 dabei $\sum_{i=1}^n m_i=0$ bzw. $\sum_{i=1}^n (m_{i1},m_{i2})=(0,0)$ gelten.\\

Durch die Einschr"ankung B$_1$ bzw $\widetilde{\rm B}_1$ alleine werden jedoch den Permutationen $g$ noch zu
 viele Freiheiten gelassen. Vergleicht man die Situation mit der im Falle der Houghton-Gruppen, so stellt man fest, da\3 dort nur eine endliche Menge
 "`willk"urlich"' abgebildet werden darf.
 Um ein "ahnliches Verhalten hinsichtlich der Permutationen aus $G_n$ bzw. $\widetilde{G}_n$ zu erreichen,
 sind weitere Forderungen n"otig, die jedoch mit B$_1$ bzw. $\widetilde{\rm B}_1$ vertr"aglich sein m"ussen.
 Diese Forderungen werden in B$_2$(a) und B$_2$(b) festgelegt.\\

Um B$_2$(a) geometrisch zu interpretieren, definieren wir f"ur ein $(x,i)\in X\times\{1,\ldots,n\}$ und $k$ aus ${\rm I}\!{\rm N}$
die Menge   \[((x,i),\ge k):\,=\{((x,y),i)\in S\mid y\ge k\}  \]
und betrachten die Mengen $((x,i),\ge y_0)$ f"ur alle $(x,i)\in X\times\{1,\ldots,n\}$.
Die Forderung B$_2$(a) bedeutet, da\3 diese Mengen durch $g$ jeweils erst ordnungserhaltend
auf eine andere Menge $((x,i)',\ge y_0)$ dieser Form abgebildet und anschlie\3end um den Wert $q_{(x,i)}$ aus 
${\sl Z}\!\!{\sl Z}$ verschoben werden (s. Abbildung 4, zur Vereinfachung wurde in diesem Beispiel $i=i'$ gew"ahlt).\\
 Die Permutation $g$ wird also auf der Menge $((x,i),\ge y_0)$ durch Angabe von $(x,i)'$ und $q_{(x,i)}$ festgelegt.
 Dabei ist zu beachten, da\3 f"ur $(x,i)\in X\times\{1,\ldots,n\}$ mit
 $x\ge x_0$ die Elemente $(x,i)'$ sowie $q_{(x,i)}$ 
durch $m_i\in{\sl Z}\!\!{\sl Z}$ im Falle $g\in G_n$ bzw.
 durch $(m_{i1},m_{i2})\in{\sl Z}\!\!{\sl Z}^2$ im Falle
 $g\in\widetilde{G}_n$ bestimmt sind:\\
 F"ur $g\in G_n$ gilt $(x,i)'=(x'\!,i')=(x+m_i,i)$ und $q_{(x,i)}=m_i\,$; im Falle
 $g\in\widetilde{G}_n$ ist dagegen $(x'\!,i')=(x+m_{i1},i)$ und $q_{(x,i)}=m_{i2}$
 (s. Abbildung 5).

 \begin{figure}[ht]
  \unitlength 0.6cm
  \begin{center}\begin{picture}(11,12)
      \put(0,0.5){\begin{picture}(10,10)
                  \put(0,0){\vector(1,0){10}}
                  \put(0,0){\vector(0,1){10}}
                  \multiput(1,0)(1,0){9}{\circle*{0.2}}
                  \multiput(0,1)(0,1){9}{\circle*{0.2}}
                  \put(10,0){\makebox(1.2,0){x}}
                  \put(0,10){\makebox(0,1.3){y}}
                  \put(1,-1.3){\makebox(0,1.3){1}}
                  \put(-1,1){\makebox(1,0){1}}
                  \thicklines
                  \put(3,3){\line(0,1){7}}\put(2.75,3){\line(1,0){0.5}}
                  \put(7,5){\line(0,1){5}}\put(6.75,5){\line(1,0){0.5}}
                  \thinlines
                  \put(3.6,3){\vector(1,0){2.8}}\put(7.4,3.25){\vector(0,1){1.5}}
                  \put(3,1.5){\makebox(0,1.5){$\scriptstyle ((x,i),\ge y_0)$}}
                  \put(7,1.5){\makebox(0,1.5){$\scriptstyle ((x,i)\:',\ge y_0)$}}
                  \multiput(7,3)(0,1){2}{\line(0,1){0.5}}
                  \put(6.75,3){\line(1,0){0.5}}
                  \put(-1.2,3){\makebox(1.2,0){$y_0$}}
                  \put(-2.8,5){\makebox(2.8,0){$\scriptstyle y_0+q_{(x,i)}$}}
                  \put(3,-1.2){\makebox(0,1.2){$x$}}
                  \put(4,-1.3){\makebox(0,1.3){$x_0$}}
                  \put(7,-1.1){\makebox(0,1.1){$x'$}}
                 \end{picture}  }
   \end{picture}\end{center}
  \caption{Veranschaulichung der Bedingung B$_2$(a) }

  \begin{center}\begin{picture}(26,15)
    \multiput(0,1)(15,0){2}{\begin{picture}(10,10)
                            \thicklines\put(0,0){\vector(1,0){10}}
                            \put(0,0){\vector(0,1){10}} 
                            \thinlines
                            \multiput(1,0)(1,0){9}{\circle*{0.2}}
                            \multiput(0,1)(0,1){9}{\circle*{0.2}}
                            \put(10,0){\makebox(1.2,0){x}}
                            \put(0,10){\makebox(0,1.3){y}}
                            \put(1,-1.3){\makebox(0,1.3){1}}
                            \put(-1,1){\makebox(1,0){1}}
                            \end{picture}} 
    \put(0,1){\begin{picture}(10,10)
              \thicklines
              \put(4,2){\line(0,1){8}}\put(3.75,2){\line(1,0){0.5}}
              \put(7,5){\line(0,1){5}}\put(6.75,5){\line(1,0){0.5}}
              \put(4.2,2.2){\vector(1,1){2.8}}\thinlines
              \put(4.4,2){\vector(1,0){2.5}}\put(7,2){\vector(0,1){3}}
              \put(2,0){\makebox(0,-1.3){$x_0$}}
              \put(4,0){\makebox(0,-1.2){$x$}}
              \put(7,0){\makebox(0,-1.2){$\scriptstyle x+m_i$}}
              \put(-1.2,2){\makebox(1.2,0){$y_0$}}
              \put(-2.3,5){\makebox(2.6,0){$\scriptstyle y_0+m_i$}}
              \end{picture} }
   \put(15,1){\begin{picture}(10,10)
               \thicklines\put(4,3){\line(0,1){7}}\put(3.75,3){\line(1,0){0.5}}
               \put(7,5){\line(0,1){5}}\put(6.75,5){\line(1,0){0.5}}
               \put(4.2,3.2){\vector(3,2){2.8}}\thinlines
               \put(4.4,3){\vector(1,0){2.5}}\put(7,3){\vector(0,1){2}}
               \put(2,-1.3){\makebox(0,1.3){$x_0$}}
               \put(4,-1.2){\makebox(0,1.2){$x$}}
               \put(7,-1.2){\makebox(0,1.2){$\scriptstyle x+m_{i1}$}}
               \put(-1.2,3){\makebox(1.2,0){$y_0$}}
               \put(-2.4,5){\makebox(2.9,0){$\scriptstyle y_0+m_{i2}$}}
              \end{picture}  }
  \end{picture}
  \begin{minipage}[t]{15cm}
       \begin{minipage}{6cm}\begin{center} $g\in G_n$ \end{center}
       \end{minipage}\hfill
       \begin{minipage}{5cm}\begin{center} $g\in\widetilde{G}_n$ \end{center}
       \end{minipage}
  \end{minipage}
 \end{center}
 \caption{ Die Bedingung B$_2$(a) f"ur Elemente $(x,i)$ mit $x\ge x_0$}
 \end{figure}
 \clearpage
 Die Bedingung B$_2$(b) l"a\3t sich analog interpretieren mit Hilfe der Mengen
     \[ ((y,i)\ge k):\,=\{ ((x,y),i)\in S\mid x\ge k\}\:,  \]
wobei $y\in Y\times\{1,\ldots,n\}$ und $k\in {\rm I}\!{\rm N}$ ist. Auf der Menge 
$((y,i),\ge x_0)$ wird die Pemutation $g$ durch $(y,i)'\in Y\times\{1,\ldots,n\}$ und 
$r_{(y,i)}\in{\sl Z}\!\!{\sl Z}$ festgelegt, wobei nun f"ur alle $y\ge y_0$ gilt:
 $(y,i)'=(y'\,i')=(y+m_i\,,i)$ und $r_{(y,i)}=m_i$ falls $g\in G_n$, dagegen $(y'\,i')=(y+m_{i2},i)$
 und $r_{(y,i)}=m_{i1}$ falls $g\in \widetilde{G}_n$.\\\\

 {\bf Bemerkung 2.1:}$\;$ F"ur $(x_1,i_1),(x_2,i_2)\in X\times\{1,\ldots,n\}$
 mit $(x_1,i_1)\not=(x_2,i_2)$ ist auch $(x_1,i_1)'\not=(x_2,i_2)'$,
 d.h. $g\in \widetilde{G}_n$ induziert  eine Permutation
    \[     \pi\,(g):\, X\times\{1,\ldots,n\}\,\longrightarrow\, X\times\{1,\ldots,n\}\]
         \[ (x,i)\,\pi(g):\,=(x,i)'\:.\]
 $\pi(g)$ beschreibt sozusagen die durch $g$ verursachte Permutation der
 Halbgeraden $((x,i),\ge y_0)$ (für $\,(x,i)\in X\times \{1,\ldots,n\}\,$) unter
 Vernachl"assigung der Verschiebung (um $q_{(x,i)}$) in Richtung der y-Achse.
 Dabei gilt f"ur alle $x\ge x_0$: $(x,i)'=(x+m_{i1})$ mit $m_{i1}\in{\sl Z}\!\!{\sl Z}$ (vgl. Abschnitt 2.1, Erl"auterungen).
 D.h. $\pi(g)$ ist ein Element der Houghtongruppe $H_n$ bez"uglich der zugrundeliegenden Menge
 $X\times\{1,\ldots,n\}$. Analog erh"alt man zu $g$ eine Permutation $\sigma(g)$ der Menge $Y\times\{1,\ldots,n\}$
      \[  \sigma(g):\, Y\times\{1,\ldots,n\}\,\longrightarrow\, Y\times\{1,\ldots,n\} \]
            \[   (y,i)\,\sigma(g):\,=(y,i)'\:,\]
 wobei $(y,i)'$ aus Bedingung B$_2$(b) stammt und f"ur $y\ge y_0$ die Beziehung
  $(y,i)'=(y+m_{i2},i)$ mit $(m_{12},\ldots,m_{n2})\in{\sl Z}\!\!{\sl Z}^n$ besteht.
  $\sigma(g)$ ist also ebenfalls aus $H_n$ (bez"uglich
 zugrundeliegender Menge $Y\times \{1,\ldots,n\}$).\\
Die Zuordnungen $g\,\longmapsto\,\pi(g)$ bzw. $g\,\longmapsto\,\sigma(g)$, die wir im weiteren mit $\pi$ bzw. $\sigma$ bezeichnen, sind daher jeweils
 Homomorphismen von der Gruppe $\widetilde{G}_n$ nach $H_n$. 
 Dabei ist schon die Einschr"ankung der Abbildung $\pi$ bzw. $\sigma$ auf die Gruppe
 $G_n$  surjektiv, d.h. die Houghton-Gruppe $H_n$ taucht als Faktorruppe von
 $G_n$ und $\widetilde{G}_n$ auf.\\

Abbildung 6 zeigt ein Beispiel f"ur eine Permutation $g$ aus $\widetilde{G}_2$. Daran wird auch das 
prinzipielle Verhalten der Permutationen aus $G_n$ deutlich, soda\3 wir auf ein zus"atzliches Beispiel
 f"ur ein $g\in G_n$ verzichten.

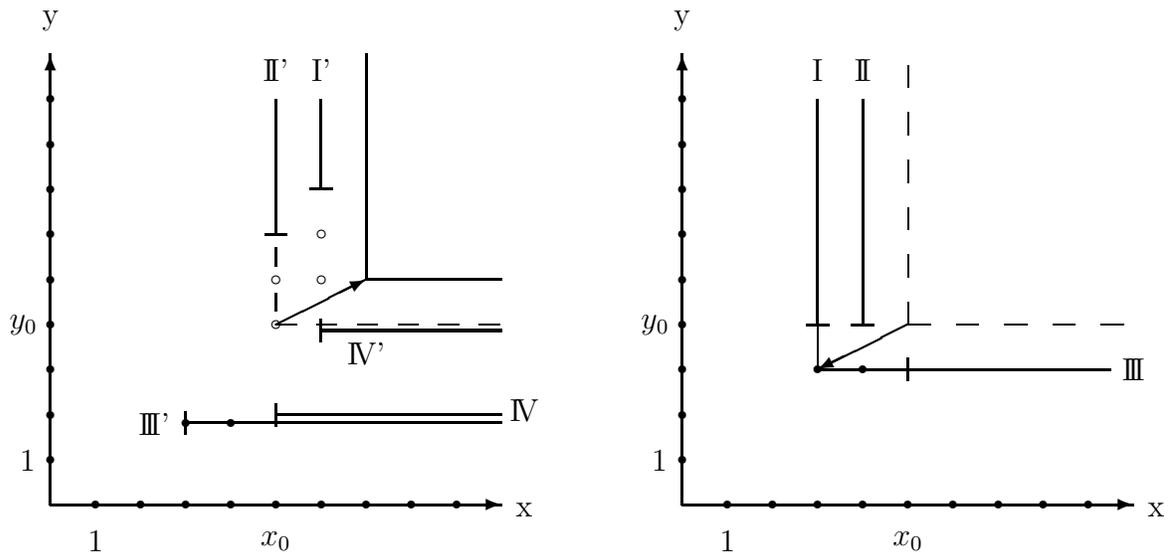
\begin{figure}[ht]\unitlength0.6cm
\begin{picture}(30,14)
\put(-2.4,0)
{\begin{picture}(30,14)
\multiput(3,2)(14,0){2}{\begin{picture}(11,12)
\put(0,0){\thicklines\vector(0,1){10}}\put(0,0){\thicklines\vector(1,0){10}}
\multiput(1,0)(1,0){9}{\circle*{0.2}}
\multiput(0,1)(0,1){9}{\circle*{0.2}}\put(10,-1.1){\makebox(1,2){x}}
\put(0,10){\makebox(0,1.5){y}}\put(0,-1.3){\makebox(2,1){1}}
\put(-1,0){\makebox(1,2){1}}\put(4,-1.3){\makebox(2,1){$x_0$}}
\put(-1.1,3){\makebox(1,2){$y_0$}}
\end{picture}}
\put(6,3.82){\thicklines\line(1,0){7}}\put(8,4){\thicklines\line(1,0){5}}
\put(9,5.87){\thicklines\line(1,0){4}}\put(8,8){\thicklines\line(0,1){3}}
\put(9,9){\thicklines\line(0,1){2}}
\put(10,7){\line(0,1){5}}\put(10,7){\line(1,0){3}}
\put(8,6){\thicklines\vector(2,1){2}}
\put(20,6){\thicklines\line(0,1){5}}
\put(21,6){\thicklines\line(0,1){5}}
\put(22,5){\thicklines\line(1,0){4.5}}
\put(22,6){\thicklines\vector(-2,-1){2}}
\put(20,5){\line(0,1){1}}\put(20,5){\line(1,0){2}}
\put(6,3.57){\thicklines\line(0,1){0.5}}
\put(8,3.75){\thicklines\line(0,1){0.5}}
\put(9,5.62){\thicklines\line(0,1){0.5}}
\put(7.75,8){\thicklines\line(1,0){0.5}}
\put(8.75,9){\thicklines\line(1,0){0.5}}
\put(19.75,6){\thicklines\line(1,0){0.5}}
\put(20.75,6){\thicklines\line(1,0){0.5}}
\put(22,4.75){\thicklines\line(0,1){0.5}}
\multiput(8,6)(0.9,0){6}{\line(1,0){0.45}}
\multiput(8,6.3)(0,1){2}{\line(0,1){0.4}}
\multiput(22,6)(1.07,0){5}{\line(1,0){0.5}}
\multiput(22,6)(0,1.05){6}{\line(0,1){0.5}}
\put(20,11){\makebox(0,1.3){I}}
\put(21,11){\makebox(0,1.3){I$\!$I}}
\put(8,11){\makebox(0,1.3){I$\!$I'}}
\put(9,11){\makebox(0,1.3){I'}}
\put(9,5.37){\makebox(2,0){I$\!$V'}}
\put(4.3,3.8){\makebox(2,0){I$\!$I$\!$I'}}\put(12.5,4.1){\makebox(2,0){I$\!$V}}
\put(26,5){\makebox(2,0){I$\!$I$\!$I}}
\multiput(8,7)(1,0){2}{\circle{0.2}}
\put(9,8){\circle{0.2}}\put(8,6){\circle{0.2}}
\multiput(6,3.82)(1,0){2}{\circle*{0.2}}
\multiput(20,5)(1,0){2}{\circle*{0.2}}
\end{picture}}\end{picture}
\parbox{15.3cm}{$p_0=(5,4)\:,\: m_1=(2,1)\:,\: m_2=(-2,-1)\:,\:q_{(3,2)}=3\:,\:
q_{(4,2)}=2\:,\:r_{(2,1)}=-1\:,\:r_{(3,2)}=-2$\\\\}
\parbox{15.3cm}{\small Alle Elemente $\scriptstyle ((x,y),i)$ mit
 $\scriptstyle (x,y)\ge p_0=(x_0,y_0)$ werden gem"a\3 der Pfeilangabe verschoben.
 Die Mengen der Form $\scriptstyle(\,(x,i)\,,\,\ge 4)$ mit 
$\scriptstyle x<5$ bzw. $\scriptstyle(\,(y,i)\,,\,\ge 5)$ mit 
$\scriptstyle y<4$, welche unter $g$ nicht identisch abgebildet
 werden, sind fettgedruckt und mit den R"omischen Ziffern  I, II, III, und IV  gekennzeichnet.
 Die zugeh"origen Bildmengen unter $g$ besitzen jeweils die Bezeichnung
 I ', II', III' und IV'. Die ausgef"ullten Kreise stellen die Elemente dar,
 denen noch ein Bildelement zugewiesen werden mu\3, dagegen markieren die
leeren Kreise die Stellen aus $\scriptstyle S$, die als Bild\-elemente
 in Frage kommen.\\} \normalsize 
\caption{ Beispiel einer Permutation $g$ aus $\widetilde{G}_2$ }
\end{figure}
\clearpage

\subsection{ $G_n$ als Normalteiler von $\protect{\widetilde{G}_n}$ }

Mit Hilfe des $n$-Tupels $((m_{11},m_{12}),\ldots,(m_{n1},m_{n2}))$ aus $({\sl Z}\!\!{\sl Z}^2)^n$,
 das nach Bedingung $\widetilde{\rm B}_1$ zu jedem $g\in \widetilde{G}_n$ existiert,
 definieren wir den Homomorphismus
    \[ \phi : \widetilde{G}_n\,\longrightarrow\,{\sl Z}\!\!{\sl Z}^n\]
            \[         \phi(g):\,= (m_{11}-m_{12},\ldots,m_{n1}-m_{n2})\:.\]
 Wegen $\sum_{i=1}^n (m_{i1},m_{i2})=(0,0)$, also auch $\sum_{i=1}^n (m_{i1}-m_{i2})=0$,
 ist das Bild von $\phi$ isomorph zur Gruppe ${\sl Z}\!\!{\sl Z}^{n-1}$. Der Kern von 
 $\phi$ besteht aus den Permutationen $g\in \widetilde{G}_n$ mit $m_{i1}=m_{i2}$ f"ur alle $i=1,\ldots,n$.
 Also ist Ker $\phi$ gerade die Gruppe $G_n$. Wir erhalten somit eine exakte Sequenz von Gruppen
     \[ 0\:\longrightarrow \:G_n\:\longrightarrow\:\widetilde{G}_n\:\longrightarrow\:{\sl Z}\!\!{\sl Z}^{n-1} \:\longrightarrow\: 0\:. \]
 Da ${\sl Z}\!\!{\sl Z}^{n-1}$ f"ur alle $n\ge 1$ vom Typ {\sc (FP)}$_{\infty}$ ist (vgl. [Br2], Kapitel 8, Abschnitt 6),
 "ubertragen sich nach Lemma 1.5 von denen in dieser Arbeit nachgewiesenen
  Endlichkeitseigenschaften der Gruppe $G_n$ (vom Typ {\sc (FP)}$_{n-1}$, nicht vom Typ {\sc (FP)}$_n$ und f"ur $n\ge 3$ endlich
 pr"asentiert) die positiven Teile auf 
 $\widetilde{G}_n$.
 \vfill
\subsection{Die Houghton-Gruppe als Untergruppe von $G_n$}
 \vfill
 Die Houghtongruppe taucht in mehrfacher Form als Untergruppe von $G_n$ (und
 auch von $\widetilde{G}_n$) auf. Es gibt sogar zu beliebigem $m$ aus 
 ${\rm I}\!{\rm N}$ eine Einbettung $H_m\hookrightarrow G_n$. Wir stellen in
 diesem Abschnitt spezielle Untergruppen $U$ von $G_n$ vor mit $H_m\hookrightarrow
 U$ f"ur ein $m\in{\rm I}\!{\rm N}$, und deren Faktorgruppe nach $H_m$ endlich ist.
 Solche Untergruppen spielen im weiteren Verlauf der Arbeit eine wichtige
 Rolle, da sie als Stabilisatoren bez"uglich der $G_n$-Operation auf dem im
 3. Abschnitt konstruierten Simplizialkomplex auftreten. \\
  \vfill
 Wir betrachten zu zwei
 positiven ganzen Zahlen $k_1$ und $k_2$ mit $k_1+k_2\ge 1$ $\;k_1$ viele 
 verschiedene Elemente $(x_1,i_1),(x_2,i_2),\ldots,(x_{k_1},i_{k_1})$ aus $X\times \{1,\ldots,n\}$
 sowie $k_2$ viele verschiedene Elemente $(y_1,j_1),(y_2,j_2),\ldots,(y_{k_2},j_{k_2})$
 aus $Y\times \{1,\ldots,n\}$, w"ahlen $m_1,\ldots,m_{k_1}$ sowie $l_1,\ldots,l_{k_2}$
 aus ${\rm I}\!{\rm N}$ und bilden die Menge\vfill 
       \[{A:\,=((x_1,i_1),\ge m_1)\cup\ldots\cup
((x_{k_1},i_{k_1}),\ge m_{k_1})\atop\quad\qquad\cup\:((y_1,j_1),\ge l_1)\cup\ldots\cup
((y_{k_2},j_{k_2}),\ge l_{k_2})\:.}\]\vfill
( Zur Definition der Mengen $((x,i),\ge m)$ bzw. $((y,j),\ge l)$ mit $m$ und
 $l$ aus ${\rm I}\!{\rm N}$ vergleiche man Abschnitt 2.1 ). Dabei seien die Mengen $((x_{\nu},i_{\nu}),\ge m_{\nu})$
 f"ur $\nu =1,\ldots,k_1$ und $((y_{\mu},j_{\mu})\ge l_{\mu})$ f"ur $\mu =1,\ldots,k_2$
 paarweise disjunkt zueinander (Abbildung 7 zeigt ein Beispiel f"ur $A$).
 Wir k"onnen also "uber die Zuordnung
 \pagebreak
\[ f:\,A\longrightarrow\,{\rm I}\!{\rm N}\times \{1,\ldots,k_1+k_2\}\]
  \[ \begin{array}{rcll}
    f\,((x_{\nu},i_{\nu}),y-1+ m_{\nu})&=&(y,\nu)\quad & \mbox{f"ur alle}\: y\in {\rm I}\!{\rm N},\,\nu=1,\ldots,k_1\\
      f\,((y_{\mu},j_{\mu}),x-1+l_{\mu})&=&(x,k_1+\mu)\quad & \mbox{f"ur alle}\:x\in {\rm I}\!{\rm N},\mu=1,\ldots,k_2
   \end{array}\]
 die Menge $A$ mit ${\rm I}\!{\rm N}\times \{1,\ldots,k_1+k_2\}$ identifizieren. Es sei nun
 $U_{k_1+k_2}$ die Untergruppe von $G_n$, welche s"amtliche Elemente
 aus $G_n$ enth"alt, die au\3erhalb von $A$ die Identit"at sind.
 Dies "ubertr"agt sich insofern auf die in Bemerkung 2.1 erkl"arten
 Abbildungen $\pi$ bzw. $\sigma$, da\3
 f"ur ein $g$ aus $U_{k_1+k_2}$ das Element $\pi\,(g)$ bzw. $\sigma\,(g)$
 au\3erhalb der endlichen Menge $\bigcup_{\nu=1}^{k_1} (x_{\nu},i_{\nu})$
 bzw. $\bigcup_{\mu=1}^{k_2} (y_{\mu},j_{\mu})$ ebenfalls die Identit"at ist.
 ($\pi(g)$ bzw. $\sigma(g)$ beschreibt die durch $g$ induzierte Permutation
 auf der Menge $X\times\{1,\ldots,n\}$ bzw. $Y\times\{1,\ldots,n\}$).
 $\pi\,(g)$ bzw. $\sigma\,(g)$ kann daher als Element von $\Sigma_{k_1}$ bzw.
 $\Sigma_{k_2}$ ( mit $\Sigma_{k}$ bezeichnen wir die Permutationsgruppe der Menge $\{1,\ldots,k\}$) aufgefa\3t werden, und wir erhalten einen Homomorphismus
 \[\begin{array}{c}\Phi\,:U_{k_1+k_2}\,\longrightarrow\,\Sigma_{k_1}\times \Sigma_{k_2}\\\\
     \Phi\,(g):\,=(\pi\,(g),\sigma\,(g))\:. \end{array}\]
 Im Kern von $\Phi$ liegen diejenigen Permutationen $g$ aus $G_n$,
 die nur noch Elemente innerhalb $A$ permutieren und dabei jede der 
 Halbgeraden $((x_{\nu},i_{\nu}),\ge m_{\nu})$ und

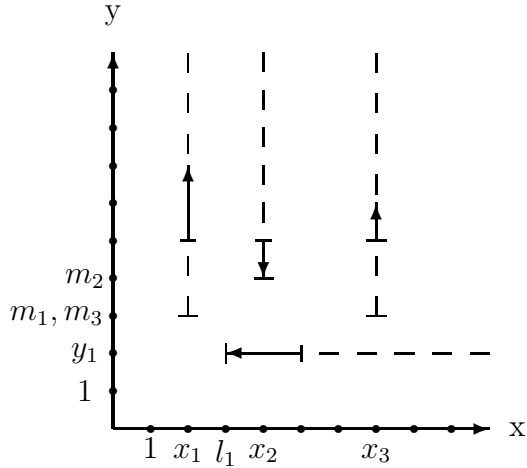
\unitlength0.5cm
\begin{figure}[ht]
\begin{center}
    \begin{picture}(12,11.5)
       \put(0,0.7){\begin{picture}(10,10)\thicklines
            \put(0,0){\vector(1,0){10}}
            \put(0,0){\vector(0,1){10}}
            \put(2,5){\vector(0,1){2}}
            \put(4,5){\vector(0,-1){1}}
            \put(7,5){\vector(0,1){1}}
            \put(5,2){\vector(-1,0){2}}
            \put(1.8,5){\line(1,0){0.4}}
            \put(3.8,5){\line(1,0){0.4}}
            \put(6.8,5){\line(1,0){0.4}}
            \put(5,1.8){\line(0,1){0.4}}
            \thinlines
            \multiput(1,0)(1,0){9}{\circle*{0.2}}
            \multiput(0,1)(0,1){9}{\circle*{0.2}}
            \multiput(2,3)(0,1.08){7}{\line(0,1){0.5}}
            \multiput(4,5.5)(0,1){5}{\line(0,1){0.5}}
            \multiput(7,3)(0,1.08){7}{\line(0,1){0.5}}
            \multiput(5.5,2)(1,0){5}{\line(1,0){0.5}}
            \put(10,0){\makebox(1.5,0){x}}
            \put(0,10){\makebox(0,2){y}}
            \put(1,0){\makebox(0,-1){1}} 
            \put(2,0){\makebox(0,-1.3){$x_1$}}
            \put(3,0){\makebox(0,-1.3){$l_1$}}
            \put(4,0){\makebox(0,-1.3){$x_2$}}
            \put(7,0){\makebox(0,-1.3){$x_3$}}
            \put(-1.5,1){\makebox(1.5,0){1}}
            \put(-1.5,2){\makebox(1.5,0){$y_1$}}
            \put(-3,3){\makebox(3,0){$m_1,m_3$}}
            \put(-1.5,4){\makebox(1.5,0){$m_2$}}
            \put(1.75,3){\line(1,0){0.5}}
            \put(3.75,4){\line(1,0){0.5}}
            \put(6.75,5){\line(1,0){0.5}}
            \put(6.75,3){\line(1,0){0.5}}
            \put(3,1.75){\line(0,1){0.5}}
             \end{picture}}
            \end{picture}\end{center}
            \parbox{15.3cm}
            {\small Zur Vereinfachung der Darstellung nehmen wir an, da\3 die Menge
            $\scriptstyle A$ sich innerhalb nur eines Quadranten $\scriptstyle
            Q_i$ befindet. Sie wird in der Abbildung durch die gestrichelten
            Linien gekennzeichnet. Die Pfeile stellen die durch $\scriptstyle
            q_1,q_2,q_3$ und $\scriptstyle r_1$ gegebene Translation dar.}
 \caption{ Beispiel f"ur ein $g\in$Ker$\:\Phi $}
\end{figure}

 $((y_{\mu},j_{\mu}),\ge l_{\mu})$,
 aus denen sich $A$ zusammensetzt, ab einer Stelle in sich verschieben (vgl. Abbildung 7).
 Anschaulich ist klar, da\3 Ker $\Phi$ isomorph zur Houghtongruppe $H_{k_1+k_2}$ 
 ist. Formal erh"alt man den Isomorphismus durch die Zuordnung
 $\:g\longmapsto f^{-1}\,gf\:$.\\

 Da die Faktorgruppe von $U_{k_1+k_2}$ nach der Untergruppe Ker$\,\Phi$ endlich
 ist, erhalten wir mit Hilfe von Lemma 1.4 folgendes Resultat:
 \begin{lemma} Mit $m:\,=k_1+k_2$ ist die Gruppe $U_{k_1+k_2}$ vom Typ $(FP)_{m-1}$,
  jedoch nicht vom Typ $(FP)_m$. F"ur $m\ge 3$ ist sie zus"atzlich endlich pr"asentiert.
 \end{lemma}
 {\bf Bemerkung 2.2} Die hier vorgestellten Gruppen $U_{k_1+k_2}$ treten in
 sp"ateren Situationen in einer etwas ver"anderten Form auf. Und zwar vergr"o\3ert
 sich dort die auf der vorigen Seite beschriebene Menge $A$ um eine
 endliche, zu $A$ disjunkte Menge $P=\{p_1,\ldots,p_r\}$ (s. Abbildung 8).
\begin{figure}[ht]\unitlength0.5cm
\begin{center}
\begin{picture}(11,11)
\thicklines
  \put(0,0){\vector(1,0){10}}\put(0,0){\vector(0,1){10}}
  \put(2,3){\line(0,1){7}}\put(1.75,3){\line(1,0){0.5}}
  \put(4,4){\line(1,0){6}}\put(4,3.75){\line(0,1){0.5}}
  \put(6,6){\line(0,1){4}}\put(5.75,6){\line(1,0){0.5}}
  \thinlines\put(2,2){\circle{0.25}}\put(5,2){\circle{0.25}}
            \multiput(4,6)(1,0){2}{\circle{0.25}}\put(4,8){\circle{0.25}}
  \multiput(1,0)(1,0){9}{\circle*{0.2}}\multiput(0,1)(0,1){9}{\circle*{0.2}}
  \put(10,0){\makebox(1.3,0){x}}\put(0,10){\makebox(0,1.5){y}}
 \end{picture}\end{center}
\caption{ Beispiel der Menge $A\cup P$ in einem Quadranten $Q_i$ }
\end{figure}
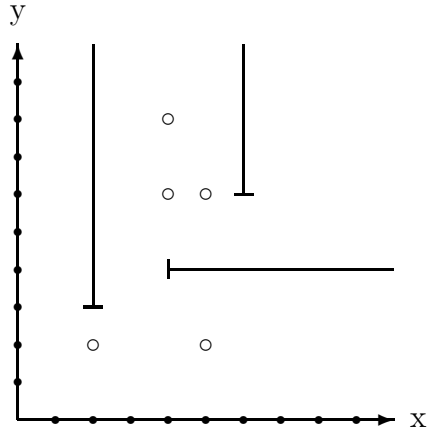
 $U'_{k_1+k_2}$ ist dann entsprechend die Gruppe derjenigen Elemente aus $G_n$, die
 au\3erhalb von $A\cup P$ die Identit"at sind. Durch eine leichte Modifikation
 der Funktion $f$ kann man die Aussage von Lemma 2.1 auf die Gruppen $\;U'_{k_1+k_2}\;$
 "ubertragen. Wir betrachten  
     \[ f\,':\,A\cup P\,\longrightarrow\,{\rm I}\!{\rm N}
                                 \times \{1,\ldots,k_1+k_2\}  \]
    \begin{eqnarray*}
      f\,'(p_j)&:\,=&(j,1)\quad\mbox{f"ur}\:j=1,\ldots,r\\
      f\,'(\,(x_1,i_1),m_1+y-1)&:\,=&(r+y,1)\quad\mbox{f"ur alle}\:y\in {\rm I}\!{\rm N}\\
       f\,'(\,(x_{\nu},i_{\nu}),m_{\nu}+y-1)&:\,=&(y,\nu)\quad\mbox{f"ur alle} \:y\in {\rm I}\!{\rm N},\,\nu=2,\ldots k_1\\
       f\,'(\,(y_{\mu},j_{\mu}),l_{\mu}+x-1)&:\,=&(x,k_1+\mu)\quad\mbox{f"ur alle}\:x\in {\rm I}\!{\rm N},\,\mu=1,\ldots k_2\:.
   \end{eqnarray*}
 Die endlich vielen Elemente $p_1,\ldots,p_r$ aus $P$ werden also lediglich der Menge
$((x_1,i_1),\ge m_1)$ vorgeschaltet. Dies ver"andert nichts an den f"ur die 
Gruppen $U_{k_1+k_2}$ gef"uhrten "Uberlegungen, d.h. Lemma 2.1 gilt ebenfalls f"ur die
 Gruppen $U'_{k_1+k_2}$.\newpage

\section{Konstruktion eines $G_n$-Komplexes}

 Ziel dieses Abschnittes ist es, einen $G_n$-Komplex zu konstruieren, der die Voraussetzungen
 von Korollar 1.3  erf"ullt. Zur Erinnerung: Wir ben"otigen einen
 $CW$-Komplex, auf dem $G_n$ in geeigneter Weise (Permutation der Zellen) 
 operiert, soda\3 der Stabilisator jeder Zelle vom Typ $(FP)_n$ und f"ur $n\ge 3$
 zus"atzlich endlich pr"asentiert ist. Dar"uberhinaus mu\3 die Zusammenziehbarkeit
 dieses Komplexes gew"ahrleistet sein.\\\\
Die einzelnen Schritte, die bei der Konstruktion eines solchen $G_n$-Komplexes
 durchlaufen werden, decken sich im wesentlichen mit denen, die im Falle der 
 Houghton-Gruppen ausgef"uhrt wurden. Auch hier ist der Ausgangspunkt das Monoid
 aller injektiven Abbildungen von $S$ nach $S$, die die Bedingungen B$_1$ und B$_2$ erf"ullen. Mit
 Hilfe eines Untermonoides wird auf der dem Monoid zugrundeliegenden Menge 
 eine Ordnung erkl"art, wodurch  ein Simplizialkomplex entsteht, auf dem 
$G_n$ simplizial operiert (vgl. Abschnitt 1.3, Schritt 1).\\
 Probleme ergeben sich dabei hinsichtlich der Stabilisatoren; es zeigt sich,
 da\3  nicht alle die oben genannten Eigenschaften besitzen. Die Situation l"a\3t
 sich jedoch retten, indem letztendlich nur ein gewisser Unterkomplex des erhaltenen 
 Simplizialkomplexes zur Fortf"uhrung des Beweises benutzt wird.

\subsection{Die partiell geordnete Menge $\cal M$}

Es sei $M$ das Monoid aller injektiven Abbildungen von $S$ nach $S$, die die in Abschnitt 2.1 
festgelegten Bedingungen B$_1$ und B$_2$ 
erf"ullen. F"ur ein $\alpha$ aus $M$, das die Elemente im $i$-ten Quadranten
 ab der Stelle $p_0$ um $m_i$ verschiebt, ist $\sum_{i=1}^n m_i\ge 0$.\\

 Ein Element $t$ aus $M$ hei\3t {\em Translation\/}, falls es ein 
 n-Tupel $(m_1,\ldots,m_n)$ aus ${\sl Z}\!\!{\sl Z}^n$ gibt, soda\3 
 \[(\,(x,y),i)\,t=(\,(x+m_i,y+m_i),i)\quad\mbox{f"ur {\em alle}}\quad(x,y)\in X\times Y\]
erf"ullt ist. Dabei mu\3 notwendigerweise $m_i\ge 0$ f"ur jedes $i\in\{1,\ldots,n\}$ gelten.\\
Mit $T$ bezeichnen wir das kommutative Untermonoid von $M$,
welches aus allen Translationen besteht. $T$ wird erzeugt von den Elementen
$t_1,\ldots,t_n$,
wobei
\[ \begin{array}{lcccl} 
      ((x,y),j)\,t_i=((x,y),j) & &\mbox{ f"ur alle $(x,y)\in (X\times Y)$,}
                             & & j\in\{1,\ldots,n\}-\{i\}\\
      ((x,y),i)\,t_i=((x+1,y+1),i) & & \mbox{ f"ur alle $(x,y)\in (X\times Y)$} & & \mbox{sonst}\:.
  \end{array}\]
$t_i$ ist also au\3erhalb des Quadranten $Q_i$ die Identit"at, und innerhalb von $Q_i$ verschiebt es s"amtliche Elemente 
jeweils um 1 in Richtung der x- und in Richtung der y-Achse (s. Abbildung 9). 
Mit Hilfe von $T$ wird nun eine Ordnung auf $M$ erkl"art. \\\\
{\bf Definition:}\quad F"ur $\alpha$ und $\beta$ aus $M$ ist $\alpha\le\beta$,
 falls es ein $t$ aus $T$ gibt mit $t\alpha=\beta$.\\
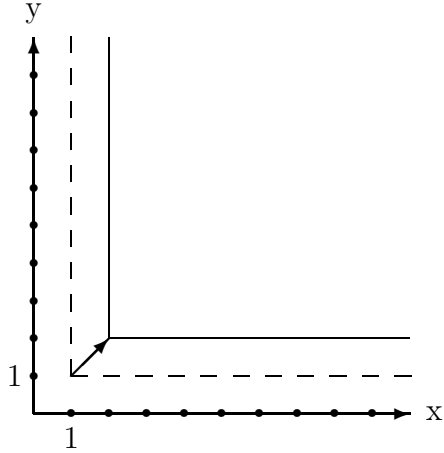
\begin{figure}[ht]
\unitlength 0.5cm
\begin{center}
\begin{picture}(12,11.5)
 \put(1,0.5){\begin{picture}(10,10)
    \thicklines\put(0,0){\vector(1,0){10}}
     \put(0,0){\vector(0,1){10}}
     \thinlines\multiput(1,1)(0.86,0){11}{\line(1,0){0.43}}
     \multiput(1,1)(0,0.86){11}{\line(0,1){0.43}}
     \multiput(1,0)(1,0){9}{\circle*{0.2}}
     \multiput(0,1)(0,1){9}{\circle*{0.2}}
     \put(1,-1.3){\makebox(0,1.3){1}}
     \put(-1,1){\makebox(1,0){1}}
     \put(10,0){\makebox(1.3,0){x}}
     \put(0,10){\makebox(0,1.3){y}} 
     \put(1,1){\thicklines\vector(1,1){1}}
  \put(2,2){\line(1,0){8}}
  \put(2,2){\line(0,1){8}}
\end{picture}}\end{picture}
  \end{center}
    \small Die gestrichelten Linien kennzeichnen den Bereich, welcher
            verschoben wird; der Pfeil gibt die Translation an, und das
            schraffierte Gebiet kennzeichnet den zugeh"origen Bildbereich.
\caption{Die Translation $t_i$ in $Q_i$}
\end{figure}
 \vfill
 Reflexivit"at und Transitivit"at dieser Relation sind unmittelbar ersichtlich.
Gilt f"ur $\alpha$,$\beta$ aus $M$ sowohl $\alpha\le \beta$ als auch $\beta\le \alpha$,
also $t\alpha=\beta$ und $t'\beta=\alpha$ mit $t,t'\in T$, so folgt $tt'\beta=\beta$.
 Daraus ergibt sich wegen der Injektivit"at von $\beta$: $tt'=$\,id$_S$. Dies ist nur m"oglich,
 wenn $t=t'=$\,id$_S$ (bei Translationen ist jedes $m_i\ge 0$), also $\alpha=\beta$.\\
Mit $\cal M$ bezeichnen wir im folgenden
die dem Monoid zugrundeliegende Menge, ausgestattet mit dieser Ordnung.
$T$ besitzt nun eine besondere Eigenschaft, die die Gerichtetheit von $\cal M$ gew"ahrleistet.
\vfill
\begin{lemma}
  Die Menge $T$ ist eine cofinale Teilmenge von $\cal M$, d.h. f"ur alle $\alpha$
  aus $\cal M$ existiert ein $t\in T$ mit $\alpha\le t$.
 \end{lemma}
\vfill
{\sc Beweis}:$\;$ Ein $\alpha$ aus $\cal M$ erf"ullt B$_1$, d.h. es gibt ein $p_0=(x_0,y_0)\in X\times Y$
 und ein $n$-Tupel $(m_1,\ldots,m_n)\in{\sl Z}\!\!{\sl Z}^n$ mit
    \[ ((x,y),i)\,\alpha=((x,y)+(m_i,m_i),i)\quad\mbox{f"ur alle}\quad (x,y)\ge p_0  \]
 und $i=1,\ldots,n$. Also entspricht f"ur alle $((x,y),i)\in S$ mit $(x,y)\ge p_0$ die Abbildung $\alpha$ der 
Translation $t_1{}^{m_1}t_2{}^{m_2}\cdot\,\ldots\,\cdot t_n{}^{m_n}$. Wir 
 konstruieren nun eine Translation $t_{\alpha}$, die die Elemente jedes
 Quadrnten soweit verschiebt, da\3 die Bildelemente $\ge p_0$ sind. Die
 Zusammensetzung $t_{\alpha}\alpha$ ist dann die gesuchte Translation $t$.\\
 Sei $l:\,=max\{x_0,y_0\}$.
Wir definieren $t_{\alpha}:\,=t_1{}^{l-1}\cdot\,\ldots\,\cdot t_n{}^{l-1}$
(s. Abbildung 10).
\pagebreak
 
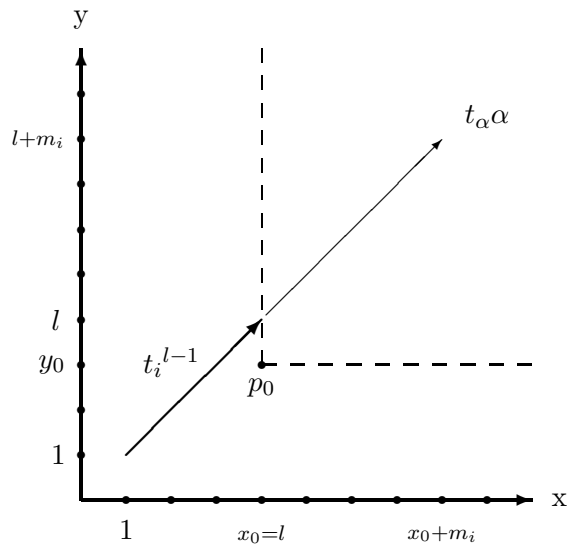
\begin{figure}[ht]\unitlength0.6cm
\small Die gestrichelte Linie begrenzt den Bereich, in dem $\scriptstyle
 \alpha$ eine Translation ist (die durch den schmalen Pfeil angegeben wird).
 
\begin{center}
\begin{picture}(12,11.2)
  \put(1,0.5){\begin{picture}(10,10)
                             \thicklines
                             \put(0,0){\vector(1,0){10}}
                             \put(0,0){\vector(0,1){10}}
                             \thinlines
                             \multiput(1,0)(1,0){9}{\circle*{0.2}}
                             \multiput(0,1)(0,1){9}{\circle*{0.2}}
                             \put(10,0){\makebox(1.2,0){x}}
                             \put(0,10){\makebox(0,1.3){y}}
                             \put(1,-1.3){\makebox(0,1.3){1}}
                             \put(-1,1){\makebox(1,0){1}}
                             \put(1,1){\thicklines\vector(1,1){3}}
   \put(1,2.3){\makebox(2,1.5){$t_i{}^{l-1}$}}
   \put(4.1,4.1){\vector(1,1){3.9}}
   \multiput(4,3)(0.63,0){10}{\line(1,0){0.315}}
   \multiput(4,3)(0,0.61){12}{\line(0,1){0.305}}
   \put(4,2){\makebox(0,1){$p_0$}}\put(4,3){\circle*{0.2}}
   \put(8,8){\makebox(2,1){$t_{\alpha}\alpha$}}
   \put(4,-1.5){\makebox(0,1.5){$\scriptstyle x_0=l$}}
   \put(8,-1.5){\makebox(0,1.5){$\scriptstyle x_0+m_i$}}
   \put(-1.3,3){\makebox(1.3,0){$y_0$}}
   \put(-1.1,4){\makebox(1.1,0){$l$}}
   \put(-2.1,8){\makebox(2.3,0){$\scriptstyle l+m_i$}}
                            \end{picture}}
 \end{picture}\end{center}
\caption{ Konstruktion von $t_{\alpha}$ in $Q_i$ mit $t_{\alpha}\alpha\in T$} 
\end{figure}

 $t_{\alpha}$ geh"ort zu $T$ und es gilt
 f"ur alle $i=1,\ldots,n$ und $(x,y)\in X\times Y$:
  \begin{eqnarray*}
    ((x,y),i)\,t_{\alpha}&=&((x,y),i)\,t_i{}^{l-1}\\
                          &=&((x+l-1,y+l-1),i)\\
                          &\ge & (l,l)\\
                          &\ge & (x_0,y_0)=p_0\:,
 \end{eqnarray*}
also ist $t_{\alpha}\alpha=t_{\alpha}(t_1{}^{m_1}\cdot\,\ldots\,\cdot t_n{}^{m_n})\in T$.
\hfill $\Box$\\

\begin{cor}
  $\cal M$ ist eine gerichtete Menge. (D.h. f"ur alle $\alpha,\beta\in\cal M$ existiert
   ein $\gamma\in\cal M$ mit $\alpha\le\gamma$ und $\beta\le\gamma$).
\end{cor}

{\sc Beweis}:$\;$ F"ur $\alpha,\beta$ aus $\cal M$ existieren nach dem letzten
 Lemma Elemente $t_1$ und $t_2$ aus $T$ mit $\alpha\le t_1$, $\beta\le t_2$.
 Daraus ergibt sich
     \begin{eqnarray*}  \alpha\le t_1\le t_2t_1& &\:\mbox{und}\\
                        \beta\le t_2 \le t_1t_2&,&
     \end{eqnarray*}
 wobei wegen der Kommutativit"at von $T$ $t_2t_1=t_1t_2$ ist. \hfill $\Box$ \\

\subsection{Untersuchung der Operation von $G_n$ auf $\protect{\mid\cal M\mid}$}

Wie bei den Houghton-Gruppen (vgl. Abschnitt 1.3, Schritt 1, (ii)) liefert
die gerichtete Menge $\cal M$ einen zusammenziehbaren $G_n$-Komplex $\mid\cal M\mid$,
 wobei die Operation von $G_n$ auf $\mid\cal M\mid$ durch
     \[ (\alpha_0,\ldots\alpha_p)\,g=(\alpha_0\,g,\ldots,\alpha_p\,g) \]
f"ur $(\alpha_0,\ldots,\alpha_p)\in \left|{\cal M}\right|$ und $g\in G_n$ gegeben ist.
Hinsichtlich der Stabilisatoren dieser Operation von $G_n$ ist die Situation
jedoch nicht so g"unstig wie im Falle der Houghton-Gruppen, wo diese alle
endlich sind. Wie wir im folgenden Abschnitt sehen werden, treten hier etliche
Stabilisatoren auf, die nur vom Typ {\sc FP}$_m$ f"ur ein $m<n$ sind.\\

Allerdings bilden  -- wie in 3.3 gezeigt wird --
 diejenigen Simplizes aus $\mid\cal M\mid$, deren Stabilisatoren
 die geforderten Endlichkeitseigenschaften aus Korollar 1.3 erf"ullen,
 einen $G_n$-invarianten und zusammenziehbaren Unterkomplex von $\mid\cal M\mid$, soda\3 es m"oglich ist,
 mit diesem den Beweis fortzusetzen.

\subsubsection{Eigenschaften der Stabilisatoren}

Bevor wir beginnen, uns mit dem Stabilisator eines Simplexes aus $\mid\cal M\mid$
auseinanderzusetzen, zeigen wir, da\3 es ausreicht, die Untersuchung auf Elemente
aus $\cal M$ zu beschr"anken.

\begin{lemma} F"ur ein p-Simplex $(\alpha_0,\ldots,\alpha_p)$ aus $\mid\cal M\mid$ ist
                \[ \mbox{\rm Stab}_{G_n}\left((\alpha_0,\ldots,\alpha_p)\right)
                    =\mbox{\rm Stab}_{G_n}(\alpha_0).\qquad\qquad\quad  \]
\end{lemma}

{\sc Beweis}:$\:$ Aus $(\alpha_0g,\ldots,\alpha_pg)=(\alpha_0,\ldots,\alpha_p)$
    folgt $\alpha_0\,g=\alpha_0$, also Stab$_{G_n}\left((\alpha_0,\ldots,\alpha_p)\right)
    \subseteq\mbox{\rm Stab}_{G_n}(\alpha_0)$.
    Wegen $\alpha_i>\alpha_0$, also $d_i\alpha_0=\alpha_i$ mit $d_i$ aus $T$
    f"ur alle $i=1,\ldots,p$, ist mit $\alpha_0 g=\alpha_0$ auch 
    $\alpha_i g=(d_i\alpha_0)g=d_i(\alpha_0 g)=d_i\alpha_0=\alpha_i$, woraus
    sich die umgekehrte Inklusion ergibt.\hfill $\Box$\\

  Sei nun $\alpha$ aus $\cal M$. Dann besteht Stab$_{G_n}(\alpha)$ genau aus
 denjenigen $g\in G_n$, deren Einschr"ankung auf die Menge $S\alpha$ die 
 Identit"at ist. Der Stabilisator von $\alpha$ wird also durch die Menge
 $S-S\alpha$ charakterisiert. Aus diesem Grund ist es notwendig,
 Aufschlu\3 "uber die Gestalt
 der Menge $S-S\alpha$ zu gewinnen. Mit den in Abschnitt 2.1 definierten Mengen
 \begin{eqnarray*}
     ((x,i),\ge m)&=&\{((x,y),i)\in S\mid y\ge m\}\quad\mbox{sowie}\\
     ((y,j),\ge\,l)&=&\{((x,y),j)\in S\mid x\ge l\}
 \end{eqnarray*}
 f"ur $(x,i)$ aus $X\times \{1,\ldots,n\}$, $(y,i)$ aus $Y\times \{1,\ldots,n\}$
 und $m,l$ aus ${\rm I}\!{\rm N}$ erhalten wir folgende Aussage.

 \begin{lemma}
    $\;$Sei $\alpha$ aus $\cal M$. Dann gibt es jeweils $k$ ($k\ge 0$)  viele paarweise verschiedene
    Elemente $(x_1,i_1),\ldots,(x_{k},i_{k})$ aus $X\times \{1,\ldots,n\}$ bzw.
    $(y_1,j_1),\ldots,(y_{k},j_{k})$ aus \mbox{$Y\times \{1,\ldots,n\}$}, sowie nat"urliche
    Zahlen $n_1,\ldots,n_{k}$ und $l_1,\ldots,l_{k}$, soda\3 $S-S\alpha$
    die disjunkte Vereinigung der Mengen
    \begin{eqnarray*}
       ((x_{\nu},i_{\nu}),\ge n_{\nu})\quad & & 1\le\nu\le k, \\
       ((y_{\mu},j_{\mu}),\ge\,l_{\mu})\quad & & 1\le\mu\le k\quad\mbox{und}\\
                    P\qquad\qquad                   & &                     
    \end{eqnarray*}
    ist, wobei $P$ eine endliche Teilmenge von $S$ ist. Dabei sind die Elemente
    $(x_1,i_1)\ldots,(x_{k},i_{k})$ bzw. $(y_1,j_1),\ldots,(y_{k},j_{k})$
     im Falle $k>0$ eindeutig bestimmt.
  \end{lemma}
   \vfill
   Abbildung 11 zeigt ein Beispiel einer solchen Menge f"ur $n=2$, also
   $S=Q_1\cup Q_2$.\\
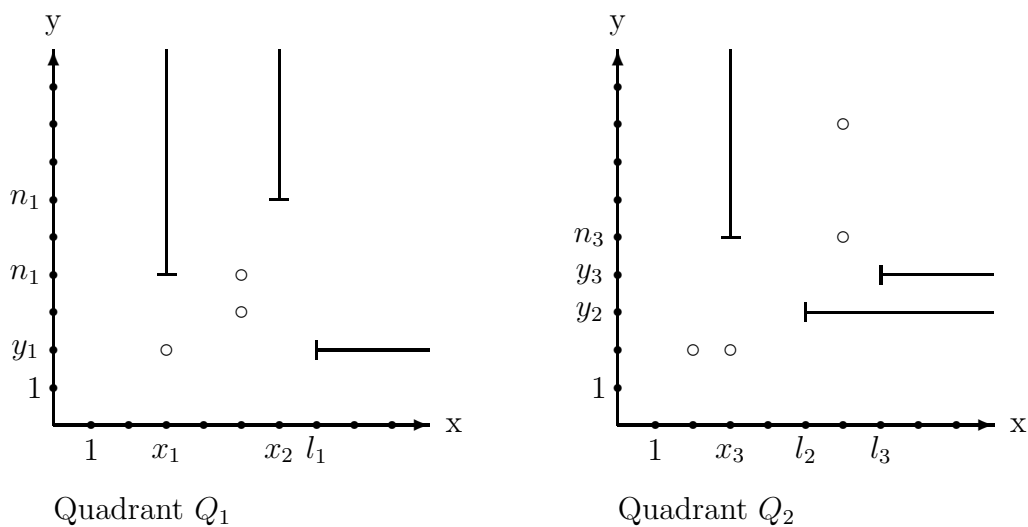
\begin{figure}[ht]
\unitlength0.5cm
   \begin{center}
   \begin{picture}(25,12.7)
     \multiput(0,2.5)(15,0){2}{\begin{picture}(10,10)
               \thicklines
               \put(0,0){\vector(1,0){10}}
               \put(0,0){\vector(0,1){10}}
               \thinlines
               \multiput(1,0)(1,0){9}{\circle*{0.2}}
               \multiput(0,1)(0,1){9}{\circle*{0.2}}
               \put(10,0){\makebox(1.3,0){x}}
               \put(0,10){\makebox(0,1.3){y}}
               \put(1,-1.3){\makebox(0,1.3){1}}
               \put(-1,1){\makebox(1,0){1}}
                              \end{picture}}
   \put(0,2.5){\begin{picture}(10,10)  
       \put(3,-1.5){\makebox(0,1.5){$x_1$}}
       \put(6,-1.5){\makebox(0,1.5){$x_2$}}
       \put(7,-1.3){\makebox(0,1.3){$l_1$}} 
       \put(-1.5,2){\makebox(1.5,0){$y_1$}}
       \put(-1.5,4){\makebox(1.5,0){$n_1$}}
       \put(-1.5,6){\makebox(1.5,0){$n_1$}}
       \put(3,2){\circle{0.3}} 
       \put(5,3){\circle{0.3}} 
       \put(5,4){\circle{0.3}} 
       \thicklines
       \put(3,4){\line(0,1){6}}\put(2.75,4){\line(1,0){0.5}}
       \put(6,6){\line(0,1){4}}\put(5.75,6){\line(1,0){0.5}}
       \put(7,2){\line(1,0){3}}\put(7,1.75){\line(0,1){0.5}}
             \end{picture}}
     \put(15,2.5){\begin{picture}(10,10)
                  \thicklines
                  \put(3,5){\line(0,1){5}}\put(2.75,5){\line(1,0){0.5}}
                  \put(5,3){\line(1,0){5}}\put(5,2.75){\line(0,1){0.5}}
                  \put(7,4){\line(1,0){3}}\put(7,3.75){\line(0,1){0.5}}
                  \put(-1.5,5){\makebox(1.5,0){$n_3$}}
                  \put(-1.5,4){\makebox(1.5,0){$y_3$}}
                  \put(-1.5,3){\makebox(1.5,0){$y_2$}}
                  \put(3,-1.5){\makebox(0,1.5){$x_3$}} 
                  \put(5,-1.3){\makebox(0,1.3){$l_2$}} 
                  \put(7,-1.3){\makebox(0,1.3){$l_3$}}\thinlines
                  \put(2,2){\circle{0.3}} 
                  \put(3,2){\circle{0.3}} 
                  \put(6,5){\circle{0.3}} 
                  \put(6,8){\circle{0.3}} 
                 \end{picture}}
     \put(0,0){Quadrant $Q_1$}
     \put(15,0){Quadrant $Q_2$}
     \end{picture}
     \end{center}
   \parbox{15.3cm}{\vspace*{0.5ex}\small Die fettgedruckten Linien entsprechen den Mengen
                     $\scriptstyle ((x_1,1),\ge n_1),((x_2,1),\ge n_2),
                     ((x_3,2),\ge n_3)$ und $\scriptstyle ((y_1,1),\ge l_1),
                     ((y_2,2),\ge l_2),((y_3,2),\ge l_3)$ dar. Die Kreise 
                     stellen die Menge $\scriptstyle P$ \\}
   \caption{ Die Menge $S-S\alpha$ f"ur ein $\alpha$ aus $\cal M$ }
   \end{figure}
   \vfill  
  {\sc Beweis}:$\:$ Im folgenden bezeichnen wir zu einer Teilmenge $A$ von $S$
     und einem $(x,y)$ aus $X\times Y$ mit $A_{<(x,y)}$ die Menge aller Elemente
     $((x'\!,y'),i)$ aus $A$ mit $(x'\!,y')<(x,y)$. Analog definiert man die
     Mengen $A_{\le (x,y)}$, $A_{>(x,y)}$ und $A_{\ge (x,y)}$.\\
     \vfill
     F"ur ein $\alpha$ aus $\cal M$ gibt es ein $p_0=(x_0,y_0)\in
    (X\times Y)$ und ein $n$-Tupel $(m_1,\ldots,m_n)\in{\sl Z}\!\!{\sl Z}^n$,
    soda\3 $((x,y),i)\alpha=((x,y)+(m_i,m_i),i)$ f"ur
    alle $(x,y)\ge p_0$ und $i=1,\ldots,n$ ist (Bedingung B$_1$).
    Es sei nun $x_i:\,=x_0+m_i$ und $y_i:\,=y_0+m_i$, $i=1,\ldots,n$.
    Dann ist (s. Abbildung 12)

 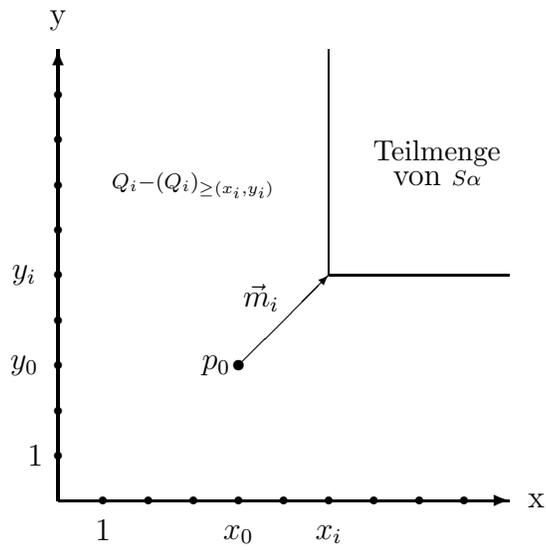
\begin{figure}[ht]
 \unitlength0.6cm
 \begin{center}
 \begin{picture}(10,11.5)
   \put(0,0.5){\begin{picture}(10,10)
            \thicklines
            \put(0,0){\vector(1,0){10}}
            \put(0,0){\vector(0,1){10}}
            \thinlines
            \multiput(1,0)(1,0){9}{\circle*{0.2}}
            \multiput(0,1)(0,1){9}{\circle*{0.2}}
            \put(10,0){\makebox(1.2,0){x}}
            \put(0,10){\makebox(0,1.3){y}}
            \put(1,-1.3){\makebox(0,1.3){1}}
            \put(-1,1){\makebox(1,0){1}}
            \put(4,-1.5){\makebox(0,1.5){$x_0$}}
            \put(6,-1.5){\makebox(0,1.5){$x_i$}}
            \put(-1.5,3){\makebox(1.5,0){$y_0$}}
            \put(-1.5,5){\makebox(1.5,0){$y_i$}}
            \put(3,3){\makebox(1,0){$p_0$}}
            \put(4,3){\vector(1,1){2}}
            \put(6,5){\line(1,0){4}}
            \put(6,5){\line(0,1){5}}
            \put(3.5,4){\makebox(2,1){$\vec{m}_i$}}
            \put(0,4.5){\makebox(6,5){$\scriptstyle Q_i-(Q_i)_{\ge (x_i,y_i)}$}}  
            \put(6.4,5){\makebox(4,5){${\mbox{\small Teilmenge}\atop\mbox{von }\scriptstyle S\alpha}$}}
            \put(4,3){\circle*{0.25}}
            \end{picture}}
    \end{picture}\end{center}
   \caption{Das schraffierte Gebiet entspricht der Menge $(Q_i)_{\ge (x_i,y_i)}$}
    \end{figure}
  \clearpage
  \[ (Q_i)_{\ge (x_i,y_i)}=(Q_i)_{\ge p_0}\alpha\subseteq 
         S\,\alpha\quad\mbox{bzw.}\quad
         \bigcup_{i=1}^n(Q_i)_{\ge (x_i,y_i)}
         \subseteq  S\,\alpha\;,\;\mbox{d.h.} \]
  \begin{eqnarray}
        S-S\alpha &=& (\bigcup_{i=1}^n Q_i)-S\alpha  \\
                  &\subseteq& (\bigcup_{i=1}^n Q_i)-(\bigcup_{i=1}^n (Q_i)_{\ge (x_i,y_i)})
                  = \bigcup_{i=1}^n (Q_i-(Q_i)_{\ge (x_i,y_i)})\:.\nonumber
  \end{eqnarray}
 Die Menge $Q_i-(Q_i)_{\ge(x_i,y_i)}$ l"a\3t sich wiederum in folgende zueinander
 disjunkte Mengen aufteilen (s. Abbildung 13):\\
 \parbox{1cm}{\begin{eqnarray}\end{eqnarray}}
 \hfill\parbox{14cm}
   {\begin{eqnarray*}
      ((\xi,i),\ge y_i) & & \quad\mbox{mit}\quad \xi\in X,1\le \xi<x_i\:, \\
      ((\eta,i),\ge x_i) & & \quad\mbox{mit}\quad \eta\in Y,1\le \eta<y_i
                              \quad\mbox{und}\\
      (Q_i)_{<(x_i,y_i)}& &.
   \end{eqnarray*}}
 
  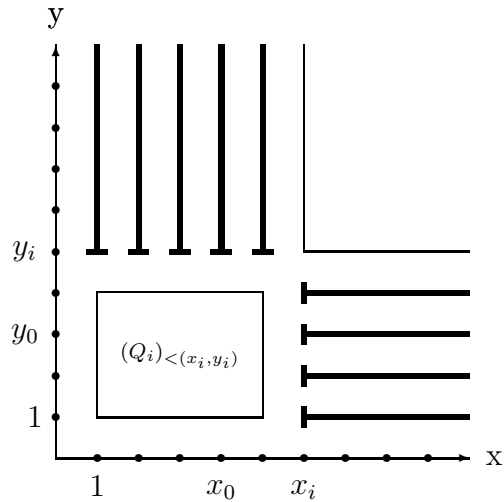
\begin{figure}[ht]
  \parbox{15.3cm}
             {\small 
              Die fettgedruckten vertikalen Linien repr"asentieren die Mengen der Form
              $\scriptstyle((\xi ,i),\ge y_i)$ f"ur $\scriptstyle
              1\le \xi\le x_i$. Entsprechend
              stellen die horizontalen
              fettgedruckten Linien die Mengen $\scriptstyle ((\eta ,i),\ge x_i)\;,1\le \eta\le y_i$ dar.}
  \unitlength0.55cm
  \begin{center}
  \begin{picture}(10,12)
 \put(0,1){\begin{picture}(10,10)
           \put(0,0){\vector(1,0){10}}
           \put(0,0){\vector(0,1){10}}
           \multiput(1,0)(1,0){9}{\circle*{0.2}}   
           \multiput(0,1)(0,1){9}{\circle*{0.2}}
           \put(10,0){\makebox(1.2,0){x}}
           \put(0,10){\makebox(0,1.3){y}}
           \put(1,-1.3){\makebox(0,1.3){1}}
           \put(4,-1.5){\makebox(0,1.5){$x_0$}}
           \put(6,-1.5){\makebox(0,1.5){$x_i$}}
           \put(0,10){\makebox(0,1.3){y}}
           \put(-1,1){\makebox(1,0){1}}
           \put(-1.5,3){\makebox(1.5,0){$y_0$}}
           \put(-1.5,5){\makebox(1.5,0){$y_i$}}
           \put(1,1){\framebox(4,3){$\scriptstyle (Q_i)_{<(x_i,y_i)}$}}
           \linethickness{0.65mm}
           \multiput(1,5)(1,0){5}{\line(0,1){5}}
           \multiput(0.75,5)(1,0){5}{\line(1,0){0.5}}
           \multiput(6,1)(0,1){4}{\line(1,0){4}}
           \multiput(6,0.75)(0,1){4}{\line(0,1){0.5}}
           \thinlines
           \put(6,5){\line(1,0){4}}
           \put(6,5){\line(0,1){5}}
           \end{picture}}
   \end{picture}\end{center}
  \caption{Aufteilung der Menge $Q_i-(Q_i)_{\ge (x_i,y_i)}$}
\end{figure}

  Wegen (3.1) befindet sich die Menge $S-S\alpha$ in den Teilmengen $Q_i-(Q_i)_{\ge (x_i,y_i)}$ 
 von $S$. Es gen"ugt also, den Schnitt von $S-S\alpha$ mit den einzelnen Mengen aus
 (3.2) zu untersuchen; daraus setzt sich dann die gesamte Menge $S-S\alpha$
 zusammen.\\

 Wir betrachten zuerst die Mengen der Form  $((\xi,i),\ge y_i)$, $1\le\xi<x_i$. Nach
 B$_2$(a) hat f"ur alle $(x,j)\in X\times \{1,\ldots,n\}$ die Bildmenge von
 $((x,j),\ge y_0)$ unter $\alpha$ die Gestalt $((x,j)^{'},\ge y_0+q)$
 f"ur ein $q\in{\sl Z}\!\!{\sl Z}$. Es gibt nun 2 F"alle:
 \begin{itemize}
   \item[1.] $\:$ Es existiert ein $(x,j)\in X\times\{1,\ldots,n\}$ mit $(x,j)^{'}=
             (\xi,i)$
   \item[2.] $\:$ F"ur alle $(x,j)\in X\times\{1,\ldots,n\}$ ist $(x,j)^{'}\not=(\xi,i)$.
 \end{itemize}
 Im 1. Fall ist $((\xi,i),\ge y_0+q)=((x,j),\ge y_0)\,\alpha\subseteq S\alpha$,
   d.h der Anteil von $S-S\alpha$ in $((\xi,i),\ge y_i)$ mu\3 aus der Menge
  $((\xi,i),\ge y_i)-((\xi,i),\ge y_0+q)$ stammen (s. Abbildung 14).
  Somit ist der Schnitt $((\xi,i),\ge y_i)\cap(S-S\alpha)$ endlich (bzw. f"ur
  $ y_0+q\le y_i$ leer).

 
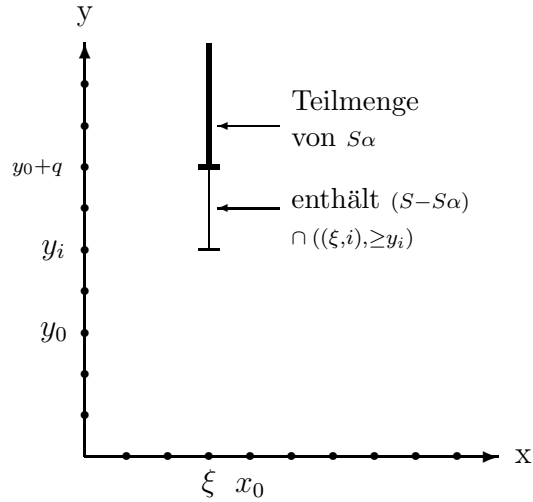
\begin{figure}[ht]
 \unitlength0.55cm
 \begin{center}
  \begin{picture}(10,11.5)
    \put(0,0.5){\begin{picture}(10,10)
           \thicklines
           \put(0,0){\vector(1,0){10}}
           \put(0,0){\vector(0,1){10}}
           \thinlines
           \multiput(1,0)(1,0){9}{\circle*{0.2}}   
           \multiput(0,1)(0,1){9}{\circle*{0.2}}
           \put(10,0){\makebox(1.2,0){x}}
           \put(0,10){\makebox(0,1.3){y}}
           \put(3,-1.3){\makebox(0,1.3){$\xi$}}         
           \put(-1.5,3){\makebox(1.5,0){$y_0$}}         
           \put(4,-1.5){\makebox(0,1.5){$x_0$}}
           \put(-1.5,5){\makebox(1.5,0){$y_i$}}
           \put(-2.3,7){\makebox(2.3,0){$\scriptstyle y_0+q$}}
           \put(3,5){\line(0,1){5}}\put(2.75,5){\line(1,0){0.5}}
           \linethickness{0.65mm}
           \put(3,7){\line(0,1){3}}\put(2.75,7){\line(1,0){0.5}}
           \thinlines
           \put(5,8){\parbox{2cm}{\small Teilmenge von $\scriptstyle S\alpha$}}
           \put(5,6){\parbox[t]{2.45cm}{\small enth"alt $\scriptstyle
                                      (S-S\alpha)\\\cap\:((\xi,i),\ge y_i)$}}
           \put(4.7,8){\vector(-1,0){1.5}}
           \put(4.7,6){\vector(-1,0){1.5}}
           \end{picture}}
   \end{picture}\end{center}
  \begin{center}
    \parbox[t]{15.3cm}{\small Die fette Linie entspricht der Menge
                     $\scriptstyle ((\xi,i),\ge y_0+q)$, die andere Linie
                     der Menge $\scriptstyle ((\xi,i),\ge y_i).$}
  \end{center}
  \caption{ 1. Fall}
  \end{figure}
 
  Im 2.Fall wird {\em kein\/} Element $((x,y),j)$ aus $S$ mit $y\ge y_0$
  unter $\alpha$
  auf \mbox{$((\xi,i),\ge y_i)$} abgebildet.
  Bleibt noch die Menge der Elemente $((x,y),j)$ mit $y<y_0$, die
  m"oglicherweise Bildelemente in $((\xi,i),\ge y_i)$ besitzt.\\ Diese Menge
  setzt sich aus den einzelnen $((\eta,j),\ge x_0)$ mit $1\le \eta <y_0$, 
  $j=1,\ldots,n$ sowie der Menge $S_{<(x_0,y_0)}$ zusammen.
  Da im Bildbereich dieser $((\eta,j),\ge x_0)$ nur y-Werte $<y_i$ vorkommen
 k"onnen (sonst w"are der Schnitt mit $(Q_i)_{\ge (x_i,y_i)}$ nicht leer),
 wird kein Element von $((\eta,j),\ge x_0)$ nach $((\xi,i),\ge y_i)$ abgebildet.
 Es bleiben daher h"ochstens die Elemente aus der endlichen Menge
 $S_{<(x_0,y_0)}$ mit m"oglichen Bildelementen in $((\xi,i),\ge x_0)$,
 d.h. $((\xi,i),\ge y_0)\cap S\alpha$ ist endlich.\\

  Es existiert daher ein $n_{(\xi,i)}\in{\rm I}\!{\rm N}$ mit
   \begin{eqnarray*}
       ((\xi,i),\ge y_i)\cap (S-S\alpha)&=& ((\xi,i),\ge y_i)-(\,((\xi,i),\ge y_i)\cap S\alpha)\;)\\
                                        &=& ((\xi,i),\ge n_{(\xi,i)})\cup E_{(\xi,i)}\:,
   \end{eqnarray*}
  wobei $E_{(\xi,i)}$ eine endliche Menge disjunkt zu $((\xi,i),\ge n_{(\xi,i)})$
 ist (s. Abbildung 15).
 \pagebreak
 
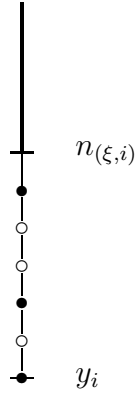
\begin{figure}[ht]
 \unitlength0.5cm
 \begin{center}
  \begin{picture}(25,11)
  \put(1.7,0){\line(1,0){0.6}}
  \multiput(2,0.2)(0,1){5}{\line(0,1){0.6}}
  \put(2,5){\line(0,1){1}}
  \put(2,6){\thicklines\line(0,1){4}}
  \put(2,0){\circle*{0.25}}
  \put(2,2){\circle*{0.25}}
  \put(2,5){\circle*{0.25}}
  \put(2,1){\circle{0.25}}
  \put(2,3){\circle{0.25}}
  \put(2,4){\circle{0.25}}
  \put(1.7,6){\thicklines\line(1,0){0.6}}
  \put(3,6){\makebox(2.5,0){$n_{(\xi,i)}$}}
  \put(3,0){\makebox(1.5,0){$y_i$}}
  \put(9,8){\parbox[t]{7cm}
            {\small \mbox{Nebenstehende Abbildung zeigt die Men-}\linebreak ge $\scriptstyle ((\xi,i),\ge y_i)$.
            Die ausgef"ullten Kreise stellen die endliche Menge
            $\scriptstyle ((\xi,i),\ge y_i)\cap S\alpha$ dar. Die unausgef"ullten
            Kreise entsprechen der Menge $\scriptstyle E_{(\xi,i)}$, und die
            fettgedruckte Linie repr"asentiert die Menge $\scriptstyle ((\xi,i),\ge n_{(\xi,i)})$.}}
 \end{picture}
 \end{center}
 \caption{ 2.Fall}
 \end{figure}

Analog zeigt man, da\3 f"ur die Mengen $((\eta,i),\ge x_i)$ mit $1\le \eta 
 <y_i$, $i=1,\ldots,n$ einer der beiden folgenden F"alle zutrifft:
 \begin{itemize} 
   \item[1.]$\:$ $((\eta,i),\ge x_i)\cap (S-S\alpha)$ ist endlich\quad oder
   \item[2.]$\:$ $((\eta,i),\ge x_i)=((\eta,i),\ge l_{(\eta,i)})\cup E^{'}_{(\eta,i)}$
 \end{itemize}
 mit $l_{(\eta,i)}\in{\rm I}\!{\rm N}$, $E^{'}_{(\eta,i)}$ endliche und zu 
 $((\eta,i),\ge l_{(\eta,i)})$ disjunkte Menge.
 Dabei tritt der 1. Fall ein,
 wenn es ein $(y,j)\in Y\times \{1,\ldots,n\}$ gibt mit der Eigenschaft
$((y,j),\ge x_0)\alpha
 =((\eta,i),\ge x_0+r)$, ansonsten der 2.Fall.\\
Vereinigt man nun s"amtliche endliche Teilmengen von $S$, die zu den 1. F"allen             
geh"oren, sowie s"amtliche Mengen $E_{(\xi,i)}$ bzw. $E^{'}_{(\eta,i)}$,die in den
2. F"allen auftreten, zu der Menge $P$, so erh"alt man die gew"unschte Darstellung
 der Menge $S-S\alpha$ durch\\
  \parbox{1cm}{\begin{eqnarray}\end{eqnarray}}\hfill
  \parbox{14cm}
    {\begin{eqnarray*}
      ((\xi,i),\ge n_{(\xi,i)})& &\quad\mbox{$(\xi,i)$ geh"ort zum 1. Fall}\:,\\
      ((\eta,i),\ge l_{(\eta,i)})& &\quad\mbox{$(\eta,i)$ geh"ort zum 1. Fall}\quad\mbox{und}\\
            P\:.    & &
     \end{eqnarray*}}
 Bleibt zu "uberlegen, wie gro\3 die Anzahl der $(\xi,i)$ (bzw. $(\eta,i)$) ist, 
 die zum 1. Fall geh"oren. Wir betrachten dazu die in 2.3 definierte Abbildung
 $ \pi\,:\,G_n\,\longrightarrow\,H_n\,$,
 die jedem $g\in G_n$ die induzierte Permutation auf der Menge $X\times\{1,\ldots,n\}$
 zuordnet.
$\pi$ l"a\3t sich fortsetzen auf $\cal M$, wobei $\pi(\alpha)$ 
 eine injektive Abbildung von $X\times\{1,\ldots,n\}$ in sich selbst ist, die
 die Elemente in $X\times\{i\}$ ab $x_0$ um $m_i$ verschiebt
 (dabei stammt $m_i$ von $\alpha$).\\
 Jedes Element $(x,j)$, das in $(X\times\{1,\ldots,n\}-(X\times\{1,\ldots,n\})\,\pi(\alpha)$
 liegt, geh"ort zu Fall 1. 
 Die Anzahl der Elemente in $(X\times\{1,\ldots,n\}-(X\times\{1,\ldots,n\})\,\pi(\alpha)$
 ist aber gerade $\sum_{i=1}^n m_i$. Analog zeigt man, da\3 es genau $\sum_{i=1}^n m_i$
 viele $(\eta,i)$ gibt, die zu Fall 1 geh"oren, d.h. mit $k:\,=\sum_{i=1}^n m_i$
 folgt die Behauptung. \hfill $\Box$\\ 

{\bf Definition:}$\;$ F"ur ein $\alpha$ aus $\cal M$ hei\3t
 die eindeutig bestimmte positive ganze Zahl $k$ aus Lemma 3.4 der
 {\em Grad\/} von $\alpha$.
 Dieser Grad wird im folgenden mit $gr(\alpha)$ bezeichnet.\\

{\bf Bemerkung:}$\;$Allgemein versteht man bez"uglich einer partiell geordneten
 Menge $M$ unter dem Grad eines Elementes $a$ aus $M$ das Supremum aller
 L"angen von Ketten aus $M$, die $a$ als gr"o\3tes Element besitzen (vgl. [Qu]).
 Wir zeigen in 3.2.2, das die aus
Lemma 3.4 stammende Zahl $k$ genau den Grad gem"a\3 dieser Auffassung
 angibt, soda\3 obige Definition gerechtfertigt ist.\\[0.8cm]
In Abschnitt 1.2 haben wir Untergruppen von $G_n$ untersucht, deren Elemente
 au\3erhalb einer Menge der in 
Lemma 3.4 beschriebenen Form die Identit"at sind. Im Hinblick auf Lemma 2.1 sowie Bemerkung 2.2
ergibt sich folgendes Resultat.

\begin{cor}
   F"ur ein $\alpha$ aus $\cal M$ mit $k:\,=gr(\alpha)\ge 1$ ist {\rm Stab}$_{G_n}(\alpha)$
   vom Typ {\sc (FP)}$_{2k-1}$ und nicht vom Typ {\sc (FP)}$_{2k}$.
   Im Fall $2k\ge 3$ ist er zus"atzlich endlich pr"asentiert.
\end{cor}

Das hei\3t f"ur Simplizes $(\alpha_0,\ldots,\alpha_p)$ aus $\left|\cal M\right|$ mit
$2\,gr(\alpha_0)\ge n+1$ ist der Stabilisator Stab$_{G_n}(\alpha_0,\ldots,\alpha_p)$ vom Typ {\sc (FP)}$_n$ und 
f"ur  $n\ge 3$ endlich pr"asentiert. Aus diesem Grunde ist es erforderlich,
 zur Konstruktion des Simplizialkomplexes nur
diejenigen $\alpha$ aus $\cal M$ zu benutzen mit $2\,gr(\alpha)\ge n+1$. 
 Wegen $2n\ge n+1$ f"ur alle $n$ aus ${\rm I}\!{\rm N}$, kann man sich auf
  Elemente beschr"anken, deren Grad mindestens $n$ betr"agt. Im Anschlu\3 
an den n"achsten Abschnitt zeigen wir, da\3 sich aus diesen Elementen ein Simplizialkomplex ergibt,
 der alle geforderten Vorraussetzungen erf"ullt.
Zuvor erl"autern wir einige Eigenschaften des Grades.

\subsubsection{Eigenschaften des Grades}

Wir zeigen in diesem Abschnitt den Zusammenhang zwischen dem in 3.2.1 definierten 
Grad eines Elementes $\alpha$ 
 und gewissen Ketten aus $\cal M$. Es stellt sich heraus, da\3 $gr(\alpha)$ die gr"o\3te ganze Zahl $k$ ist,
 soda\3 es eine Kette der Form $\:\alpha=\alpha_k>\alpha_{k-1}>\ldots >\alpha_0$ in $\cal M$ gibt.\\\\
Die maximale L"ange einer solchen Kette ist erst erreicht, wenn f"ur alle 
 $i=1,\ldots,k$ ein $t_{j_i}\in \{1,\ldots,n\}$ ($t_1,\ldots,t_n$
 sind die Erzeugenden von $T$) existiert mit
   \begin{equation}
          \alpha_i= t_{j_i}\alpha_{i-1}\:.
  \end{equation}
 (Ansonsten k"onnte durch Einschieben neuer
Elemente die gegebene Kette verl"angert werden, zum Beispiel l"a\3t sich
 aus $t_it_j\beta>\beta$ die Kette $t_it_j\beta>t_j\beta>\beta$ bilden).\\

 Wir untersuchen daher zun"achst das Verhalten des Grades hinsichtlich
 zweier Elemente aus $\cal M$, die sich nur durch Vorschalten einer Translation $t_i$
 voneinander unterscheiden.

\begin{lemma}
   F"ur $\alpha$ und $\beta$ aus $\cal M$
   mit $t_i\beta=\alpha$ ist $gr(\beta)=gr(\alpha)-1$.
\end{lemma}
\vspace*{0.2cm}
{\sc Beweis}:$\;$ Aus $t_i\beta=\alpha$ folgt $S(t_i\beta)=S\alpha$.
    Wir zeigen nun
      \begin{equation}
        S-S\alpha=(S-S\beta)\cup(S-St_i)\beta
      \end{equation}
  \begin{quote}
     Bew.:$\:$ Sei $s\in S$, $s\notin S\alpha=St_i\beta$. Da $St_i\subseteq
     S$, also $St_i\beta\subseteq S\beta$, ist dies "aquivalent zu: $s\notin
     S\beta$ (also $s\in S-S\beta$) oder $s\in S\beta-St_i\beta$.
     Wegen der Injektivit"at von $\beta$ ist $S\beta-St_i\beta=
     (S-St_i)\beta$, d.h. entweder ist $s\in S-S\beta$ oder 
     $s\in (S-St_i)\beta$.
  \end{quote}
 Nun ist $S-St_i=\{((1,y),i)\mid y\in{\rm I}\!{\rm N}\}\cup
                 \{((x,1),i)\mid x\in{\rm I}\!{\rm N}\}$  (s. Abbildung 16).

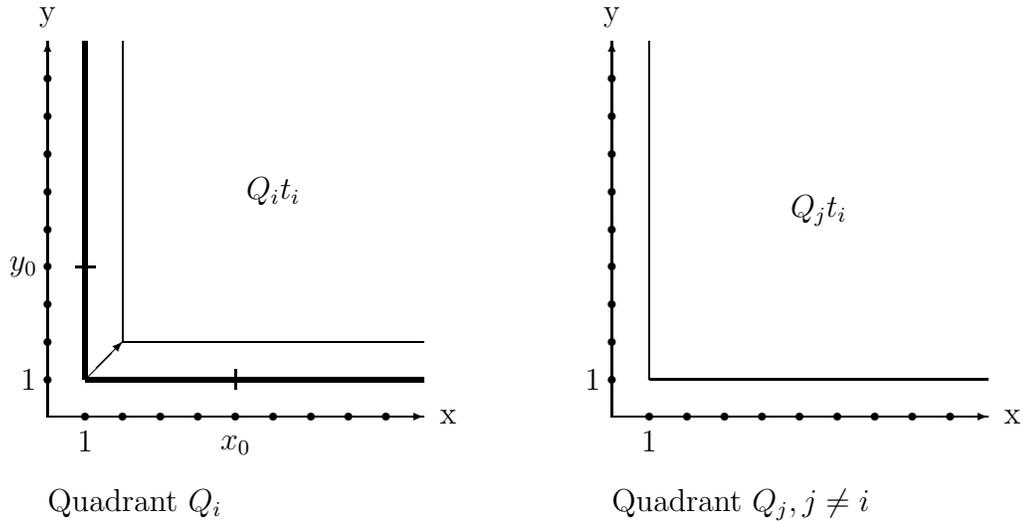
\begin{figure}[ht]
  \unitlength0.5cm
   \begin{center}
   \begin{picture}(25,13)
     \multiput(0,2)(15,0){2}{\begin{picture}(10,10)
               \put(0,0){\vector(1,0){10}}
               \put(0,0){\vector(0,1){10}}
               \multiput(1,0)(1,0){9}{\circle*{0.2}}
               \multiput(0,1)(0,1){9}{\circle*{0.2}}
               \put(10,0){\makebox(1.3,0){x}}
               \put(0,10){\makebox(0,1.3){y}}
               \put(1,-1.3){\makebox(0,1.3){1}}
               \put(-1,1){\makebox(1,0){1}}
                              \end{picture}}
     \put(0,2){\begin{picture}(10,10)
               \put(2,2){\line(1,0){8}}
               \put(2,2){\line(0,1){8}}
               \linethickness{0.65mm}
               \put(1,1){\line(0,1){9}}
               \put(1,1){\line(1,0){9}}
               \thicklines
               \put(0.75,4){\line(1,0){0.5}}
               \put(5,0.75){\line(0,1){0.5}}
               \thinlines
               \put(-1.3,4){\makebox(1.3,0){$y_0$}}
               \put(5,-1.5){\makebox(0,1.5){$x_0$}}
               \put(2,2){\makebox(8,8){$Q_it_i$}}
               \put(1,1){\vector(1,1){1}}
               \end{picture}}
    \put(15,2){\begin{picture}(10,10)
               \put(1,1){\line(1,0){9}}
               \put(1,1){\line(0,1){9}}
               \put(1,1){\makebox(9,9){$Q_jt_i$}}
               \end{picture}}
    \put(0,-0.5){Quadrant $Q_i$}
    \put(15,-0.5){Quadrant $Q_j,j\not=i$}
   \end{picture}  
  \end{center}
  \caption{Die fettgedruckten Linien entsprechen der Menge $S-St_i=Q_i-Q_it_i$}
  \end{figure}
 \vspace*{0.2cm}
 Die Bildmenge von $\{((x,1),i)\mid x\ge x_0\}$ unter $\beta$ ist von der Form
 (vgl. B$_2$(a)) 
            \[  ((x,1)',\ge m) \]
 f"ur ein $(x,1)'\in X\times\{1,\ldots,n\}$, $m\in{\rm I}\!{\rm N}$, und
 die Bildmenge von $\{((1,y),i)\mid y\ge y_0\}$ l"a\3t sich schreiben als
 (vgl. B$_2$(b))
            \[  ((y,1)',\ge l)  \]
 f"ur ein $(y,1)'\in Y\times \{1,\ldots,n\}$, $l\in {\rm I}\!{\rm N}$.
 Somit ist 
   \begin{equation}
         (S-St_i)\beta=((x,1)',\ge m)\cup ((y,1)',\ge l)\cup E\:,
   \end{equation}
 wobei $E:\,=\{((x,1),i)\mid 1\le x<x_0\}\beta
      \cup\{((1,y),i)\mid 1\le y<y_0\}\beta$
 eine endliche und zu $((x,1)',\ge m)\cup((y,1)',\ge l)$ disjunkte 
 Menge ist. Sei nun $gr(\beta)=k$ und 
   \begin{eqnarray}
      S-S\beta&=&((x_1,i_1),\ge n_1)\cup\ldots\cup((x_k,i_k),\ge n_k)\cup\\
              & & ((y_1,j_1),\ge l_1)\cup\ldots\cup((y_k,j_k),\ge l_k)\cup\nonumber\\
              & & P \nonumber
   \end{eqnarray}
  eine wie in Lemma 3.4 beschriebene Zerlegung der Menge $S-S\beta$.
  Nach (3.8) gilt $S-S\alpha=(S-S\beta)\cup(S-St_i)\beta$, wobei
  $(S-St_i)\beta\subseteq S\beta$, also disjunkt zu $S-S\beta$ ist. Schreiben
  wir $(S-St_i)\beta$ wie in (3.9) und $S-S\beta$ wie in (3.10), so erhalten
  wir $S-S\alpha$ als disjunkte Vereinigung der Mengen
      \begin{eqnarray*}
     ((x_1,i_1),\ge n_1)\cup\ldots\cup((x_k,i_k),\ge n_k)\cup((x,1)',\ge m)&\cup&\\
     ((y_1,j_1),\ge l_1)\cup\ldots\cup((y_k,j_k),\ge l_k)\cup((y,1)',\ge l)&\cup&\\
      P\cup E\:,
      \end{eqnarray*}
  d.h. $gr(\alpha)=k+1=gr(\beta)+1$.\hfill $\Box$\\
 
 Weiterhin ben"otigen wir

\begin{lemma}
    Sei $i\in \{1,\ldots,n\}$. Es existiert genau dann ein $\beta$ aus
    $\cal M$ mit $t_i\beta=\alpha$, wenn $gr(\alpha)>0$ ist.
 \end{lemma}

{\sc Beweis}:$\;$Die Notwendigkeit von $gr(\alpha)>0$ ergibt sich aus dem
   vorigen Lemma und der Eigenschaft $gr(\gamma)\ge 0$ f"ur alle $\gamma
   \in \cal M$.\\
   F"ur die umgekehrte Richtung setzen wir $k:\,=gr(\alpha)$, $k>0$. Dann
   existieren in $S-S\alpha$ mindestens zwei zueinander disjunkte Mengen
   der Form
    \begin{eqnarray*}
        ((y_1,j_1),\ge l_1)\quad\mbox{mit}& & (y_1,j_1)\in Y\times \{1,\ldots,n\}, l_1\in {\rm I}\!{\rm N}\quad\mbox{sowie}\\
        ((x_1,i_1),\ge n_1)\quad\mbox{mit}& & (x_1,i_1)\in X\times \{1,\ldots,n\}, n_1\in {\rm I}\!{\rm N}\:.
    \end{eqnarray*}
 Wir definieren nun $\beta$ auf der Menge $St_i$ als $t_i{}^{-1}\alpha\,$;
 und auf der Menge $S-St_i=\{((x,1),i)\mid x\in {\rm I}\!{\rm N}\}\cup
 \{((1,y),i)\mid y\in {\rm I}\!{\rm N}\}$ setzen wir (s. Abbildung 17)\\
   \parbox{1cm}{\begin{eqnarray}\end{eqnarray}}\hfill
   \parbox{14cm}{\begin{eqnarray*}
                   ((1,y),i)\,\beta:&=&((x_1,y+n_1-1),i_1)\quad\mbox{f"ur alle}\quad y\in{\rm I}\!{\rm N},\\
                   ((x,1),i)\,\beta:&=&((x+l_1-2,y_1),j_1)\quad\mbox{f"ur alle}\quad x\in{\rm I}\!{\rm N},
                     x\ge 2.
                   \end{eqnarray*}}

  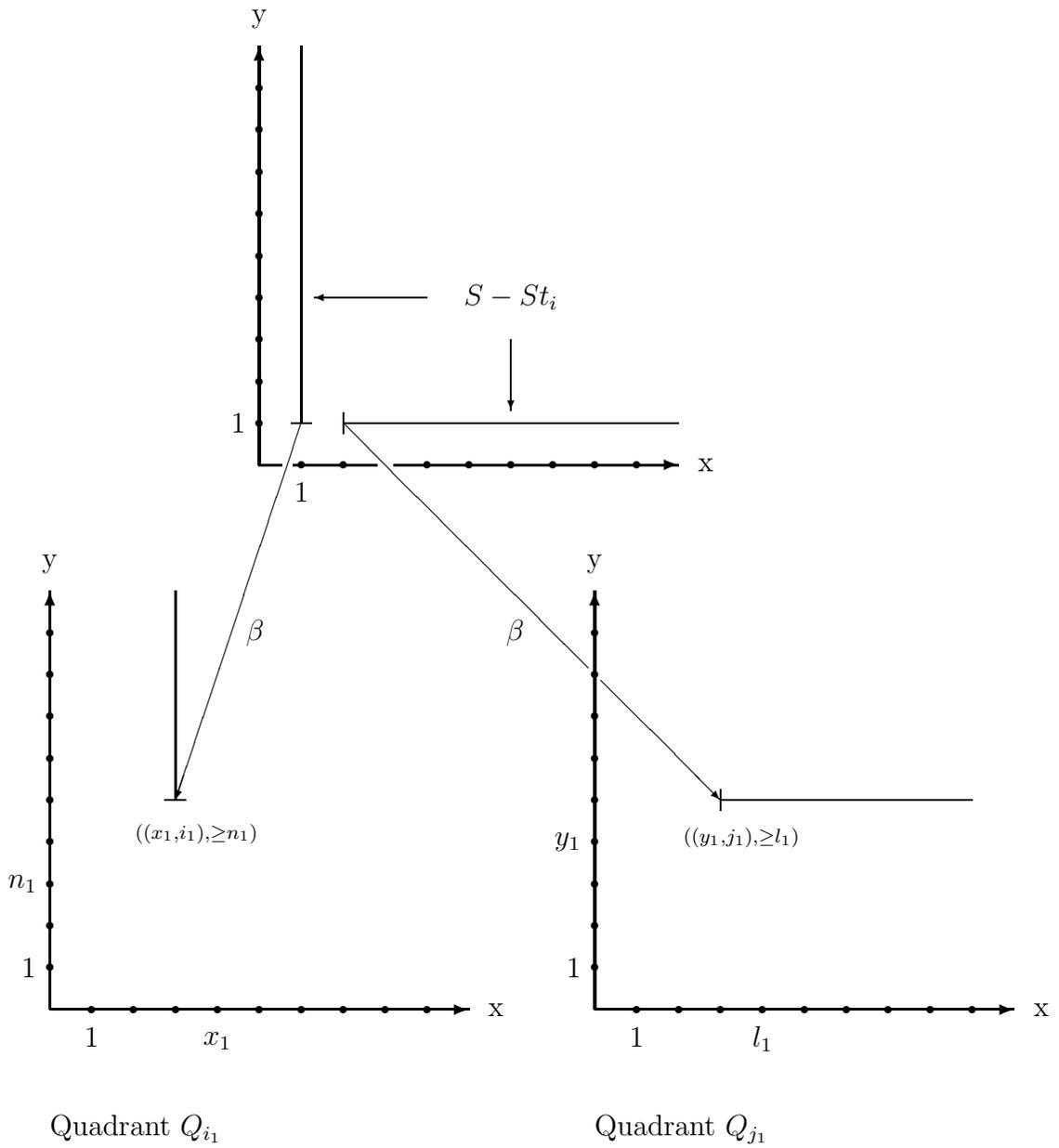
\begin{figure}[ht]
  \unitlength0.6cm
   \begin{center}
   \begin{picture}(27,27)
     \put(7,16){\begin{picture}(10,10)
               \thicklines
               \put(0,0){\line(1,0){0.3}}
               \put(0.8,0){\line(1,0){2}}
               \put(3.2,0){\vector(1,0){6.8}}
               \put(0,0){\line(1,0){0.55}}
               \put(0,0){\vector(0,1){10}}
               \thinlines
               \multiput(1,0)(1,0){2}{\circle*{0.2}}
               \multiput(4,0)(1,0){6}{\circle*{0.2}}
               \multiput(0,1)(0,1){9}{\circle*{0.2}}
               \put(10,0){\makebox(1.3,0){x}}
               \put(0,10){\makebox(0,1.3){y}}
               \put(1,-1.3){\makebox(0,1.3){1}}
               \put(-1,1){\makebox(1,0){1}}
               \put(1,1){\line(0,1){9}}\put(0.75,1){\line(1,0){0.5}}
               \put(2,1){\line(1,0){8}}\put(2,0.75){\line(0,1){0.5}}
               \put(4,3){\makebox(4,2){$S-St_i$}}
               \put(4,4){\vector(-1,0){2.7}}
               \put(6,3){\vector(0,-1){1.7}}
                              \end{picture}}

     \multiput(2,3)(13,0){2}{\begin{picture}(10,10)
               \thicklines
               \put(0,0){\vector(1,0){10}}
               \put(0,0){\vector(0,1){10}}
               \thinlines
               \multiput(1,0)(1,0){9}{\circle*{0.2}}
               \multiput(0,1)(0,1){9}{\circle*{0.2}}
               \put(10,0){\makebox(1.3,0){x}}
               \put(0,10){\makebox(0,1.3){y}}
               \put(1,-1.3){\makebox(0,1.3){1}}
               \put(-1,1){\makebox(1,0){1}}
                              \end{picture}}
     \put(2,3){\begin{picture}(10,10)
               \put(3,5){\line(0,1){5}}\put(2.75,5){\line(1,0){0.5}}
               \put(-1.3,3){\makebox(1.3,0){$n_1$}}
               \put(4,-1.5){\makebox(0,1.5){$x_1$}}
               \put(1.5,3.7){\makebox(4,1){$\scriptstyle((x_1,i_1),\ge n_1)$}}
               \end{picture}}
    \put(15,3){\begin{picture}(10,10)
               \put(3,5){\line(1,0){6}}
               \put(3,4.75){\line(0,1){0.5}}
               \put(-1.3,4){\makebox(1.3,0){$y_1$}}
               \put(4,-1.4){\makebox(0,1.4){$l_1$}}
               \put(2,3.6){\makebox(3,1){$\scriptstyle((y_1,j_1),\ge l_1)$}}
               \end{picture}}
    \put(8,17){\vector(-1,-3){3}}
    \put(9,17){\line(1,-1){5.85}}
    \put(15.15,10.85){\vector(1,-1){2.85}}
    \put(6.4,12){\makebox(1,0){$\beta$}}
    \put(12.6,12){\makebox(1,0){$\beta$}}
    \put(2,0){Quadrant $Q_{i_1}$}
    \put(15,0){Quadrant $Q_{j_1}$}
    \end{picture}
  \end{center}
  \caption{ Definition von $\beta$ auf $S-St_i$}
  \end{figure}

  Damit $\beta$ in $\cal M$ liegt, mu\3 es injektiv sein und die Bedingungen
  B$_1$, B$_2$ erf"ullen.\\

  {\em Injektivit"at\/}:
  Offensichtlich ist $\beta$ auf den Mengen $St_i$ und $S-St_i$ jeweils
  injektiv. Dabei ist $(St_i)\beta=S\alpha$ und $(S-St_i)\beta=((y_1,j_1),\ge l_1)
  \cup ((x_1,i_1),\ge n_1)\subseteq (S-S\alpha)$. Also gilt
      \[  (St_i)\beta\cap(S-St_i)\beta\subseteq S\alpha\cap(S-S\alpha)=
           \emptyset \,,         \]             
  d.h. $\beta$ ist auf ganz $S$ injektiv.\\

  {\em Bedingung\/ {\rm B}$_1$}:
  Au\3erhalb von $Q_i$ stimmt $\beta$ mit $\alpha$ "uberein; wir m"ussen 
  also nur noch das Verhalten von $\beta$ in $Q_i$ untersuchen.
  Verschiebt $\alpha$ die Elemente aus $Q_i$ ab $p_0=(x_0,y_0)$ um den
  Vektor $(m_i,m_i)$, so verschiebt $\beta$ ab $p_0^{'}:\,=(x_0+1,y_0+1)$
  um $(m_i-1,m_i-1)$:\\
   \parbox{9.5cm}
   {\begin{eqnarray*}
       ((x,y),i)\beta&=&((x,y),i)t_i{}^{-1}\alpha\\
                     &=&((x-1,y-1),i)\alpha\\
                     &=&((x+m_i-1,y+m_i-1),i)
    \end{eqnarray*}}
   \hfill
   \parbox{5.5cm}
     {$\forall\;(x,y)\ge (x_0+1,y_0+1)\:.$}
  {\em Bedingung\/ {\rm B$_2$(a)}}:
  Wir identifizieren $((x,y),i)$ aus $S$ mit $((x,i),y)$.
  F"ur $(x,i)\in X\times\{1,\ldots,n\}$ gilt
    \[ ((x,i),y)\alpha=((x,i)^{'},y+q_{(x,i)})\quad\mbox{f"ur alle }\quad
                                              y\ge y_0  \:,  \]
 wobei $(x,i)^{'}\in X\times \{1,\ldots,n\}$ und $q\in {\sl Z}\!\!{\sl Z}$ ist.
  Daraus erh"alt man f"ur alle $x\ge 2$:\\
     \parbox{9.5cm}
     {\begin{eqnarray*}
          ((x,i),y)\beta&=&((x,i),y)t_i{}^{-1}\alpha\\
                        &=&((x-1,i),y-1)\alpha\\
                        &=&((x-1,i)^{'},y+q_{(x-1,i)}-1)
      \end{eqnarray*}}\hfill
     \parbox{5.5cm}{$ \forall\; y\ge y_0+1\:.$}
  F"ur $(1,i)\in X\times\{1,\ldots,n\}$ erf"ullt $\beta$ B$_2$(a) wegen
  (3.11). Die Bedingung B$_2$(b) zeigt man analog. \hfill $\Box$\\
\clearpage   
 {\bf Bemerkung 3.2$\:$:}
      $\:$Das Element $\beta$ aus dem letzten Lemma
  ist nicht eindeutig bestimmt. Wie aus dem Beweis des Lemmas ersichtlich wird,
  gibt es beliebig viele M"oglichkeiten f"ur die Wahl eines solchen $\beta$'s.
  Insbesondere existiert im Falle $gr(\alpha)=1$ ein $\beta$ aus $G_n$
  mit $t_i\beta=\alpha$. Da dann  $S-S\alpha$ von der 
  Form 
     \[  S-S\alpha=((x_1,i_1),\ge n_1)\cup((y_1,j_1,\ge l_1)\cup
                    (P=\{p_1,\ldots,p_r\})                     \]
   ist, k"onnen wir
   die Menge $S-St_i$ {\em surjektiv\/} (s. Abbildung 18) auf $S-S\alpha$ abbilden:
   \[\begin{array}{rcll}
     ((x,1),i)\,\beta&:\,=&((x+l_1-1,y_1),j_1)& \quad\mbox{f"ur alle $x\in{\rm I}\!{\rm N}$}\\
    ((1,y+r+1),i)\beta&:\,=&((x_1,y+n_1-1),i_1)&\quad\mbox{f"ur alle $y\in {\rm I}\!{\rm N}$ und}\\
     ((1,y),i)\beta&:\,=&p_y&\quad\mbox{f"ur }1\le y\le r\,.\\
                   & &                           
       \end{array} \] 

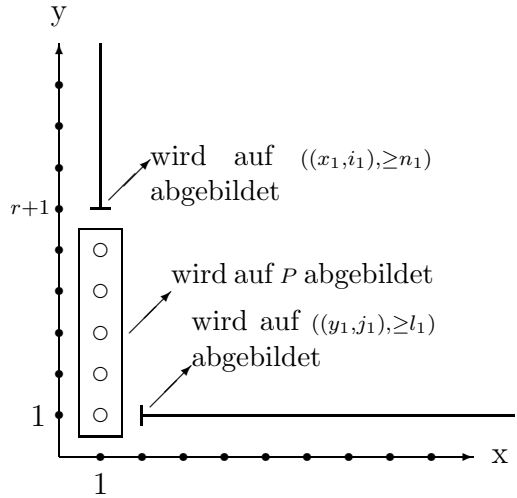
\begin{figure}[ht]
  \unitlength0.55cm
   \begin{center}
   \begin{picture}(11,11.5)
     \put(0,0.5){\begin{picture}(10,10)
               \put(0,0){\vector(1,0){10}}
               \put(0,0){\vector(0,1){10}}
               \multiput(1,0)(1,0){9}{\circle*{0.2}}
               \multiput(0,1)(0,1){9}{\circle*{0.2}}
               \put(10,0){\makebox(1.3,0){x}}
               \put(0,10){\makebox(0,1.3){y}}
               \put(1,-1.3){\makebox(0,1.3){1}}
               \put(-1,1){\makebox(1,0){1}}
                              \end{picture}}
     \put(0,0.5){\begin{picture}(10,10)
               \thicklines
               \put(1,6){\line(0,1){4}}\put(0.75,6){\line(1,0){0.5}}
               \put(2,1){\line(1,0){9}}\put(2,0.75){\line(0,1){0.5}}
               \thinlines
               \multiput(1,1)(0,1){5}{\circle{0.33}}
               \put(0.5,0.5){\framebox(1,5){}}
               \put(-1.4,6){\makebox(1.4,0){$\scriptstyle r+1$}}
               \put(2.2,6.2){\parbox[b]{3.75cm}{\small wird auf $\scriptstyle
                          ((x_1,i_1),\ge n_1)$ abgebildet}}
               \put(2.7,4.2){\parbox[b]{3.5cm}{\small wird auf $\scriptstyle P$
                                         abgebildet}}
               \put(3.2,2.2){\parbox[b]{3.25cm}{\small wird auf $\scriptstyle
                           ((y_1,j_1),\ge l_1)$ abgebildet }}
               \put(1.2,6.2){\vector(1,1){1}}
               \put(1.7,3){\vector(1,1){1}}
               \put(2.2,1.2){\vector(1,1){1}}
               \end{picture}}
   \end{picture}
   \end{center}
   \caption{ Definition von $\beta$ im Falle $gr(\alpha)=1$}
   \end{figure}

  Mit den vorangegangenen beiden Lemmata erhalten wir nun:

  \begin{cor}
     F"ur $\alpha$ aus $\cal M$ ist $gr(\alpha)$ die gr"o\3te
     ganze Zahl $k$, so da\3 es eine Kette der Form \mbox{$\alpha=\alpha_k>\alpha_{k-1}>\ldots>\alpha_0$}
     in $\cal M$ gibt.
 \end{cor}

{\sc Beweis}:$\:$ Sei $\alpha=\alpha_k>\alpha_{k-1}>\ldots>\alpha_0$ eine Kette ausgehend von $\alpha$ mit 
      maximaler L"ange $k$. Dann gilt f"ur alle $i=1,\ldots,k$ (vgl. (3.7))
      \[ \alpha_i=t_{j_i}\alpha_{i-1} \:. \]
      Mit Lemma 3.6 folgt daraus
      $\, gr(\alpha_i)=gr(\alpha_{i-1})+1 \,, $
    und somit 
       \[  gr(\alpha)=gr(\alpha_k)
                   =gr(\alpha_0)+k \:.\]
  Wegen Lemma 3.7 mu\3 $gr(\alpha_0)=0$
    sein (sonst lie\3e sich die gegebene Kette verl"angern), d.h.
    $gr(\alpha)=k$.\hfill $\Box$\\

 {\bf Bemerkung 3.3$\:$:}$\;$Die Funktion $\:gr:\cal M \rightarrow {\rm I}\!{\rm N}$$_0$, die jedem $\alpha$ aus 
   $\cal M$ den Grad $gr(\alpha)$ zuordnet, ist wegen dem letzten Korollar ordnungserhaltend
  (bzw sogar streng monoton wachsend).

\subsection{Der Simplizialkomplex $\protect{\mid{\cal M}^n\mid}$}

Es sei $\;{\cal M}^n:\,=\{\,\alpha\in{\cal M}\mid gr(\alpha)\ge n\,\}\;$.
 ${\cal M}^n$ ist als
Teilmenge von $\cal M$ ebenfalls eine partiell geordnete Menge, und der zu ${\cal M}^n$ geh"orige 
Simplizialkomplex $\mid{\cal M}^n\mid$ ist ein Unterkomplex von $\mid\cal M\mid$.
 Wir zeigen, da\3 $\mid{\cal M}^n\mid$ invariant 
 unter der Operation von $G_n$ ist:

\begin{lemma} F"ur alle $\alpha$ aus $\cal M$ sowie $g$ aus $G_n$ ist
               {$gr(\alpha g)=gr(\alpha)\,.$}
\end{lemma}

{\sc Beweis}:$\:$ Jede Kette der Form 
     \begin{itemize}
       \item[(a)]$\;\alpha=\alpha_k>\alpha_{k-1}>\ldots>\alpha_0\quad\mbox{mit}
               \quad t_{i_j}\alpha_{j-1}=\alpha_j $
     \end{itemize}
               f"ur alle $j=1,\ldots,k$ liefert eine Kette ($G_n$ operiert
               ordnungserhaltend auf $\cal M$)
     \begin{itemize}
       \item[(b)]$\;\alpha g=\alpha_k g>\alpha_{k-1}g>\ldots>\alpha_0 g
             \quad\mbox{mit}\quad t_{i_j}\alpha_{j-1}=\alpha_j$
     \end{itemize}
    f"ur $j=1,\ldots,k$. Umgekehrt erh"alt man aus jeder Kette ausgehend von $\alpha g$ 
    durch Multiplikation mit $g^{-1}$ eine Kette ausgehend von $\alpha$
    derselben L"ange $k$. Aus
    Korollar 3.8 somit folgt die Behauptung.\hfill $\Box$\\

 {\bf Bemerkung 3.4:}$\;$ Da f"ur $\alpha_1,\alpha_2\in \cal M$ mit $\alpha_1
  <\alpha_2$ auch $\alpha_1\beta<\alpha_2\beta$ ist ($t\alpha_1=\alpha_2
  \Rightarrow t(\alpha_1\beta)=(t\alpha_1)\beta=\alpha_2\beta$), gilt
  $gr(\alpha)\le gr(\alpha\beta)$ f"ur alle $\alpha,\beta$ aus $\cal M$.\\

 Wegen Lemma 3.9 ist $\mid{\cal M}^n\mid$  ein $G_n$-Komplex, wobei
 nach den Ergebnissen aus 3.2.1 (Korollar 3.5) die Stabilisatoren jedes
 p-Simplexes aus $\mid{\cal M}^n\mid$ mindestens vom Typ $(FP)_n$ und f"ur
  $n\ge 3$ (eigentlich schon f"ur $n\ge 2$) endlich pr"asentiert sind.\\
 Es bleibt zu zeigen, da\3 der Komplex $\mid{\cal M}^n\mid$ ebenfalls
  zusammenziehbar ist. Wir wissen nach Lemma 3.1,
 da\3 f"ur alle $\alpha$ aus ${\cal M}^n$ ein $t\in T$ existiert
  mit $\alpha\le t$. Nun gilt ($gr$ ist ordnungserhaltend, s. Bemerkung 3.3) 
    \[ gr(t)\ge gr(\alpha)\ge n \,,  \]
 d.h. $t$ liegt in ${\cal M}^n$. Daraus ergibt sich wie im Beweis von Korollar 3.2,
 da\3 ${\cal M}^n$ eine gerichtete Menge und somit $\mid{\cal M}^n\mid$ zusammenziehbar ist.\\

{\bf Bemerkung 3.5$\:$:}$\:$ Die maximal m"ogliche L"ange
 einer Kette der Form $\alpha=\alpha_k>\alpha_{k-1}>\ldots>\alpha_0$ 
mit $\alpha_i\in {\cal M}^n$ f"ur alle $i=0,\ldots,k$ ist gerade $gr(\alpha)-n$.\newpage

\section{Filtrierung von $\protect{\mid{\cal M}^n\mid}$}

"Ahnlich wie bei den Houghtongruppen filtrieren wir den $G_n$-Komplex $\mid{\cal M}
 ^n\mid$ mit Hilfe der Teilmengen
     \[ {\cal M}^{[n,k]}:\,=\{\alpha\in {\cal M}\mid n\le gr(\alpha)\le k\} \]
von ${\cal M}^n$, wobei $k$ eine nat"urliche Zahl $\ge n$ ist. Die Familie
 der Unterkomplexe $\{\mid{\cal M}^{[n.k]}\mid\}_{k\ge n}$ bildet eine 
 Filtrierung von $\mid{\cal M}^n\mid$, wobei sich die $G_n$-Invarianz der
 einzelnen $\mid{\cal M}^{[n,k]}\mid$ aus der Beziehung $\,gr(\alpha g)
 =gr(\alpha)\,$ (vgl. Lemma 3.9) ergibt.\\
 Um Korollar 1.3 anwenden zu k"onnen, bleibt zu zeigen, da\3 jeder dieser
 Unterkomplexe $\mid{\cal M}^{[n,k]}\mid$ mod $G_n$ endlich ist, und da\3 
 f"ur hinreichend gro\3e $k$ $\mid{\cal M}^{[n,k+1]}\mid$ bis auf Homotopie
 aus $\mid{\cal M}^{[n,k]}\mid$ durch Hinzuf"ugen von $n$-Zellen entsteht.
 Dies wollen wir in den n"achsten beiden Abschnitten 4.1 und 4.2 tun.

 \subsection{ Die Unterkomplexe der Filtrierung mod $G_n$ }

F"ur ein $p$-Simplex $(\alpha_0,\ldots,\alpha_p)$ aus $\mid{\cal M}^n\mid$
 gibt es zu jedem $i\in \{1,\ldots,p\}$ ein $t_{\alpha_i}$ aus $T$ mit
     \[  \alpha_i=t_{\alpha_i}\alpha_{i-1}\: .  \]
 Dabei ist $t_{\alpha_i}=\alpha_i(\alpha_{i-1})^{-1}$, also eindeutig bestimmt.
 Wir werden nun die Bahn von $(\alpha_0,\ldots,\alpha_p)$ unter der Operation
 von $G_n$ mit Hilfe des $p$-Tupels $(t_{\alpha_1},\ldots,t_{\alpha_p})\in
 T^p$ charakterisieren.

 \begin{lemma}
    Seien $(\alpha_0,\ldots,\alpha_p)$ und $(\beta_0,\ldots,\beta_p)$
    zwei $p$-Simplizes aus $\mid {\cal M}^n\mid$. Dann sind folgende
    Aussagen "aquivalent:
   \begin{itemize}
     \item[{\rm (i)}]$\;$ Es gibt ein $g\in G_n\;$ mit $\;(\alpha_0,\ldots,\alpha_p)g
                                              =(\beta_0,\ldots,\beta_p)$
     \item[{\rm (ii)}]$\;$ $gr(\alpha_0)=gr(\beta_0)\,$ und 
            $\;(t_{\alpha_1},\ldots,t_{\alpha_p})
            =(t_{\beta_1},\ldots,t_{\beta_p})$.
   \end{itemize}
 \end{lemma}
 \vspace{0.3mm}
 {\sc Beweis}:$\;$(i)$\,\Rightarrow\,$(ii)\\
 F"ur alle $i=0,\ldots,p$ ist $\alpha_ig=\beta_i$.
  Wegen  Lemma 3.9 gilt somit $gr(\alpha_0)=gr(\alpha_0g)=gr(\beta_0)$. Weiterhin 
 erh"alt man f"ur alle $i=1,\ldots,p$
   \begin{eqnarray*}
      t_{\alpha_i}\beta_{i-1}&=&t_{\alpha_i}(\alpha_{i-1}g)\\
                             &=&(t_{\alpha_i}\alpha_{i-1})g\\
                             &=& \alpha_ig=\beta_i\:,
   \end{eqnarray*}
 woraus wegen der Eindeutigkeit von $t_{\beta_i}$ folgt: $t_{\alpha_i}
  =t_{\beta_i}$.\\
 (ii)$\,\Rightarrow\,$(i):$\;$
 Wir konstruieren ein $g\in G_n$  mit $\alpha_0g=\beta_0$. Daraus folgt
 wegen $(t_{\alpha_1},\ldots,t_{\alpha_p})=(t_{\beta_1},\ldots,t_{\beta_p})$
 f"ur alle $i=1,\ldots,p$:
    \begin{eqnarray*}
       \alpha_ig&=&t_{\alpha_i}t_{\alpha_{i-1}}\cdot\,\ldots\,\cdot t_{\alpha_1}\alpha_0g\\
                 &=&t_{\beta_i}t_{\beta_{i-1}}\cdot\,\ldots\,\cdot t_{\beta_1}\beta_0\\
                 &=&\beta_i\:.
    \end{eqnarray*}
 Konstruktion von $g$:$\:$Sei $k:=gr(\alpha_0)=gr(\beta_0)$. Wir wenden
 nun ausgehend von $\alpha_0$ bzw. $\beta_0$ $k$-mal Lemma 3.7
 mit $i=1$ an und erhalten jeweils eine Kette der Form
   \begin{eqnarray*}
      \alpha_0&=&\alpha_{0,k}>\alpha_{0,k-1}>\ldots >\alpha_{0,1}>\alpha_{0,0}
                 \quad\mbox{bzw.}\\
      \beta_0&=&\beta_{0,k}>\beta_{0,k-1}>\ldots >\beta_{0,1}>\beta_{0,0}\:,
   \end{eqnarray*}
  wobei $\alpha_{0,j}=t_1\alpha_{0,j-1}$ bzw. $\beta_{0,j}=t_1\beta_{0,j-1}$
  f"ur alle $j=1,\ldots,k$ ist. Nach Lemma 3.6 ist $gr(\alpha_{0,j})=j$,
  also insbesondere $gr(\alpha_{0,1})=1$, weshalb es m"oglich ist, $\alpha_{0,0}$
  aus $G_n$ zu w"ahlen (s. Bemerkung 3.2). Dasselbe gilt f"ur 
  $\beta_{0,0}$. Dann ist 
     \[ g:\,=(\alpha_{0,0})^{-1}\beta_{0,0}  \]
 aus $G_n$, und wegen $\alpha_0=t_1^k\alpha_{0,0}$ bzw. $\beta_0=t_1^k
  \beta_{0,0}$ erhalten wir
    $\,  \alpha_0g=t_1^k\alpha_{0,0}g=t_1^k\alpha_{0,0}(\alpha_{0,0})^{-1}
                  \beta_{0,0}=t_1^k\beta_{0,0}=\beta_0.$   
  \hfill $\Box$\\
 
 \begin{beh}
   Die Menge $T_{\le k}:\,=\{t\in T\mid gr(t)\le k\}$ ist endlich.
 \end{beh}
 {\sc Beweis}:$\;$ $T$ wird erzeugt von $t_1,\ldots,t_n$, d.h. jedes $t\in
 T$ l"a\3t sich schreiben als
     \[ t=d_1\cdot\,\ldots\,\cdot d_m\quad\mbox{mit}\quad d_i\in\{t_1,\ldots,t_n\}\:. \]
 Daraus ergibt sich eine Kette der Form
   $ t=d_1\cdot\,\ldots\,\cdot d_m>d_2\cdot\,\ldots\,\cdot d_m>\ldots>d_m>id\,,$
 woraus nach Korollar 3.8 folgt:$\quad gr(t)\ge m$. Somit ist
    \[ T_{\le k}\subseteq \{ d_1\cdot\,\ldots\,\cdot d_m\in T \mid d_i\in 
                             \{t_1,\ldots,t_n\},i=1,\ldots,m\quad\mbox{und}
                              \quad m\le k \}\:.\]
 Letztere Menge ist endlich, d.h. auch $T_{\le k}$ ist endlich. \hfill $\Box$\\
    
 \begin{cor}
    Der Simplizialkomplex $\mid{\cal M}^{[n,k]}\mid$ ist endlich mod $G_n$.
 \end{cor}

 {\sc Beweis}:$\;$ Wegen Bemerkung 3.5 ist die Dimension von $\mid
  {\cal M}^n\mid$ endlich. Es gen"ugt daher zu zeigen, da\3 die Menge aller
  $p$-Simplizes (f"ur festes $p$) endlich ist mod $G_n$.\\
  Wir ordnen dazu der Bahn eines $p$-Simplexes $(\alpha_0,\ldots,\alpha_p)$
  das Paar 
     \[ (\,gr(\alpha_0),(t^{(1)},\ldots,t^{(p)})\,)\in [n,k]\times (T^p) \]
  zu ($\,[n,k]=\{j\in {\rm I}\!{\rm N}\mid n\le j\le k\}\,$),
  wobei $t^{(i)}\alpha_{i-1}=\alpha_i$, $1\le i\le p$ ist. Nach dem vorigen
  Lemma ist diese Zuordnung wohldefiniert und injektiv. Nun gilt
     \begin{eqnarray*}
          gr(t^{(i)})&\le& gr(t^{(i)}\alpha_{i-1})\quad\mbox{(vgl. Bemerkung 3.4)}\\
                 &=& gr(\alpha_i)\\
                 &\le & k \quad\mbox{f"ur alle}\quad i=1,\ldots,p\:,
    \end{eqnarray*}
   d.h. $(t^{(1)},\ldots,t^{(p)})\in (T_{\le k})^p$, also $(\,gr(\alpha_0),(t^{(1)},
   \ldots,t^{(p)})\,)\in [n,k]\times(T_{\le k})^p$. Nach der letzten Behauptung
   ist $(T_{\le k})^p$ und damit auch $[n,k]\times (T_{\le k})^p$ endlich.
   Wegen der Injektivit"at obiger Zuordnung  kann es demnach nur endlich viele verschiedenen Bahnen geben.
   \hfill $\Box$\\

 \subsection{Homotopieeigenschaften von $\protect{\left|{\cal M}^{[n,k]}
             \right|\subseteq\left|{\cal M}^{[n,k+1]}\right|}$}

Wir ben"otigen im folgenden zu gegebenem $\alpha$ aus $\cal M$ die Mengen
   \begin{eqnarray*}
      {\cal M}_{<\alpha}&:\,=&\{\beta\in {\cal M}\mid \beta<\alpha\}\quad\mbox{bzw.} \\
      {\cal M}_{\le \alpha}&:\,=&\{ \beta\in {\cal M}\mid \beta \le \alpha\}
             \quad\mbox{und} \\
      {\cal M}_{<\alpha}^n&:\,=& {\cal M}^n\cap{\cal M}_{<\alpha}\:.
   \end{eqnarray*}
 Im letzten Abschnitt dieser Arbeit zeigen wir, da\3 f"ur ein $\alpha$ aus
 $\cal M$ mit $gr(\alpha)\ge 2n$ der Simplizialkomplex $\mid{\cal M}_{<\alpha}
 ^n\mid$ homotopie"aquivalent zu einem $(n-1)$-dimensionalen Sph"arenbouquet
 ist.\\
 Mit derselben Argumentation wie bei den Houghton-Gruppen ($\mid{\cal M}^{[n,k+1]}
  \mid$ entsteht aus $\mid{\cal M}^{[n,k]}\mid$ durch Hinzuf"ugen eines 
  Kegels "uber ${\cal M}_{<\alpha}^n$ f"ur jedes $\alpha$ mit $gr(\alpha)=k+1$,
 vgl. Abschnitt 1.3, Schritt 2), ergibt sich daraus, da\3 man $\mid{\cal M}^
 {[n,k+1]}\mid$ bis auf Homotopie aus $\mid {\cal M}^{[n,k]}\mid$ durch
 Hinzuf"ugen von $n$-Zellen erh"alt.\\
 Der Nachweis dieser Homotopie"auquivalenz geschieht in zwei Schritten.
 Zuerst wird eine "Uberdeckung von $\mid {\cal M}_{<\alpha}^n\mid$ konstruiert,
 deren Nerv homotopie"aquivalent zu $\mid{\cal M}_{<\alpha}^n\mid$ ist.
 Anschlie\3end zeigen wir, da\3 dieser Nerv einer speziellen Art kombinatorischer
 Simplizialkomplexe angeh"ort, die vom Homotopietyp eines $(n-1)$-dimensionalen
 Sph"arenbouquets sind.

\subsubsection{Eine "Uberdeckung von $\protect{\mid{\cal M}_{<\alpha}^n\mid}$}

 Wir betrachten die maximalen Elemente in ${\cal M}_{<\alpha}^n$, d.h.
 solche Elemente $\beta$ aus $\cal M$ mit der Eigenschaft
      \[  t_i\beta=\alpha  \]
 f"ur ein $i\in\{1,\ldots,n\}$. $\mid{\cal M}_{<\alpha}^n\mid$ kann
  durch die Unterkomplexe $\mid{\cal M}_{\le \beta}^n\mid$ "uberdeckt werden,
 wobei $\beta$ die Menge aller maximalen Elemente von ${\cal M}_{<\alpha}^n$
 durchl"auft. Den Nerv dieser "Uberdeckung bezeichnen wir mit 
 ${\cal N}_{\alpha}$.\\
 Um zu zeigen, da\3 ${\cal N}_{\alpha}$ homotopie"aquivalent zu $\mid{\cal M}
 _{<\alpha}^n\mid$ ist, gen"ugt es nach \mbox{Korollar 1.8}, zu jeder endlichen Menge
 $\{\beta_1,\ldots,\beta_p\}$ bestehend aus maximalen Elementen von ${\cal M}
 _{<\alpha}^n$, die in ${\cal M}^n$ nach unten beschr"ankt ist, die Existenz
 einer {\em gr"o\3ten\/} unteren Schranke in ${\cal M}^n$ nachzuweisen.\\
 Zu diesem Zweck beschreiben wir das Vorhandensein einer unteren Schranke der
 $\{\beta_1,\ldots,\beta_p\}$ durch eine dazu "aquivalente Eigenschaft
 der maximalen Elemente $\beta_i$.\\

 {\bf Definition:}$\:$ Zu $\beta$ aus $\cal M$
 mit $t_i\beta=\alpha$ sei $\bar{\beta}_i$ die Einschr"ankung von $\beta$ auf die Menge $S-St_i$.
 B($\bar{\beta}_i):\,=(S-St_i)\bar{\beta}_i$ sei die Bildmenge von $\bar{\beta}_i$.\\

 Wegen $(St_i)\beta=
 S\alpha$ mu\3 
  \begin{equation} {\rm B}(\bar{\beta}_i)\subseteq(S-S\alpha)
  \end{equation}
 sein. Mit Hilfe von $\bar{\beta}_i$ und B($\bar{\beta}_i$) erhalten wir
 folgendes Kriterium:\\
 \begin{lemma}
    Sei $\alpha$ aus $\cal M$ und $\{\beta_1,\ldots,\beta_p\}$ eine Menge
     maximaler Elemente aus ${\cal M}_{<\alpha}^n$. F"ur $gr(\alpha)\ge 2n$
    sind folgende Aussagen "aquivalent:
  \begin{itemize}
    \item[{\rm (i)}]$\;$ $\{\beta_1,\ldots,\beta_p\}$ ist in ${\cal M}^n$
                    nach unten beschr"ankt
    \item[{\rm (ii)}]$\;$ F"ur $\bar{\beta}_{i_1},\ldots,\bar{\beta}_{i_p}$ und
                    {\rm B}$(\bar{\beta}_{i_1}),\ldots,{\rm B}(\bar{\beta}_{i_p})$ gilt
                    f"ur alle $\nu,\mu\in\{1,\ldots,p\}$, $\:\nu\not=\mu$:
                       \[ i_{\nu}\not=i_{\mu}\quad\mbox{ und }\quad
                     {\rm B}(\bar{\beta}_{i_{\nu}})
                    \cap {\rm B}(\bar{\beta}_{i_{\mu}})=\emptyset.  \]
    \end{itemize}
  \end{lemma}
 {\sc Beweis}:$\;$(i)$\,\Rightarrow\,$(ii):$\;$
  Sei $\delta\in {\cal M}^n$ die gemeinsame untere Schranke der $\beta_i$'s,
  also
         \[  d_i\delta=\beta_i\quad\mbox{mit}\quad d_i\in T  \]
  f"ur alle $i=1,\ldots,p$. Angenommen $t_{i_{\nu}}\beta_{\nu}=\alpha$
  und $t_{i_{\nu}}\beta_{\mu}=\alpha$. Dann ist
    \[  \begin{array}{ccrcll}
          & & t_{i_{\nu}}\beta_{\nu}&=&t_{i_{\nu}}\beta_{\mu}& \\
      &\Leftrightarrow& t_{i_{\nu}}d_{\nu}\delta&=&t_{i_{\nu}}d_{\mu}\delta\\
      &\Leftrightarrow& t_{i_{\nu}}d_{\nu}&=&t_{i_{\nu}}d_{\mu}
                         &\mbox{($\delta$ injektiv)}\\
      &\Leftrightarrow& d_{\nu}t_{i_{\nu}}&=&d_{\mu}t_{i_{\nu}}
                         &\mbox{($T$ kommutativ)}\\
      &\Leftrightarrow& d_{\nu}&=&d_{\mu}
                         &\mbox{($t_{i_{\nu}}$ injektiv)}\\
      &\Leftrightarrow& \beta_{\nu}&=&\beta_{\mu}&\mbox{bzw}
                        \quad \nu=\mu\:.
    \end{array}  \]

 Somit ist f"ur $\nu\not=\mu$ auch $i_{\nu}\not=i_{\mu}$. Wir betrachten
  nun die Quadranten $Q_j$, die jeweils $T$-invariante Teilmengen von $S$
  sind. Es gilt also f"ur $\nu\not=\mu$
   
    \[ \begin{array}{ccrcl}
     & &  Q_{i_{\nu}}d_{\nu}&\cap& Q_{i_{\mu}}d_{\mu}=\emptyset\\
     &\Rightarrow&  Q_{i_{\nu}}d_{\nu}\delta&\cap& Q_{i_{\mu}}d_{\mu}\delta
                    =\emptyset\quad\mbox{($\delta$ injektiv)}\\
    &\Leftrightarrow& Q_{i_{\nu}}\beta_{\nu}&\cap& Q_{i_{\mu}}\beta_{\mu}=\emptyset\:.
     \end{array}\]
 Nun ist B($\bar{\beta}_{i_{\nu}})=(S-St_{i_{\nu}})\beta_{\nu}\subseteq
 Q_{i_{\nu}}\beta_{\nu}$ (da $\,(S-St_{i_{\nu}})\subseteq Q_{i_{\nu}}$, vgl. Abbildung 16),
 und analog B($\bar{\beta}_{i_{\mu}})\subseteq Q_{i_{\mu}}\beta_{\mu}$,
 d.h. B($\bar{\beta}_{i_{\nu}})\cap$B($\bar{\beta}_{i_{\mu}})=\emptyset$.\\

 (ii)$\,\Rightarrow\,$ (i):
 Wir definieren die gesuchte untere Schranke $\delta$ mit Hilfe von $\alpha$ und
 den einzelnen $\bar{\beta}_i$'s.\\
 Ist $t_{i_{\nu}}\beta_{\nu}=\alpha$ f"ur $\nu=1,\ldots,p$, so sei $\delta$
 auf der Menge $St_{i_1}\cdot\,\ldots\,\cdot t_{i_p}$  die Abbildung
  $t_{i_1}^{-1}\cdot\,\ldots\,\cdot t_{i_p}^{-1}\alpha$. Die verbleibende Menge
  $S-St_{i_1}\cdot\,\ldots\,\cdot t_{i_p}$ ist die disjunkte Vereinigung der
  $S-St_{i_{\nu}}$, $\nu=1,\ldots,p$.
 ($S=Q_1\cup\ldots\cup Q_n$. Au\3erhalb von $Q_{i_1}\cup\ldots\cup Q_{i_p}$
  ist $t_{i_1}\cdot\,\ldots\,\cdot t_{i_p}$ die Identit"at, d.h. die Quadranten
  $Q_j$, $j\in \{1,\ldots,n\}-\{i_1,\ldots,i_p\}$ sind in $St_{i_1}\cdot
  \ldots\cdot t_{i_p}$ enthalten. Somit ist 
     \begin{eqnarray*}
       S-St_{i_1}\cdot\,\ldots\,\cdot t_{i_p}
             &=&(Q_{i_1}\cup\ldots\cup Q_{i_p})-
              (Q_{i_1}\cup\ldots\cup Q_{i_p})t_{i_1}\cdot\,\ldots\,\cdot t_{i_p}\\
             &=&(Q_{i_1}\cup\ldots\cup Q_{i_p})
                -(Q_{i_1}t_{i_1}\cup\ldots\cup Q_{i_p}t_{i_p})\\
             & & \mbox{($t_{i_j}$ ist au\3erhalb von $Q_{i_j}$ die Identit"at)}\\
             &=& (Q_{i_1}-Q_{i_1}t_{i-1})\cup\ldots\cup(Q_{i_p}-Q_{i_p}t_{i_p})\\
             & & \mbox{( die einzelnen $Q_{i_j}$'s sind jeweils $t_{i_j}$-
                         invariant)}\:,
     \end{eqnarray*}
     wobei $Q_{i_{\nu}}-Q_{i_{\nu}}t_{i_{\nu}}=S-St_{i_{\nu}}$ ist (vgl.
    Abbildung 16)).
  Wir definieren $\delta$ auf $S-St_{i_{\nu}}$ als $\bar{\beta}_{i_{\nu}}$,
  $\nu=1,\ldots,p$.\\

  {\em Behauptung\/} 1:$\;\; \delta\in {\cal M}^n$\\
  {\em Beweis\/}:$\:$  $\delta$ ist injektiv auf $S-St_{i_1}\cdot\,\ldots\,\cdot t_{i_p}$,
  da die einzelnen $t_{i_{\nu}}$'s und $\alpha$ injektiv sind, und $\delta$ ist
  injektiv auf $S-St_{i_{\nu}}$ f"ur $\nu=1,\ldots,p$, da die $\bar{\beta}_
  {i_{\nu}}$'s injektiv sind. Wegen
     \[  (St_{i_1}\cdot\,\ldots\,\cdot t_{i_p})\delta=S\alpha\quad\mbox{und} \]
     \[  (S-St_{i_1}\cdot\,\ldots\,\cdot t_{i_p})\delta=\bigcup_{\nu=1}^p 
                  (S-St_{i_{\nu}})\bar{\beta}_{i_{\nu}}
              \subseteq S-S\alpha\quad\mbox{(s. (4.12))}\:,\]
  ist $(St_{i_1}\cdot\,\ldots\,\cdot t_{i_p})\delta\cap (S-St_{i_1}\cdot\ldots
  \cdot t_{i_p})\delta=\emptyset$, d.h. $\delta$ ist injektiv auf ganz $S$.
  Das $\delta$ die Bedingungen B$_1$ und B$_2$ erf"ullt, zeigt man "ahnlich
  wie im Beweis von Lemma 3.7. Es gilt also:$\:\delta\in {\cal M}$.
  Nun ist 
    \begin{eqnarray*}
      gr(\delta)&=& gr(t_{i_1}\cdot\,\ldots\,\cdot t_{i_p}\delta)-p\quad
                    \mbox{(nach Lemma 3.6)}\\
                &=& gr(\alpha)-p\\
                &\ge& 2n-p\\
                &\ge& 2n-n=n\quad \mbox{(wegen $i_{\nu}\not=i_{\mu}$ ist $p\le n$)}\:,
    \end{eqnarray*}
 d.h. $\delta\in {\cal M}^n$.\\\\
 {\em Behauptung\/} 2:$\;\delta$ ist gemeinsame untere Schranke der $\beta_1,\ldots,
                   \beta_p$\\
 {\em Beweis\/}:$\:$Es sei $\delta_j:\,=t_{i_1}\cdot\,\ldots\,\cdot t_{i_{j-1}}t_{i_{j+1}}
       \cdot\,\ldots\,\cdot t_{i_p}\delta$, $j=1,\ldots,p$. Wir zeigen:
       $\delta_j=\beta_j$.\\
     F"ur $\delta_j$ gilt:
       \[ t_{i_j}\delta_j=t_{i_1}\cdot\,
                          \ldots\,\cdot t_{i_p}\delta=\alpha=t_{i_j}\beta_j\:,\]
    d.h. $\delta_j$ und $\beta_j$ stimmen auf der Menge $St_{i_j}$ "uberein.
    Auf der Menge $S-St_{i_j}\subseteq Q_{i_j}$ ist $\delta_j$ gerade $\delta$
    ($\;t_{i_1}\cdot\,\ldots\,\cdot t_{i_{j-1}}t_{i_{j+1}}\cdot\,\ldots\,\cdot t_{i_p}$
      ist eingeschr"ankt auf $Q_j$ die Identit"at), und $\delta$ ist nach
    Definition $\bar{\beta}_{i_j}$. \hfill $\Box$\\
 \vfill  
  Im Falle (ii) aus Lemma 4.4 l"a\3t sich nun eine gr"o\3te gemeinsame
  untere Schranke der $\beta_i$'s angeben.
 \vfill
\begin{lemma}
     Sind die Voraussetzungen von Lemma 4.4 sowie {\rm (ii)} erf"ullt, so gibt
     es in ${\cal M}^n$ eine gr"o\3te gemeinsame untere Schranke der
     $\beta_1,\ldots,\beta_p$.
 \end{lemma}
 \vfill
 {\sc Beweis}:$\;$ Wir zeigen, da\3 die im Beweis von Lemma 4.4 konstruierte
  untere Schranke $\delta$ schon die gr"o\3te untere Schranke von $\{\beta_1,
  \ldots,\beta_p\}$ ist. F"ur $j=1,\ldots,p$ war
     \begin{equation}
       t_{i_1}\cdot\,\ldots\,\cdot t_{i_{j-1}}t_{i_{j+1}}\cdot\,\ldots\,\cdot 
          t_{i_p}\delta=\beta_j  \;.
     \end{equation}
 \pagebreak
 Sei nun $\gamma$ eine weitere untere Schranke von $\{\beta_1,\ldots,\beta_p\}$.
 Dann existieren $w_1,\ldots,w_p\in T$ mit
       \[  w_j\gamma=\beta_j\quad\mbox{f"ur}\quad j=1,\ldots,p\:.  \]
 Wir definieren ein $t\in T$ durch
     \[ ((x,y),i)t:\,=\left
          \{\begin{array}{rcl}
            ((x,y),i)w_1 & \quad\mbox{f"ur}\quad & i\in\{1,\ldots,n\}-
                                                        \{i_1,\ldots,i_p\}\\
            ((x,y),i)w_j & \quad\mbox{f"ur}\quad &\mbox{$ j=1,\ldots,p\:.$}
           \end{array}\right.  \]                                   
 Dann gilt f"ur $i\in\{1,\ldots,n\}-\{i_1,\ldots,i_p\}$:
    \begin{eqnarray*}
        ((x,y),i)t\gamma&=&((x,y),i)w_1\gamma\\
                        &=&((x,y),i)\beta_1\\
                        &=&((x,y),i)t_{i_2}\cdot\,\ldots\,\cdot t_{i_p}\delta
                           \quad\mbox{(nach (4.13))} \\
                        &=&((x,y),i)\delta\quad\mbox{(wegen
                                                $i\notin\{i_1,\ldots,i_p\}$)}
                           \:,
    \end{eqnarray*}
 und f"ur $i_1,\ldots,i_p$:
    \begin{eqnarray*}
       ((x,y),i_j)t\gamma&=&((x,y),i_j)w_j\gamma\\
                         &=&((x,y),i_j)\beta_j\\
                         &=&((x,y),i_j)t_{i_1}\cdot\,\ldots\,\cdot t_{i_{j-1}}
                            t_{i_{j+1}}\cdot\,\ldots\,\cdot t_{i_p}\delta
                            \quad\mbox{(nach (4.13))}\\
                         &=&((x,y),i_j)\delta \quad\mbox{( wegen $i_j\notin
                             \{i_1,\ldots,i_{j-1},i_{j+1},\ldots,i_p\}$)}\:.
   \end{eqnarray*}
 Somit ist $t\gamma=\delta$ bzw. $\gamma\le \delta$.\hfill $\Box$\\\\
Also ist ${\cal N}_{\alpha}$ homotopie"aquivalent zu $\mid{\cal M}^n_{<\alpha}\mid$.
  \vfill
\subsubsection{Homotopietyp des Nerves}
  \vfill
Zur Bestimmung des Homotopietypes ist eine andere Beschreibung von
 ${\cal N}_{\alpha}$ notwendig.
   \vfill
\begin{lemma}
   F"ur $gr(\alpha)\ge 2n$ ist ${\cal N}_{\alpha}$ hom"oomorph zu folgendem
   Simplizialkomplex $\Sigma_{\alpha}$:$\;$ Die Eckenmenge von $\Sigma_
   {\alpha}$ ist 
        \[A:\,=A_1\cup A_2\cup\ldots\cup A_n\,,\] wobei f"ur $i\in\{1,
   \ldots,n\}$ $A_i$ die Menge aller Abbildungen $\varphi:(S-St_i)
   \longrightarrow (S-S\alpha)$ ist, die sich zu einer Abbildung von $S$ nach
   $S$ aus $\cal M$ fortsetzen lassen.
   Bezeichnet man zu gegebenem $a_i\in A_i$ mit {\rm B}$(a_i)$ die Menge $(S-St_i)
   a_i$, so bestehen die Simplizes von $\Sigma_{\alpha}$ aus den endlichen
   Teilmengen
     \[ \{a_{i_0},\ldots,a_{i_p}\}\subseteq A \]
  von $A$, f"ur die $i_{\nu}\not=i_{\mu}\;$ und $\;{\rm B}(a_{i_{\nu}})\cap
      {\rm B}(a_{i_{\mu}})=\emptyset$
 f"ur alle $0\le \nu,\mu\le p$, $\nu\not=\mu$ ist.
 \end{lemma}
   \pagebreak
{\sc Beweis}:$\;$Wir betrachten die Abbildung
    \[   f:\;\{\beta\in {\cal M}\mid t_i\beta=\alpha\quad\mbox{f"ur ein $i\in
         \{1,\ldots,n\}$} \}\,\longrightarrow\,A  \]
          \[                f(\beta):\,=\bar{\beta}_i\:,\]
 wobei $\bar{\beta}_i$ die Einschr"ankung von $\beta$ auf die Menge $S-St_i$
 ist. Wegen (4.12) gilt $(S-St_i)\bar{\beta}_i\subseteq S-S\alpha$, d.h. $\bar{\beta}_i
 \in A$. Nun gibt es zu jedem $a_i$ aus $A$ ein eindeutig bestimmtes maximales
 Element $\beta$ von ${\cal M}^n_{<\alpha}$ mit $\bar{\beta}_i=a_i$. (Auf der
 Menge $St_i$ sei $\beta$ definiert als $t_i^{-1}\alpha$, auf der Menge 
 $S-St_i$ sei $\beta$ die Abbildung $a_i$). $f$ ist somit eine Bijektion.
 Weiterhin gilt:
  \begin{eqnarray*}
     & &\{\beta_0,\ldots,\beta_p\}\in {\cal N}_{\alpha}\\
     &\Leftrightarrow& \mbox{f"ur $\bar{\beta}_{i_0},\ldots,\bar{\beta}_{i_p}$
                             ist $i_{\nu}\not=i_{\mu}$ und B($\bar{\beta}_{i_{\nu}})
                             \cap$B($\bar{\beta}_{i_{\mu}})=\emptyset$}\\
     & & (s.\mbox{ Lemma 4.4)}\\
     &\Leftrightarrow& \{\bar{\beta}_{i_0},\ldots,\bar{\beta}_{i_p}\}\in
                       \Sigma_{\alpha}\:.
 \end{eqnarray*}
 \hfill $\Box$\\
 
 Wir betrachten nun zun"achst folgende allgemeine Situation. Sei $\Gamma$ 
 ein Graph mit Eckenmenge $V$ und Kantenmenge $E\subseteq \{\,\{v,w\}\mid
 v,w\in V\,\}$. Aus $\Gamma$ l"a\3t sich ein Simplizialkomplex ${\cal K}
 (\Gamma)$ konstruieren:$\;$
 Die Eckenmenge von ${\cal K}(\Gamma)$ ist ebenfalls $V$, und f"ur $p\ge 1$ 
 sind die $p$-Simplizes von ${\cal K}(\Gamma)$ diejenigen Teilmengen
 $\{v_0,\ldots,v_p\}\subseteq V$ mit
          \[  \{v_i,v_j\}\in E\quad \forall i\not= j,\;0\le i,j\le p\:. \]
  Das $1$-Ger"ust von ${\cal K}(\Gamma)$ ist gerade $\Gamma$. Besitzt 
 $\Gamma$ gewisse Eigenschaften, so ist ${\cal K}(\Gamma)$ homotopie"aquivalent
 zu einem $(n-1)$-dimensionalen Sph"arenbouquet.
 
\begin{lemma}
  Sei $\Gamma_n$ ein $n$-gef"arbter Graph, $V$ die Eckenmenge und $E$ die
  Kantenmenge von $\Gamma_n$. Erf"ullt f"ur jedes $i\in\{1,\ldots,n\}$ die
  Eckenmenge $V_i$ der Farbe $i$ von $\Gamma_n$ die Bedingungen
   \begin{itemize}
     \item[{\rm 1.}] \#$\,V_i\ge 2$
     \item[{\rm 2.}] Zu je $2(n-1)$ vielen Ecken $v^1,\ldots,v^{2(n-1)}$ aus $V-V_i$
               gibt es $2$ Ecken $v_i,w_i$ aus $V_i$ mit
               \[ \{v_i,v^j\},\{w_i,v^j\}\in E\quad\mbox{f"ur alle}\quad
                                             1\le j\le 2(n-1) \:, \]
    \end{itemize}
  so ist ${\cal K}(\Gamma_n)$ homotopie"aquivalent zu einem $(n-1)$-dimensionalen
  Sph"arenbouquet.
 \end{lemma}

 {\sc Beweis}:$\;$ Der Nachweis des behaupteten Homotopietypes von ${\cal K}
   (\Gamma_n)$ geschieht durch vollst"andige Induktion nach $n$. Wir benutzen
   dazu folgendes Kriterium aus [Qu], Abschnitt 8:\\

    F"ur einen endlichdimensionalen Simplizialkomplex $K$ mit $d=$dim $K$ gilt:
     \begin{eqnarray}
        & & K\simeq\vee_{j\in J} S\,^d \quad\mbox{($S\,^d$ ist die $d$-Sph"are)
                                  mit $J\not=\emptyset$}\\
        &\Leftrightarrow& K\mbox{ ist $(d-1)$-zusammenh"angend und nicht
                                   zusammenziehbar.}\nonumber
     \end{eqnarray}
  (Dabei bedeutet $-1$-zusammenh"angend: nicht leer).\\\\
   \parbox{2cm}{$n=1$:}\hfill
   \parbox[t]{13.3cm}
   { Es ist dim ${\cal K}(\Gamma_1)=0$ und $\#\,V=\,\#\,V_1\ge 2$ (s. 
   Vorraussetzung 1.), d.h es ist ${\cal K}(\Gamma_1)\not=\emptyset$ und nicht
   zusammenziehbar. Aus (4.14) folgt ${\cal K}(\Gamma_1)\simeq\vee S\,^0$.}\\\\
   \parbox{2cm}{$n=2$:}\hfill
   \parbox[t]{13.3cm}{ Wir zeigen zun"achst: ${\cal K}(\Gamma_2)$ ist
      wegzusammenh"angend. Seien $v,w\in V=V_1\cup V_2$. Es gibt 2 F"alle:
       \begin{itemize}
       \item[a)]$\;$ $v,w$ stammen aus demselben $V_i$, $i\in \{1,2\}$
       \item[b)]$\;$ $v\in V_1$ und $w\in V_2$.
       \end{itemize}}
   Im Fall a) existiert wegen Voraussetzung 2. ein $z\in V_j,j\not=i$ mit
   $\{v,z\},\{w,z\}\in E$, d.h. $(v,z,w)$ ist ein Weg von $v$ nach $w$.
   Im Fall b) gibt es wieder wegen 2. zum einen 
  \begin{figure}[ht]
  \unitlength0.5cm
  \begin{center}
  \begin{picture}(20,8)
     \put(0,2){\begin{picture}(5,5)
        \multiput(0,0)(0,5){2}{\multiput(0,0)(5,0){2}{\circle*{0.32}}}
        \multiput(0,0)(5,0){2}{\line(0,1){5}}
        \put(0,5){\line(1,-1){5}}
        \put(0,0){\makebox(-1.5,0){$w$}}
        \put(0,5){\makebox(-1.5,0){$v'$}}
        \put(5.2,0){\makebox(1.5,0){$w'$}}
        \put(5,5){\makebox(1.5,0){$v$}}
           \end{picture}}
     \put(15,2){\begin{picture}(5,5)
        \multiput(0,0)(0,5){2}{\multiput(0,0)(5,0){2}{\circle*{0.32}}}
        \multiput(0,0)(5,0){2}{\line(0,1){5}}
        \put(0,5){\line(1,-1){5}}
        \put(0,0){\line(1,1){5}}
        \put(0,0){\makebox(-1.5,0){$w$}}
        \put(0,5){\makebox(-1.5,0){$v$}}
        \put(5.2,0){\makebox(1.5,0){$w'$}}
        \put(5,5){\makebox(1.5,0){$v'$}}
                \end{picture}}
   \put(2.5,0){\makebox(0,1){${\cal K}(\Gamma_2)$ zusammenh"angend (Fall b))}}
   \put(17.5,0){\makebox(0,1){${\cal K}(\Gamma_2)$ nicht zusammenziehbar}}
   \end{picture}
   \end{center}
   \caption{}
   \end{figure}
    ein $w'\in V_2$ mit  $\{v,w'\}\in E$, und zum anderen 
   ein $v'\in V_1$ mit $\{v',w\},\{v',w'\}\in E$ (s. Abbildung 19, linker Teil),
   d.h. $(v,w',v',w)$ verbindet $v$ und $w$. Da es zu $v,v'\in V_1$ $2$ Elemente
   $w,w'\in V_2$ gibt,die beide zu $v$ und $v'$ benachbart sind (s. Abbildung 19, rechter Teil), 
   ist ${\cal K}(\Gamma_2)$ nicht zusammenziehbar (dim ${\cal K}(\Gamma_2)=1$).
   Aus (4.3) folgt: ${\cal K}(\Gamma_2)\simeq\vee S\,^1.$ \\\\
  $n\ge 3$:$\;$ Wir nehmen an, die Behauptung gilt f"ur $n-1$. Der Komplex
  ${\cal K}(\Gamma_n)$ wird schrittweise aufgebaut. Wir beginnen mit dem
  zusammenziehbaren Unterkomplex
           \[   K_0:\,=st(v_1)\quad\mbox{(Stern von $v_1$)}  \]
  f"ur eine Ecke $v_1\in V_1$. Anschlie\3end betrachten wir nach und nach f"ur
  $i=1,\ldots,n$ den von $K_{i-1}$ und den Ecken $w_i$ aus $V_i$ mit $\{v_1,w_i\}
  \notin E$ erzeugten vollen Unterkomplex $K_i$ von ${\cal K}(\Gamma_n)$:
      \[ K_i:\,=U_{\rm voll}\,(K_{i-1},\mbox{ alle $w_i\in V_i$ mit
                                        $\{v_1,w_i\}\notin E$})\:.  \]
  Wir erhalten also 
       \[ K_0\subseteq K_1\subseteq \ldots\subseteq K_n \:, \]
  wobei f"ur $j=1,\ldots,n$ in $K_j$ alle Ecken aus $V_1\cup\ldots
 \cup V_{j-1}\cup V_j$ enthalten sind, aus $V_{j+1}\cup\ldots\cup V_n$ jedoch
 nur die Ecken $v\in V$ mit $\{v_1,v\}\notin E$. Insbesondere ist $K_n={\cal K}
  (\Gamma_n)$ ($K_n$ enth"alt $V=V_1\cup\ldots\cup V_n$).\\
 Jedes Simplex aus $K_i$ ($i\in \{1,\ldots,n\}$), das nicht in $K_{i-1}$ liegt,
 bezieht {\em genau eine\/} der Ecken $w_i$ aus $V_i$ mit $\{v_1,w_i\}\notin E$ mit
 ein. $K_i$ entsteht somit aus $K_{i-1}$, indem man f"ur jede dieser Ecken $w_i$
 einen Kegel "uber Link$_{K_{i-1}}(w_i)$ (den Link von $w_i$ in $K_{i-1}$)
 hinzuf"ugt. Wir zeigen nun mit Hilfe der Induktionsannahme:\\\\
 {\em Behauptung\/} 1$\;$
   F"ur alle $i\in \{1,\ldots,n\}$ und f"ur jede Ecke $w_i\in V_i$ mit
   $\{v_1,w_i\}\notin E$ ist Link$_{K_{i-1}}(w_i)$ homotopie"aquivalent zu einem
   $(n-2)$-dimensionalen Sph"arenbouquet.\\\\
 {\em Beweis\/}:$\;$Wir betrachten die F"alle $i=1$ und $i\in\{2,\ldots,n\}$.\\
 a)$\: i=1$:$\;$ Sei $w_1\in V_1-\{v_1\}$. Eine Ecke $v\in V$ ist genau dann in 
  Link$_{K_0}(w_1)$, wenn sie sowohl zu $v_1$ als auch zu $w_1$ benachbart ist
  (insbesondere ist $v\notin V_1$). Mit 
      \[  V_j':\,=\{ v\in V_j\mid \{w_1,v\},\{v_1,v\}\in E \}  \]
  f"ur $j=2,\ldots,n$ ist somit $V':\,=V_2'\cup\ldots\cup V_n'$ die Eckenmenge
  von Link$_{K_0}(w_1)$. D.h. Link$_{K_0}(w_1)={\cal K}(\Gamma_{n-1})$, wobei
  $\Gamma_{n-1}$ der $(n-1)$-gef"arbte Untergraph von $\Gamma_n$ mit Eckenmenge
  $V'$ und Kantenmenge $E'=\{\,\{v,w\}\in E\mid v,w\in V'\,\}$ ist. Nach 
  Vorraussetzung  gibt es zu $v_1,w_1$ ($\#\,\{v_1,w_1\}=2\le 2(n-1)$ f"ur alle $n\ge 2$)
  in jedem $V_j,j\not= 1$ mindestens 2 Elemente, die sowohl zu $v_1$ als auch zu $w_1$
  benachbart sind. Also ist $\#\,V_j'\ge 2$. Weiterhin gibt es zu $2(n-2)=2(n-1)-2$
   vielen Elementen $v^1,\ldots,v^{2(n-2)}$ aus $V'-V_j'\subseteq V-V_j$ 
   (ebenfalls nach Vor. 2.) 2 Elemente $v_j,w_j\in V_j$
   mit
   \begin{eqnarray*}
     & &\{v_j,v^k\},\{w_j,v^k\}\in E\quad\mbox{f"ur alle $1\le k\le 2(n-2)$}\\
     & & \{v_j,v_1\},\{v_j,w_1\},\{w_j,v_1\},\{w_j,w_1\}\in E\:.\nonumber 
   \end{eqnarray*}
  Daraus folgt $v_j,w_j\in V_j'$ und $\{v_j,v^k\},\{w_j,v^k\}
  \in E'$ f"ur $1\le k\le 2(n-2)$. Somit ist $\Gamma_{n-1}$ ein $(n-1)$-gef"arbter
  Graph, der die Voraussetzungen aus Lemma 4.7 erf"ullt. Nach Induktionsannahme
  folgt: Link$_{K_0}(w_1)={\cal K}(\Gamma_{n-1})\simeq\vee S\,^{(n-2)}$.\\

 b)$\:i\in \{2,\ldots,n\}$:$\;$ $K_{i-1}$ enth"alt alle Ecken aus $V_1\cup
  \ldots\cup V_{i-1}$, jedoch f"ur $i\le j\le n$ nur diejenigen $v_j\in V_j$, 
  die zu st$(v_1)$ geh"oren. Die Eckenmenge von Link$_{K_{i-1}}(w_i)$ ist somit
      \[ V'=V_1'\cup\ldots\cup V_{i-1}\!'\cup V_{i+1}\!'\cup\ldots\cup V_n'\,, 
                 \quad\mbox{wobei}  \]
   \begin{eqnarray*}
     V_j'&:\,=&\{ v_j\in V_j\mid \{w_i,v_j\}\in E\}\quad
                           \mbox{f"ur $1\le j\le i-1$}  \\
     V_j'&:\,=&\{ v_j\in V_j\mid \{v_1,v_j\},\{w_i,v_j\}\in E\}\quad
                           \mbox{f"ur $i+1\le j\le n$}\:.
   \end{eqnarray*}
  Sei wieder $\Gamma_{n-1}$ der von der Eckenmenge $V'$ erzeugte Untergraph von $\Gamma_n$.
  $\Gamma_{n-1}$ ist ein $(n-1)$-gef"arbter Graph, und analog zu Teil a)
  "uberlegt man sich, da\3 er die Voraussetzungen von Lemma 4.7 erf"ullt.
  (F"ur $i+1\le j\le n$ ist die Situation v"ollig analog zu Teil a), f"ur
  $1\le j\le i-1$ vereinfacht sie sich sogar). Somit ist auch hier Link$_{K_{i-1}}(w_i)
  \simeq\vee S\,^{(n-2)}$.\\
 
 Wir wollen nun die Homologiegruppen von ${\cal K}(\Gamma_n)=K_n$ bestimmen
 und zeigen dazu:\\[0.46cm]
 {\em Behauptung\/} 2$\;$ Sei $X$ ein Simplizialkomplex, $K$ ein Unterkomplex und
 $w$ eine Ecke aus $X-K$. Sind die Homologiegruppen $H_i\,(K)=0$ f"ur $1\le i<n-1$
 und ist Link$_K(w)$ homotopie"aquivalent zu einem $(n-2)$-dimensionalen
 Sph"arenbouquet, so gilt auch
    \[ H_i(\,K\cup C(w)\,)=0\quad 1\le i<n-1\,,  \]
 wobei $C(w)$ der Kegel "uber Link$_K(w)$ ist.\\

 {\em Beweis\/}:$\;$ Es ist $K\cap C(w)=$Link$_K(w)$. Aus dem Homotopietyp von Link$_K(w)$
 folgt $H_j(\,K\cap C(w)\,)=0$ f"ur $1\le j<n-2$. Betrachten wir die
 Mayer-Vietoris-Sequenz f"ur $1\le i<n-1$
         \[ \ldots\rightarrow\,\underbrace{H_i\,(K)}_{=0}
      \oplus\underbrace{H_i\,(C(w))}_{=0}
      \,\rightarrow\, H_i(\,K\cup C(w)\,)\rightarrow
      \underbrace{H_{i-1}(\,K\cap C(w)\,)}_{=0}
      \,\rightarrow\,\ldots\:,   \]
 so folgt die Behauptung.\\
  
 Gilt nun f"ur ein $j\in \{0,\ldots,n-1\}$ $H_i(K_j)=0$ f"ur $1\le i<n-1$,
 so l"a\3t sich mit Hilfe von Behauptung 2 zeigen, da\3 f"ur diese $i$
 auch $H_i(K_{j+1})=0$ ist:\\
 $K_{j+1}$ entsteht aus $K_j$ durch Hinzuf"ugen eines Kegels $C(w)$ "uber
 Link$_{K_j}(w)$ f"ur jede Ecke $w\in W:\,=\{v\in V_{j+1}\mid \{v_1,v\}\notin E\}$.
 Wir betrachten nun f"ur eine endliche Teilmenge $T\subseteq W$ den Unterkomplex
    \[ U_T=K_j\cup\,\bigcup\limits _{w\in T}C(w).  \]
Da nach Behauptung 1 Link$_{K_j}(w)\simeq\vee S^{n-2}$ ist, folgt durch
 sukzessives Hinzuf"ugen der $C(w)$ f"ur $w\in T$ mit Hilfe von
 Behauptung 2
       \[H_i(U_T)=0\quad \mbox{f"ur alle $1\le i<n-1$}. \]
 Nun ist
       \[   K_{j+1}=\lim_{\to T}\,U_T\:,  \]
 wobei $T$ die Menge aller endlichen Teilmengen von $W$ durchl"auft, d.h.
 $H_i(K_{j+1})=0$ f"ur alle $1\le i<n-1$.\\
 Da $K_0=st(v_1)$ zusammenziehbar ist, also $H_i(K_0)=0$ f"ur alle $i\ge 1$,
 kann man diesen Proze\3 beginnend bei $K_0$ f"ur $j=0,\ldots,n-1$ iterieren
 und erh"alt
   \vfill
   \begin{equation}
        H_i({\cal K}(\Gamma_n))=H_i(K_n)=0\quad\mbox{f"ur alle $1\le i<n-1$}\,.
   \end{equation}
  \vfill
 Um Informationen "uber die $(n-1)$-dimensionalen Homologiegruppen von ${\cal K}
 (\Gamma_n)$ zu gewinnen, betrachten wir zu irgendeiner Ecke $w_1\in V_1-
 \{v_1\}$ die Mayer-Vietoris-Sequenz bez"uglich des Unterkomplexes $K_0\cup C$
 mit $C=C(w_1,{\rm Link}_{K_0}(w_1))$:
   \vfill
  \[ \ldots\,\rightarrow\ H_{n-1}(K_0\cup C)\,\rightarrow\,
     H_{n-2}(\underbrace{K_0\cap C}_{{\rm Link}_{K_0}(w_1)})
     \,\rightarrow\,\underbrace{ H_{n-2}(K_0)\oplus H_{n-2}(C)}_{=0}
     \,\rightarrow\,\ldots\:. \]
    \vfill
 Nach Behauptung 1 ist $H_{n-2}({\rm Link}_{K_0}(w_1))=\bigoplus\,{\sl Z}\!\!
 {\sl Z}$, d.h. $H_{n-1}(K_0\cup C)\not=0$. Wegen dim ${\cal K}(\Gamma_n)=
 n-1$ ist damit auch
    \vfill
     \begin{equation}
         H_{n-1}({\cal K}(\Gamma_n))\not=0.
     \end{equation}
    \vfill
 Ist ${\cal K}(\Gamma_n)$ f"ur $n\ge 3$ einfachzusammenh"angend, so folgt
 aus (4.15) (vgl. zum Beispiel [Sp], Hurewicz-Theorems):
     \vfill
    \begin{eqnarray*}
      & & {\cal K}(\Gamma_n) \mbox{ ist $(n-2)$-zusammenh"angend und}\\
      & & \hspace*{1cm} 
          \pi_{n-1}({\cal K}(\Gamma_n))\cong H_{n-1}({\cal K}(\Gamma_n)) 
        \quad\mbox{f"ur}\quad 1\le i\le n-1 \:.
   \end{eqnarray*}
     \vfill
  Wegen (4.16) ist $\pi_{n-1}({\cal K}(\Gamma_n))\not=0$. Also ist ${\cal K}
 (\Gamma_n)$ ist $(n-2)$-zusammenh"angend und nicht zusammenziehbar. Aus (4.14)
 folgt wieder:$\:{\cal K}(\Gamma_n)\simeq\vee S\,^{n-1}$. Es bleibt zu zeigen:\\\\
 {\em Behauptung\/} 3:$\;$ F"ur $n\ge 3$ ist ${\cal K}(\Gamma_n)$ einfach zusammenh"angend.\\\\
 {\em Beweis\/}:$\;$ Wir betrachten zuerst den Fall
 $n=3$:$\;$ Ein geschlossener
  \begin{figure}[ht]
  \unitlength 0.5cm
  \begin{center}
  \begin{picture}(4.5,3)
    \put(0,0.25){\begin{picture}(5,2.5)
              \put(0,0){\line(1,0){5}}
              \put(0,0){\line(1,1){2.5}}
              \put(2.5,2.5){\line(1,-1){2.5}}
              \multiput(0,0)(5,0){2}{\circle*{0.3}}
              \put(2.5,2.5){\circle*{0.3}}
              \put(0,0){\makebox(-1.3,0){$v_0$}}
              \put(2.5,2.5){\makebox(0,1.6){$v_1$}}
              \put(5,0){\makebox(1.3,0){$v_2$}}
              \end{picture}}
  \end{picture}\end{center}
  \end{figure}
  Kantenweg $(v_0,v_1,v_2,v_0)$, der 3 verschiedene Eckpunkte
 $v_0,v_1,v_2\in V$ durchl"auft, ist Rand eines $2$-Simplexes aus ${\cal K}(\Gamma_n)$ 
 und somit nullhomotop. F"ur einen geschlossenen Kantenweg $(v_0,v_1,\ldots,
  v_m,v_0)$ mit $m+1\ge 4$ verschiedenen Eckpunkten aus $V$ zeigen wir, da\3 dieser homotop
  zu einem Kantenweg $(v_0',v_1',\ldots,v'\!_{m-1},v_0')$ mit $m$ Eckpunkten
  ist:\\
   \vfill
 Seien $v_0,v_1,v_2,v_3$ die ersten 4 Eckpunkte von $(v_0,v_1,\ldots,v_m,v_0)$.
 Geh"ort $v_0$ zur Farbe $F_1$, $\,v_1$ zur Farbe $F_2$, so gibt es f"ur $v_2$ 
 und $v_3$ nur die 4 M"oglichkeiten
     \pagebreak
      \begin{itemize}
      \item[1.] $v_2$ geh"ort zur Farbe $F_1$, $v_3$ zur Farbe $F_2$
      \item[2.] $v_2$ geh"ort zur Farbe $F_1$, $v_3$ zur Farbe $F_3$
      \item[3.] $v_2$ geh"ort zur Farbe $F_3$, $v_3$ zur Farbe $F_2$
      \item[4.] $v_2$ geh"ort zur Farbe $F_3$, $v_3$ zur Farbe $F_1$.
      \end{itemize}
  \vspace{0.2cm}
 Zu 1.:$\;$ Zu den 4 Elementen $v_0,v_1,v_2,v_3$ au\3erhalb der Farbe $F_3$
 \begin{figure}[ht]
  \unitlength0.5cm
  \begin{center}
  \begin{picture}(20,7)
  \put(0,1.5){\begin{picture}(5,5)
              \linethickness{0.65mm}
              \multiput(0,0.2)(5,0){2}{\line(0,1){4.6}}
              \put(0.2,5){\line(1,0){4.6}}\thinlines
              \multiput(0,5)(5,-5){2}{\circle*{0.4}}
              \multiput(0,0)(5,5){2}{\circle{0.4}}
              \put(2.5,2.5){\makebox(0,0){.}}
              \put(2.5,2.5){\circle{0.4}}
              \put(0.14,0.14){\line(1,1){2.2}}
              \put(2.64,2.64){\line(1,1){2.2}}
              \put(0.14,4.86){\line(1,-1){2.2}}
              \put(2.64,2.36){\line(1,-1){2.2}}
              \put(0,0){\makebox(-1.3,0){$v_0$}}
              \put(0,5){\makebox(-1.3,0){$v_1$}}
              \put(5,5){\makebox(1.3,0){$v_2$}}
              \put(5,0){\makebox(1.3,0){$v_3$}}
              \end{picture}}
   \put(8,4){\vector(1,0){4}}
   \put(15,1.5){\begin{picture}(5,3)
              \multiput(0,0)(5,0){2}{\circle{0.4}}
              \put(2.5,2.5){\circle{0.4}}
              \put(2.55,2.5){\makebox(0,0){.}}
              \thicklines
              \put(0.14,0.14){\line(1,1){2.2}}
              \put(2.64,2.36){\line(1,-1){2.2}}
              \put(0,0){\makebox(-1.3,0){$v_0$}}
              \put(2.5,2.5){\makebox(0,1.3){$v$}}
              \put(5,0){\makebox(1.3,0){$v_3$}}
              \end{picture}}
   \end{picture}
   \begin{minipage}{15.3cm}\unitlength0.5cm\begin{center}
       \begin{picture}(0.3,0)\put(0,0.25){\circle{0.4}}\end{picture}
        = Farbe $F_1$,
       \begin{picture}(0.8,0)\put(0.5,0.25){\circle*{0.4}}\end{picture}
       = Farbe $F_2$,
       \begin{picture}(0.8,0)\put(0.5,0.25){\circle{0.4}}
               \put(0.5,0.25){\makebox(0,0){.}}\end{picture}
       = Farbe $F_3$\end{center}
  \end{minipage}
  \end{center}\end{figure}
  gibt es ein $v\in V$  \vspace*{-0.2cm}
  der Farbe $F_3$, das zu jedem dieser 4 Elemente benachbart
  ist (Voraussetzung 2.). Da somit $\{v_0,v_1,v\},\{v_1,v_2,v\}$ und $\{v_2,v_3,v\}$
 jeweils $2$-Simplizes aus ${\cal K}(\Gamma_3)$ sind, ist der Kantenzug
 $(v_0,v_1,v_2,v_3)$ homotop zu $(v_0,v,v_3)$.\\
 In den restlichen 3 F"allen kann man mit "ahnlichen "Uberlegungen zeigen, da\3
 $(v_0,v_1,v_2,v_3)$ homotop zu einem Kantenweg der Form $(v_0,v,v_3)$ ist.
 Somit ist jeder geschlossene Kantenweg aus ${\cal K}(\Gamma_3)$ nullhomotop
 ( iteriere das Verfahren bis nur noch 3 Eckpunkte durchlaufen werden), d.h.
 $\pi_1({\cal K}(\Gamma_3))=0$.\\
 $n\ge 4$:$\;$ Man "uberlegt sich, da\3 jeder geschlossene Kantenweg $(v_0,v_1,\ldots,v_m,v_0)$aus
 ${\cal K}(\Gamma_n)$ homotop ist zu einem Weg $(v_0',v_1',\ldots,v_{m'}',v_0')$, dessen
 Eckpunkte h"ochstens 3 verschiedene Farben besitzen. Damit ist der Fall zur"uckgef"uhrt
 auf $n=3$. \hfill $\Box$\\\\
Mit Hilfe des letzten Lemmas l"a\3t sich der Homotopietyp von $\Sigma_{\alpha}$
 bestimmen, denn ist $\Gamma_{\alpha}$ das $1$-Ger"ust von $\Sigma_{\alpha}$, so
 gilt $\Sigma_{\alpha}={\cal K}(\Gamma_{\alpha})$.
 \begin{lemma}
    $\Gamma_{\alpha}$ erf"ullt die Voraussetzungen von Lemma 4.7, d.h.
    ${\cal K}(\Gamma_{\alpha})=\Sigma_{\alpha}$ ist homotopie"aquivalent zu einem $(n-1)$-dimensionalen
   Sph"arenbouquet.
  \end{lemma}
 {\sc Beweis}:$\;$ Ordnet man der Eckenmenge $A_i$ die Farbe $i$ zu, so ist 
 $\Gamma_{\alpha}$ ein $n$-gef"arbter Graph. Da schon jeweils {\em eine \/}
 Menge der Form $((x_1,i_1),\ge m_1)$ und $((y_1,j_1),\ge l_1)$ aus $S-S\alpha$
 (und es gibt jeweils $gr(\alpha)\ge 2n$ viele) die Konstruktion unendlich vieler
  $\varphi: (S-St_i)\longrightarrow (S-S\alpha)$ erm"oglicht, die zu einer 
 Abbildung von $S$ nach $S$ aus $\cal M$
  \begin{figure}[ht]
  \unitlength 0.6cm
  \begin{center}
  \begin{picture}(23,13.5)
   \multiput(0,2.5)(13,0){2}{\begin{picture}(10,10)
        \multiput(1,0)(1,0){9}{\circle*{0.2}}
        \multiput(0,1)(0,1){9}{\circle*{0.2}}
        \thicklines
        \put(0,0){\vector(1,0){10}}
        \put(0,0){\vector(0,1){10}}
        \thinlines
        \put(1,0){\makebox(0,-1.3){1}}
        \put(0,1){\makebox(-1.3,0){1}}
        \put(10,0){\makebox(1.3,0){x}}
        \put(0,10){\makebox(0,1.4){y}}
        \end{picture}}
   \put(0,2.5){\begin{picture}(10,10)
             \thicklines
              \put(1,1){\line(1,0){9}}\put(1,0.75){\line(0,1){0.5}}
              \put(1,2){\line(0,1){8}}\put(0.75,2){\line(1,0){0.5}}
              \thinlines
              \put(1,2.25){\line(1,0){11.35}}
              \put(1,2){\line(4,1){11.6}}
              \end{picture}}
    \put(13,2.5){\begin{picture}(10,10)
               \put(0.4,2){\vector(1,0){3.25}}
               \put(0.4,5){\vector(4,1){3.25}}
               \put(4.5,2){\makebox(3,0){$\scriptstyle ((x_1,i_1),\ge m_1)$}}
               \put(4.5,6){\makebox(3,0){$\scriptstyle ((x_1,i_1),\ge m')$}}
               \thicklines
               \put(4,2){\line(0,1){8}}\put(3.75,2){\line(1,0){0.5}}
               \put(3.8,6){\line(0,1){4}}\put(3.55,6){\line(1,0){0.5}}
               \thinlines
               \put(1.5,2){\makebox(0,1.5){$\varphi$}}
               \put(1.5,5.3){\makebox(0,1.5){$\varphi$'}}
               \multiput(2.5,7)(0,1){3}{\circle*{0.2}}
               \end{picture}}
   \put(0,0){$S-St_i\subseteq$ Quadrant $Q_i$}
    \put(13,0){Quadrant $Q_{i_1}$}
    \end{picture}\end{center}
   \parbox{15.3cm}{\small Schon f"ur die Menge $\scriptstyle
      \{((1,y),i)\mid y\in {\rm I}\!{\rm N}\}
   \subseteq S-St_i$ gibt es f"ur die Wahl von $\scriptstyle \varphi $ beliebig
   viele M"oglichkeiten, indem man die Verschiebung in Richtung der y-Achse 
   variiert.}
   \caption{ Anzahl der Elemente in $A_i$}
   \end{figure}
      
 fortgesetzt werden k"onnen (s. Abbildung 20), ist
 f"ur alle $i=1,\ldots,n$ \[  \#\,A_i=\infty\ge 2\:. \]
 Seien nun $a^1,\ldots,a^{2(n-1)}$ $2(n-1)$ viele Elemente aus $A-A_i$. Wir
  m"ussen zeigen, da\3
  es in $A_i$ mindestens zwei Elemente $a_i$ und $b_i$  gibt mit 
    \[  {\rm B}(a_i)\cap{\rm B}(a^k)={\rm B}(b_i)\cap{\rm B}(a^k)=\emptyset\]
  f"ur alle $k=1,\ldots,2(n-1)$.
  Dazu weisen wir nach, da\3 in der Menge $(S-S\alpha)-\left({\rm B}(a^1)\cup
  \ldots\cup{\rm B}(a^{2(n-1)})\right)$ noch gen"ugend "`Platz"' vorhanden ist, um 
  sogar beliebig viele Abbildungen aus $A_i$ zu finden, deren Bildmenge
 in $(S-S\alpha)-\left({\rm B}(a^1)\cup
  \ldots\cup{\rm B}(a^{2(n-1)})\right)$ enthalten ist.
 
 Wir wissen nach Lemma 3.4:
    \[ S-S\alpha=\left(\,\bigcup_{\nu=1}^{gr(\alpha)}((x_{\nu},i_{\nu}),\ge m_{\nu})
        \cup((y_{\mu},j_{\mu}),\ge l_{\mu})\,\right)\,\cup P \:.\]
 Die Bildmenge eines $a^k$ "`verbraucht"' von $S-S\alpha$ h"ochstens eine der Mengen
 $((x_{\nu_k},i_{\nu_k}),\ge m_{\nu_k})$, eine der Mengen $((y_{\mu_k},j_{\mu_k}),\ge l_{\mu_k})$
 sowie eine endliche Menge $E_k$. Somit ist
    \begin{eqnarray}
      & &{\rm B}(a^1)\cup\ldots\cup{\rm B}(a^{2(n-1)})\\
      &\subseteq& ((x_{\nu_1},i_{\nu_1}),\ge m_{\nu_1})\cup\ldots\cup((x_{\nu_{2(n-1)}},
       i_{\nu_{2(n-1)}}),\ge m_{\nu_{2(n-1)}})\cup\nonumber\\
    & & ((y_{\mu_1},j_{\mu_1}),\ge l_{\mu_1})\cup\ldots\cup
         ((y_{\mu_{2(n-1)}},j_{\mu_{2(n-1)}}),\ge l_{\mu_{2(n-1)}})
        \cup\nonumber\\
   & & E=E_1\cup\ldots\cup E_{2(n-1)}\nonumber\:.
   \end{eqnarray}
 Im schlechtesten Fall sind dabei die $\:(x_{\nu_1},i_{\nu_1}),\ldots,
(x_{\nu_{2(n-1)}},i_{\nu_{2(n-1)}})\:$,
 sowie die $\:(y_{\mu_1},j_{\mu_1}),
 \ldots,(y_{\mu_{2(n-1)}},j_{\mu_{2(n-1)}})$
 alle verschieden. Da nach Voraussetzung 
 $gr(\alpha)\ge 2n$, also $gr(\alpha)-2(n-1)\ge 2$ ist, gibt
 es ein $\nu\in\{1,\ldots,gr(\alpha)\}-\{\nu_1,\ldots,\nu_{2(n-1)}\}$ und ein
  $\mu\in\{1,\ldots,gr(\alpha)\}-\{\mu_1,\ldots,\mu_{2(n-1)}\}$.\\
  F"ur ein solches $\nu$ bzw. $\mu$ ist $(x_{\nu},i_{\nu})\not=(x_{\nu_k},i_{\nu_k})$ und
  $(y_{\mu},j_{\mu})\not=(y_{\mu_k},j_{\mu_k})$ f"ur $k=1,\ldots,2(n-1)$.
 Da Mengen der
  Form $((x,i),\ge m)$ und $((x,i)',\ge m')$ im Falle $(x,i)\not=(x,i)'$ zueinander
 disjunkt sind (dasselbe gilt f"ur die Mengen der Form $((y,j),\ge l)$), folgt
 daher aus (4.17)
   \[ \left(\,((x_{\nu},i_{\nu}),\ge m_{\nu})\cup((y_{\mu},j_{\mu}),\ge l_{\mu})\,\right)
     \bigcap\left(\,{\rm B}(a^1)\cup\ldots\cup{\rm B}(a^{2(n-1)})\,\right)
        \subseteq E \:,\]
 wobei $E$ eine endliche Menge ist. Wir w"ahlen nun $m_{\rm neu}$,
  $l_{\rm neu}\in {\rm I}\!{\rm N}$, soda\3 
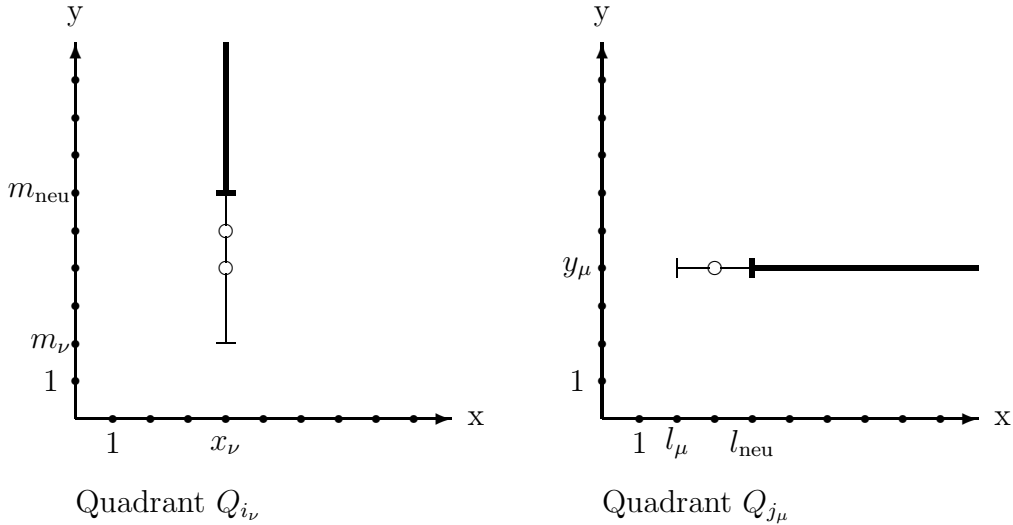
  \begin{figure}[ht]
  \unitlength0.5cm
  \begin{center}
  \begin{picture}(24,13)
   \multiput(0,2.5)(14,0){2}{\begin{picture}(10,10)
        \multiput(1,0)(1,0){9}{\circle*{0.2}}
        \multiput(0,1)(0,1){9}{\circle*{0.2}}
        \thicklines
        \put(0,0){\vector(1,0){10}}
        \put(0,0){\vector(0,1){10}}
        \thinlines
        \put(1,0){\makebox(0,-1.3){1}}
        \put(0,1){\makebox(-1.3,0){1}}
        \put(10,0){\makebox(1.3,0){x}}
        \put(0,10){\makebox(0,1.4){y}}
        \end{picture}}
   \put(0,2.5){\begin{picture}(10,10)
             \linethickness{0.65mm}
              \put(4,6){\line(0,1){4}}\put(3.75,6){\line(1,0){0.5}}
              \thinlines
              \put(4,2){\line(0,1){1.85}}\put(3.75,2){\line(1,0){0.5}}
              \put(4,4.15){\line(0,1){0.7}}
              \put(4,5.15){\line(0,1){0.85}}
              \multiput(4,4)(0,1){2}{\circle{0.32}}
              \put(0,6){\makebox(-1.9,0){$m_{\rm neu}$}}
              \put(0,2){\makebox(-1.3,0){$m_{\nu}$}}
              \put(4,0){\makebox(0,-1.4){$x_{\nu}$}}
              \end{picture}}
    \put(14,2.5){\begin{picture}(10,10)
               \put(2,4){\line(1,0){0.85}}\put(2,3.75){\line(0,1){0.5}}
               \put(3.15,4){\line(1,0){0.85}}
               \linethickness{0.65mm}
               \put(4,4){\line(1,0){6}}
               \put(4,3.75){\line(0,1){0.5}}
               \thinlines
               \put(3,4){\circle{0.32}}
               \put(0,4){\makebox(-1.3,0){$y_{\mu}$}}
               \put(2,0){\makebox(0,-1.3){$l_{\mu}$}}
               \put(4,0){\makebox(0,-1.3){$l_{\rm neu}$}}
               \end{picture}}
   \put(0,0){Quadrant $Q_{i_{\nu}}$}
    \put(14,0){Quadrant $Q_{j_{\mu}}$} 
   \end{picture}\end{center}
   \parbox{15.3cm}{\small Die Kreise stellen die endliche Menge  
    $\scriptstyle (\:((x_{\nu},i_{\nu}),\ge m_{\nu})\cup((y_{\mu},j_{\mu}),\ge l_{\mu})\:)
    \cap (\:{\rm B}(a^1)\cup\ldots\cup{\rm B}(a^{2(n-1)})\:)$ dar. Die
    fettgedruckten Linien entsprechen den Mengen $\scriptstyle ((x_{\nu},i_{\nu}),
    \ge m_{\rm neu})$ bzw $\scriptstyle ((y_{\mu},j_{\mu}),\ge l_{\rm neu})$.}
    \caption{ Wahl von $m_{\rm neu}$ und $l_{\rm neu}$}
   \end{figure}

        \[ (\,((x_{\nu},i_{\nu}),\ge m_{\rm neu})\cup((y_{\mu},j_{\mu}),
           \ge l_{\rm neu})\,)\cap(\,{\rm B}(a^1)\cup\ldots\cup{\rm B}(a^{2(n-1)})\,)
           =\emptyset   \]
   ist (s. Abbildung 21). Wie bei dem Nachweis von $\#\,A_i=\infty$ zu Beginn
   des Beweises "uberlegt man sich, da\3 es 
   beliebig viele $a_i\in A_i$ gibt, deren
   Bildmenge B($a_i)$ in $((x_{\nu},i_{\nu}),\ge m_{\rm neu})\cup((y_{\mu},j_{\mu}),
   \ge l_{\rm neu})$ enthalten ist
     und daher leeren Schnitt mit ${\rm B}(a^1)\cup\ldots
   \cup{\rm B}(a^{2(n-1)})$ besitzt. \hfill $\Box$\\

 Wir haben somit gezeigt:

 \begin{satz}
    Die Gruppe $G_n$ ist vom Typ {\sc (FP)}$_{n-1}$ und nicht vom Typ
    {\sc (FP)}$_n$. F"ur $n\ge 3$ ist $G_n$ endlich pr"asentiert.
  \end{satz}
 Weiterhin ergibt sich aus diesem Satz (wie in Abschnitt 2.2 erl"autert):
 \begin{cor}
    Die Gruppe $\widetilde{G}_n$ ist vom Typ {\sc (FP)}$_{n-1}$ und f"ur 
    $n\ge 3$ endlich pr"asentiert.
  \end{cor}

\end{document}